\documentclass{gtart_h}  

\def\ifplaintex{\expandafter\ifx\csname documentclass\endcsname\relax}

\ifplaintex 
\else
\fi

\expandafter\ifx\csname beginpicture\endcsname\relax
\expandafter\ifx\csname documentclass\endcsname\relax
\input pictex \else
\input prepictex \input pictex \input postpictex \fi\fi

\def\gt{{\mathsurround=0pt\it $\cal G\mskip-2mu$eometry \&\ 
$\cal T\!\!$opology}}        

\def\gtp{{\mathsurround=0pt\it $\cal G\mskip-2mu$eometry \&\ 
$\cal T\!\!$opology $\cal P\!$ublications}}  

\def\lognumber#1{\def\thelognumber{#1}}
\def\volumenumber#1{\def\thevolumenumber{#1}}
\def\papernumber#1{\def\thepapernumber{#1}}
\def\volumeyear#1{\def\thevolumeyear{#1}}

\def\pagenumbers#1#2{\def\startpage{#1}\def\finishpage{#2}}
\def\published#1{\def\publishdate{#1}}
\def\proposed#1{\def\theproposer{#1}}
\def\seconded#1{\def\theseconders{#1}}
\def\received#1{\def\receiveddate{#1}}

\def\accepted#1{\def\accepteddate{#1}}

\def\asciiaddress#1{\def\theasciiaddress{#1}}
\def\asciiemail#1{\def\theasciiemail{#1}}

\long\def\asciiabstract#1{\long\def\theasciiabstract{#1}}

\def\shortauthors#1{\def\theshortauthors{#1}}

\let\\\par\let\thelognumber\relax
\let\thevolumenumber\relax\let\thepapernumber\relax
\let\thevolumeyear\relax\let\thesamplenumber\relax\let\startpage\relax
\let\finishpage\relax\let\publishdate\relax\let\receiveddate\relax
\let\reviseddate\relax\let\accepteddate\relax\let\theasciititle\relax
\let\theasciiauthors\relax\let\theasciiaddress\relax
\let\theasciiabstract\relax
\let\theasciiemail\relax\let\theshortauthors\relax\let\theshorttitle\relax

\long\def\maketitlep{   

\count0=\startpage

\gt\hfill      
\beginpicture
\setcoordinatesystem units <0.33truein, 0.33truein> point at 2.2 0.9
\setplotsymbol ({$\cal G$})
\plotsymbolspacing=9truept
\circulararc 315 degrees from 0 1 center at 0 0
\setplotsymbol ({$\cal T$})
\circulararc 315 degrees from 1 -1 center at 1 0
\endpicture
\break
{\small\ifx\thesamplenumber\relax 
Volume \else Sample
\fi\thevolumenumber\ (\thevolumeyear)
\startpage--\finishpage\nl
Published: \publishdate}
\vglue 0.5truein plus 0.4fil minus 0.1truein

{\parskip=0pt\leftskip 0pt plus 1fil\def\\{\par\smallskip}{\ifplaintex\large
\else\Large\fi\bf\thetitle}\par\medskip}   

\vglue 0pt plus 0.1fil 

{\parskip=0pt\leftskip 0pt plus 1fil\def\\{\par}{\sc\theauthors}
\par\medskip}

\vglue 0pt plus 0.1fil 

{\small\parskip=0pt\let\newline\\
{\leftskip 0pt plus 1fil\def\\{\par}{\sl\theaddress}\par}
\expandafter\ifx\theemail\relax    
\relax\else\vglue 5pt plus 0.02fil minus 2pt\def\\{\stdspace{\rm 
and}\stdspace} 
\cl{Email:\stdspace\tt\theemail}\fi
\ifx\theurl\relax                  
\relax\else\vglue 5pt plus 0.02fil minus 2pt\def\\{\stdspace{\rm 
and}\stdspace}
\cl{URL:\stdspace\tt\theurl}\fi\par}

\vglue 7pt plus 0.3fil minus 3pt

{\bf Abstract}
\vglue 5pt plus 0.1fil minus 2pt

\theabstract

\vglue 7pt plus 0.3fil minus 3pt

{\bf AMS Classification numbers}\quad Primary:\quad \theprimaryclass

Secondary:\quad \thesecondaryclass

\vglue 5pt plus 0.3fil minus 2pt

{\bf Keywords:}\quad \thekeywords

\vglue 10pt plus 0.5fil minus 5pt

{\small  Proposed: \theproposer\hfill Received: \receiveddate\nl
Seconded: \theseconders\hfill 
\ifx\reviseddate\relax                         
Accepted: \accepteddate                        
\else
Revised: \reviseddate                          
\fi}
\eject
}       

\let\maketitlepage\maketitlep
\let\maketitle\maketitlepage

\font\phead=cmsl9 scaled 950
\font=cmsl9 scaled 1050
\font\pnum=cmbx10 scaled 913
\font\lnum=cmbx10 
\font\pfoot=cmsl9 scaled 950
\font=cmsl9 scaled 1050
\ifplaintex
\headline{\vbox to 0pt{\vskip -4.5mm\line{\small\phead\ifnum
\count0=\startpage ISSN 1364-0380 (on line)
1465-3060 (printed) \hfill {\pnum\folio}\else\ifodd\count0\def\\{ }%
\ifx\theshorttitle\relax\thetitle\else\theshorttitle\fi\hfill{\pnum\folio}
\else\def\\{ and }{\pnum\folio}\hfill\ifx\theshortauthors\relax\theauthors
\else\theshortauthors\fi\fi\fi}\vss}}
\footline{\vbox to 0pt{\vglue 0mm\line{\small\pfoot\ifnum\count0=\startpage
\copyright\ \gtp\hfill\else
\gt, Volume \thevolumenumber\ (\thevolumeyear)\hfill\fi}\vss
}}
\else
\makeatletter
\def\@oddhead{{\small\lhead\ifnum\count0=\startpage ISSN 1364-0380 (on line)
1465-3060 (printed) \hfill {\lnum\number\count0}\else\ifodd\count0
\def\\{ }\ifx\theshorttitle\relax \thetitle \else\theshorttitle\fi\hfill
{\lnum\number\count0}\else\def\\{ and }{\lnum\number\count0}
\hfill\ifx\theshortauthors\relax 
\theauthors\else\theshortauthors\fi\fi\fi}}\def\@evenhead{\@oddhead}
\def\@oddfoot{\small\lfoot\ifnum\count0=\startpage\copyright\ \gtp\hfill\else
\gt, Volume \thevolumenumber\ (\thevolumeyear)\hfill\fi}
\def\@evenfoot{\@oddfoot}
\makeatother
\fi

\newwrite\gtoutfile
\long\gdef\makeheadfile{  
{\def\\{, }\def\s{ }
\immediate\openout\gtoutfile head.xxx
\immediate\write\gtoutfile{Proxy-for: \ifx\theasciiauthors\relax
\theauthors\else\theasciiauthors\fi\s<\ifx\theasciiemail\relax\theemail\else\theasciiemail\fi>}
\immediate\write\gtoutfile{\noexpand\\}
\immediate\write\gtoutfile{Authors: \ifx\theasciiauthors\relax
\theauthors\else\theasciiauthors\fi}
{\def\\{ }\immediate\write\gtoutfile{Title: \ifx\theasciititle\relax
\thetitle\else\theasciititle\fi}}
\immediate\write\gtoutfile{Subj-class: GT or SG or MG etc}
\immediate\write\gtoutfile{MSC-class: \theprimaryclass\ifx\thesecondaryclass\relax\else, \thesecondaryclass\fi}
\immediate\write\gtoutfile{Journal-ref: Geom. Topol. \thevolumenumber
(\thevolumeyear) \startpage-\finishpage}
\immediate\write\gtoutfile{Comments: Published by Geometry and Topology at}
\immediate\write\gtoutfile{\s\s http://www.maths.warwick.ac.uk/gt/GTVol\thevolumenumber/paper\thepapernumber.abs.html}
\immediate\write\gtoutfile{\noexpand\\}
\immediate\write\gtoutfile{}
\ifx\theasciiabstract\relax
\immediate\write\gtoutfile{\theabstract}\else
\immediate\write\gtoutfile{\theasciiabstract}\fi
\immediate\write\gtoutfile{}
\immediate\write\gtoutfile{\noexpand\\}
\immediate\write\gtoutfile{}
\immediate\closeout\gtoutfile}}  

\def\maketitlepage{\maketitlep\makeheadfile}
\let\maketitle\maketitlepage

\lognumber{511}
\received{1 November 2004}
\volumenumber{9}\papernumber{24}\volumeyear{2005}
\pagenumbers{1043}{1114}   
\published{1 June 2005}
\accepted{30 May 2005}
\proposed{Robion Kirby}
\seconded{Dieter Kotschick, Ronald Stern}

\usepackage{amsfonts, pstricks}

\def\S{Section }
\newtheorem{thm}{Theorem}
\newtheorem{hyp}{Hypothesis}
\newtheorem{prop}{Proposition}
\newtheorem{defn}{Definition}
\newtheorem{lem}{Lemma}
\newcommand{\bR}{{\mathbb{R}}}
\newcommand{\bC}{{\mathbb{C}}}
\newcommand{\bCP}{{\mathbb{CP}}}
\newcommand{\bZ}{{\mathbb{Z}}}
\newcommand{\db}{\overline{\partial}}
\newcommand{\ux}{\underline{x}}
\newcommand{\ur}{\underline{r}}
\newcommand{\uSigma}{\underline{\Sigma}}
\newcommand{\uR}{\underline{R}}

\newcommand{\un}{\underline{n}}
\newcommand{\uC}{\underline{C}}
\newcommand{\ukappa}{\underline{\kappa}}
\newcommand{\tkappa}{\tilde{\kappa}}
\newcommand{\tF}{\tilde{F}}
\newcommand{\dtheta}{\frac{\partial}{\partial \theta}}
\newcommand{\dH}{\frac{\partial}{\partial H}}
\newcommand{\dt}{\frac{\partial}{\partial t}}
\newcommand{\dQ}{\frac{\partial}{\partial Q}}

\begin{document}

\title{Singular Lefschetz pencils}
\authors{Denis Auroux\\Simon K Donaldson\\Ludmil Katzarkov}
\shortauthors{D Auroux, S\,K Donaldson and L Katzarkov}
\address{Department of Mathematics, Massachusetts Institute of Technology\\
Cambridge, MA 02139, USA}
\address{Department of Mathematics, Imperial College\\ London SW7 2BZ, 
United Kingdom}
\address{Department of Mathematics, University of Miami,
Coral Gables, FL 33124, USA\\ and Department of Mathematics, UC Irvine, Irvine, CA 92612, USA}

\address{\medskip\\\tt\mailto{auroux@math.mit.edu}, \mailto{s.donaldson@imperial.ac.uk}, 
\mailto{lkatzark@math.uci.edu}}

\asciiaddress{Department of Mathematics, Massachusetts Institute of 
Technology\\Cambridge, MA 02139, USA\\Department of Mathematics, 
Imperial College\\London SW7 2BZ, 
United Kingdom\\
Department of Mathematics, University of Miami,
Coral Gables, FL 33124, USA\\and Department of Mathematics, UC Irvine, 
Irvine, CA 92612, USA}
\asciiemail{auroux@math.mit.edu, s.donaldson@imperial.ac.uk, lkatzark@math.uci.edu}


\begin{abstract}
We consider structures analogous to symplectic Lefschetz pencils in the
context of a closed 4--manifold equipped with a ``near-symplectic'' structure
(ie, a closed 2--form which is symplectic outside a union of circles where
it vanishes transversely). Our main result asserts that, up to blowups, every
near-symplectic 4--manifold $(X,\omega)$ can be decomposed into (a) two symplectic Lefschetz
fibrations over discs, and (b) a fibre bundle over $S^1$ which relates
the boundaries of the Lefschetz fibrations to each other via a sequence of
fibrewise handle additions taking place in a neighbourhood of the zero
set of the 2--form. Conversely, from such a decomposition one can recover a
near-symplectic structure.
\end{abstract}
\asciiabstract{%
We consider structures analogous to symplectic Lefschetz pencils in
the context of a closed 4-manifold equipped with a `near-symplectic'
structure (ie, a closed 2-form which is symplectic outside a union of
circles where it vanishes transversely). Our main result asserts that,
up to blowups, every near-symplectic 4-manifold (X,omega) can be
decomposed into (a) two symplectic Lefschetz fibrations over discs,
and (b) a fibre bundle over S^1 which relates the boundaries of the
Lefschetz fibrations to each other via a sequence of fibrewise handle
additions taking place in a neighbourhood of the zero set of the
2-form.  Conversely, from such a decomposition one can recover a
near-symplectic structure.}

\maketitle


\section{Introduction}

The classification of smooth 4--manifolds remains mysterious, but that of
{\it symplectic\/} 4--manifolds is perhaps a little clearer. The purpose of
this article is to extend some of the techniques which have been developed
in the symplectic case to more general 4--manifolds.

Let $X$ be a smooth, oriented,
 $4$--manifold and let $\omega$ be a closed $2$--form on $X$. Then
 $\omega$ is a symplectic structure, compatible with the given orientation,
  if and only if $\omega^{2}>0$ everywhere on $X$. We are interested in
  relaxing this condition. Any form $\omega$ has, at each point of $X$, a
  rank which is $0,2$ or $4$. We consider forms with $\omega^{2}\geq 0$ and
  which do not have rank $2$ at any point: thus $\omega^{2}=0$  only at the
  set $\Gamma\subset X$ of  points
  where $\omega$ vanishes. The nature of this condition becomes clearer
  if we recall that the wedge-product defines a quadratic form of signature
  $(3,3)$ on $\Lambda^{2} \bR^{4}$. Locally we can regard a $2$--form as
  a map into $\Lambda^{2} \bR^{4}$ and the condition is that the image
  of the map only meets the null-cone at the origin. Suppose $\omega$ satisfies
  this condition and let $x$ be a point of the zero-set $\Gamma$. Thus
  there is an intrinsically defined derivative $\nabla \omega_{x}\co TX_{x}
  \rightarrow \Lambda^{2} T^{*}X_{x}$. The rank of $\nabla\omega_{x}$ can
  be at most $3$, since the wedge product form is nonnegative on the image.
  
  \begin{defn}
  A closed $2$--form on $X$ is a {\it near-symplectic structure} if $\omega^{2}\geq
  0$, if  $\omega$
  does not have rank $2$ at any point and if the rank of $\nabla \omega_{x}$
  is $3$ at each point $x$ where $\omega$ vanishes. \end{defn}
   
   It follows from this definition that the zero set $\Gamma$ of a near-symplectic
   form is a $1$--dimensional submanifold of $X$. The point of this notion
   is that, on the one hand, the form defines a {\it bona fide} symplectic
   structure outside this \lq\lq small'' set, while on the other hand these
   near-symplectic structures exist in abundance. 

   \begin{prop} \label{prop:ns-sdharm}
    Suppose $\omega$ is a near-symplectic form on $X$. Then there is a
    Riemannian metric $g$ on $X$ such that $\omega$ is a self-dual harmonic
    form with respect to $g$. Conversely, if $X$ is compact and 
    $b_{2}^{+}(X)\geq 1$ then for generic
    Riemannian metrics on $X$ there is a self-dual harmonic form which
    defines a near-symplectic structure. Moreover there is a dense subset
    of metrics on $X$ for which we can choose $\omega$ such that the cohomology
    class $[\omega]$ is the reduction of a rational class. \end{prop}

   This is essentially a standard result,
   and we give the proof in Section \ref{sec:converse}. It is also worth
   mentioning another existence result for near-symplectic forms, recently
   obtained by Gay and Kirby, in which the 2--form is constructed explicitly
   from the handlebody decomposition induced by a Morse function on $X$
   \cite{GK}. In any case, the point we
   wish to bring out, in formulating things the
   way we have, is that the near-symplectic condition has a meaning independent
   of Riemannian geometry. Indeed one
   can see this as the first case of a hierarchy of conditions, for a closed
   $2$--form on a $2n$--manifold, in which one imposes constraints on the way
   in which the form meets the different strata, by rank, of 
   $\Lambda^{2}\bR^{2n}$.
   
   Given the abundance of near-symplectic structures, it is natural
  to try to extend techniques from symplectic geometry to this more
  general situation. This is, of course, the starting point for
   Taubes' programme, studying the Seiberg--Witten equations and
   pseudo-holomorphic curves \cite{Ta1,Ta2}. This article runs entirely
   parallel
   to Taubes' programme, our aim being to extend some of the 
   \lq\lq approximately
   holomorphic'' techniques developed in \cite{Do1,Do3} to the near-symplectic
   case. More specifically, recall that any compact symplectic 
   $4$--manifold $(X,\omega)$ 
   (with rational class $[\omega]$) 
   admits a symplectic Lefschetz pencil. That is, there are disjoint,
   finite sets $A,B\subset X$  and a map  $f\co X\setminus
   A\rightarrow S^{2}$ which conforms to the following local models,
   in suitable oriented (complex) co-ordinates about each point $x\in X$.
   \begin{itemize}
   \item If $x\in A$ the model is $(z_{1},z_{2})\mapsto z_{1}/z_{2}$;
   \item If $x\in B$ the model is $(z_{1},z_{2})\mapsto z_{1}^{2}
   + z_{2}^{2}$; 
   \item For all other $x$ the model is $(z_{1},z_{2})\mapsto z_{1}$.
   \end{itemize}
   Although the map $f$ is not defined at $A$ (the ``base points'' of the
   pencil), the fibres
   $f^{-1}(p)$ can naturally be regarded as closed subsets of $X$
   by adjoining the points of $A$. The connection with the symplectic
   form $\omega$ is that these fibres are symplectic subvarieties,
   Poincar\'e dual to $k\omega$, for large $k$.
   
   Conversely, under mild conditions, a $4$--manifold  which admits
   such a Lefschetz pencil is symplectic \cite{Go2}. The main aim of this
   paper is to generalise these results to the near-symplectic case.
   To formulate our result, let $Y$ be any oriented $4$--manifold
   and let $\Delta \subset Y$ be a $1$--dimensional submanifold.
   We say that a map $f\co Y\rightarrow S^{2}$ has {\it indefinite
   quadratic singularities}
   along $\Delta$ if around each point of $\Delta$ we can choose local
   co-ordinates $(y_{0}, y_{1}, y_{2}, t)$ such that $\Delta$ is
   given by $y_{i}=0$ and the map $f$ is represented in suitable
   local co-ordinates on $S^{2}$ by
   $$    (y_{0},y_{1},y_{2},t)\mapsto y_{0}^{2}- \frac{1}{2}(y_{1}^{2}+
   y_{2}^{2}) + i t. $$
   \begin{defn}
   A singular Lefschetz pencil on $Y$, with singular set $\Delta$,
   is given by a finite set $A\subset Y\setminus \Delta$ and a
   map $f\co Y\setminus A\rightarrow S^{2}$ which has indefinite quadratic singularities
   along $\Delta$ and which is a Lefschetz pencil on $Y\setminus
   \Delta.$\end{defn}
   Given such a singular Lefschetz pencil we define the fibre over
   a point $p$ in $S^{2}$ in the obvious way, adjoining the points of
   $A$. Any such fibre is homeomorphic to the space obtained from
   a disjoint union of compact oriented surfaces by identifying
   a finite number of disjoint pairs of points. We refer to the
    image of one
   of these surfaces under the composite of the homeomorphism and
   the identification map as a {\it component} of the fibre. We
   can now state our main result.
   \begin{thm}\label{thm:general}
   Suppose $\Gamma$ is a $1$--dimensional submanifold of a compact oriented $4$--manifold $X$. Then
   the following two conditions are equivalent.
   \begin{itemize}
   \item There is a near-symplectic form $\omega$ on $X$,  with zero set $\Gamma$,
   \item There is a singular Lefschetz pencil $f$ on $X$ which has quadratic
   singularities along $\Gamma$, with the property  that there is a class $h\in
   H^{2}(X)$ such that $h(\Sigma)>0$ for every component $\Sigma$ of
   every fibre of $f$.  
    \end{itemize}
   \end{thm}
   This is a somewhat simplified statement, we actually prove rather
   more, in both directions. The general drift is, roughly, that
   there is  a correspondence between these two kinds of objects:
   near-symplectic forms and singular pencils. To state a more precise
   result, in one direction, we recall a result of Honda \cite{Ho2}. 
   Take $\bR^{4}$ with co-ordinates $(x_{0},x_{1}, x_{2}, t)$ and
   consider the 2--form 
   $$\Omega = dQ\wedge dt + *\, (dQ\wedge dt), $$
   where $Q(x_{0}, x_{1}, x_{2})= 
   x_{0}^{2} - \frac{1}{2}(x_{1}^{2}+x_{2}^{2})$ and $*$ is the
   standard Hodge $*$--operator on $\Lambda^{2}\bR^{4}$.  
  Let $\sigma_{-}\co \bR^{3}\rightarrow \bR^{3}$ be the map
  $\sigma_{-}(x_{0}, x_{1}, x_{2})=(-x_{0}, x_{1}, -x_{2})$. Define
  $\overline{\sigma}_{+}\co \bR^{4}\rightarrow \bR^{4}$ to be the translation
   $$\overline{\sigma}_{+}(\ux,t)= (\ux, t+2\pi)$$
   and let $\overline{\sigma}_{-}$ be the map
   $$ \overline{\sigma}_{-}(\ux,t)= (\sigma_{-}(\ux), t+2\pi)
  . $$
   The maps $\overline{\sigma}_{\pm}$ preserve the form $\Omega$ so we
   get induced forms on the quotient spaces. Let $N_{\pm}$ be the
   quotients of the tube $B^{3}\times \bR$ by $\overline{\sigma}_{\pm}$ with
   the induced near-symplectic forms. Then, according to Honda, if
   $\omega$ is any near-symplectic form on a $4$--manifold $X$ with
   zero set $\Gamma$  there is a Lipschitz homeomorphism $\phi$ of
   $X$ --- equal to the identity on $\Gamma$,  smooth outside $\Gamma$
   and supported in an arbitrarily small neighbourhood of $\Gamma$ --- such
   that $\phi^{*}(\omega)$ agrees with one of the two models $N_{\pm}$
   in suitable trivialisations of tubular neighbourhoods of each
   component of $\Gamma$. Replacing $\omega$ by $\phi^{*}(\omega)$
   we may suppose for most purposes that the form agrees with the standard models
   in these tubular neighbourhoods. Let $f_{\pm}\co N_{\pm}\rightarrow
   \bR\times S^{1}$ be the maps defined by $(Q,t)$ in the obvious
   way. 
   
   Suppose now that $\omega$ is a near-symplectic form with $[\omega]$
   an integral class in $H^{2}(X)$. Thus we may choose a complex
   line bundle $\mathcal{L}$ with connection over $X$ having curvature $-i\omega$.
   Given the choice of this connection we get, for each component of the  singular
   set, a holonomy in $U(1) \subset \bC$. It will be convenient
   to suppose  that all these
   holonomies are equal to $-1$. The more precise result we prove
   in one direction is: 
   \begin{thm}\label{thm:main}
   Suppose that $\omega$ is a near-symplectic form on $X$ equal
   to one of the standard models in neighbourhoods of the zero set
   $\Gamma$. Suppose that  
   $[\omega]= c_{1}(\mathcal{L})$ is integral and  that $\mathcal{L}$ has holonomy
   $-1$ around each component of\/ $\Gamma$. Then for all sufficiently
   large odd integers $k$ there is a singular
   Lefschetz pencil on $X$ such that \begin{itemize}\item
   the fibres are symplectic with respect to $\omega$;
   \item the fibres are  in the homology
   class dual to $kc_{1}(\mathcal{L})$;
   \item in sufficiently small 
   neighbourhoods
   of the components of the singular set, the map is equal to the composite of 
   one of the maps $f_{\pm}$ with a diffeomorphism taking $(-\delta, \delta)
   \times S^{1}$ to a neighbourhood of the standard equator
   in $S^{2}$.\end{itemize}    
   \end{thm}

\noindent
In the last part of the statement,
the diffeomorphism taking $(-\delta,\delta)\times S^1$ to a
neighbourhood of the equator is essentially the same for every component of
$\Gamma$, as will be clear from the proof. Hence, each
component of $\Gamma$ is mapped bijectively to the equator, and
there are well-defined ``positive'' and ``negative'' sides
of the equator, corresponding to $Q>0$ and $Q<0$ in a consistent manner
for all components.

   It is easy to deduce one half of Theorem \ref{thm:general} from Theorem
   \ref{thm:main}. Given any near-sympletic form we use Honda's result to get
   a new one  compatible with the standard models. Making a
   small deformation away from $\Gamma$ we can suppose that $[\omega]$
   is a rational class and then multiplying by a suitable integer
   we obtain an integral class, associated to a line bundle with
   connection. Making a further small deformation we can suppose
   that each of the holonomies around the components of $\Gamma$
   is a root of $z^{n}=-1$, for some large $n$. Then again replacing
   the line bundle by its $n$th power we fit into the hypotheses
   of Theorem \ref{thm:main}.
    
   The more precise result in the converse direction is the following:
 \begin{thm} \label{thm:converse}
   Let $X$ be a compact oriented 4--manifold, and let $f\co X\setminus A\to S^2$
   be a singular Lefschetz pencil with singular set\/ $\Gamma$ (ie, a smooth
   map described by the above local models in oriented local co-ordinates).
   If there exists a cohomology class $h\in H^2(X)$ such that $h(\Sigma)>0$
   for every component $\Sigma$ of every fibre of $f$, then $X$ carries a
   near-symplectic form $\omega$, with zero set $\Gamma$, and which makes
   all the fibres of $f$ symplectic outside of their singular points.
   Moreover, these properties determine a deformation class of
   near-symplectic forms canonically associated to $f$.
 \end{thm}

   In particular, if every component of every fibre of $f$ contains at least
   one base point, then the cohomological assumption automatically holds.
   In that case we can require $[\omega]$ to be Poincar\'e dual to the
   homology class of the fibre.
  
  The topology of singular Lefschetz pencils is made quite complicated by
  the presence of the singular locus $\Gamma$. Nonetheless, Theorem
  \ref{thm:main} leads to an interesting structure result for
  near-symplectic 4--manifolds. Namely, given a near-symplectic 4--manifold
  $(X,\omega)$ with $\omega^{-1}(0)=\Gamma$ and a singular Lefschetz pencil
  $f\co X\setminus A\to S^2$ such that $\Gamma$ maps to the equator as in
  Theorem \ref{thm:main}, after blowing up
  the base points we can decompose the manifold into:
  
  \begin{itemize}
  \item two
  symplectic Lefschetz fibrations over discs $f_\pm\co X_\pm\to D^2$,
  obtained by restricting $f$ to the preimages of
  two open hemispheres not containing the equator $f(\Gamma)$;
  \item the preimage $W$ of a neighbourhood of the equator.
  \end{itemize}

  The 4--manifold $W$ is a fibre bundle over $S^1$, whose fibre $Y$ defines a
  cobordism between the fibres $\Sigma_+$ and $\Sigma_-$ of $f_\pm$ (note
  that these need not be connected a priori),
  consisting of a succession of handle additions. Hence the cobordism
  $W$ relates the boundaries of $X_+$ and $X_-$ to each other via a sequence
  of fibrewise handle additions (or ``round handle'' additions), one for each
  component of $\Gamma$.

  The topology of $f$ can be described combinatorially in terms of
  (a) the monodromies of the Lefschetz
  fibrations $f_\pm$, which are given by products of positive Dehn twists
  in the relative mapping class groups of $(\Sigma_\pm,A)$, and (b)
  gluing data, which can be expressed eg, in terms of a coloured link
  on the boundary of one of the Lefschetz fibrations (see Section
  \ref{sec:topology}). This information determines $f$ completely if
  the identity components in $\mathrm{Diff}(\Sigma_\pm,A)$ are simply
  connected (eg, if $\Sigma_\pm$ both have genus at least $2$).

The paper is organised in the following way. 
Sections \ref{sec:ahtheory}--\ref{sec:completion} are devoted to
the proof of Theorem \ref{thm:main}. The proof rests on techniques of
approximately holomorphic
geometry: roughly speaking, the construction of maps which have the same
topological properties as holomorphic maps but in a context where the
underlying almost-complex structure is not integrable. 
In Section~\ref{sec:ahtheory} we
develop the techniques from this theory that we need, encapsulated into a
general result (Theorem \ref{thm:glob}), which may have other applications.
(As an aside
here, we mention that it would be interesting to compare our results with
the methods developed by Presas \cite{Presas} for symplectic manifolds with
contact boundary.)
The core of the paper lies in Sections \ref{sec:sec3}--\ref{sec:completion}.
Here we show that a 4--manifold
with a near symplectic form can be endowed with the geometrical structures
required to apply Theorem \ref{thm:glob}. Almost all of the work is
devoted to the geometry in a
standard model around the zero set, and we make extensive and explicit
calculations here. Again, these geometrical constructions could conceivably
be of interest in other contexts. (One can also compare with the detailed
study by Taubes of
other geometrical phenomena in the same local model \cite{Ta2}.)
In Section \ref{sec:converse} we prove the converse result,
Theorem \ref{thm:converse}, together with
Proposition \ref{prop:ns-sdharm} above. 
Section \ref{sec:topology}  begins the exploration of the topological
aspects  of singular Lefschetz pencils and their monodromy data.
  
\rk{Acknowledgements}

DA was partially supported by NSF grant DMS-0244844.
LK was partially supported by NSF grant DMS-9878353 and NSA grant
H98230-04-1-0038.

   \section{Approximately holomorphic theory}\label{sec:ahtheory}

In the symplectic case, the construction of Lefschetz pencils 
relies on the existence of various structures which are the basic building
blocks of ``approximately holomorphic geometry'', which we now review
briefly and informally (precise statements are given below in a more general
context). Let $(X,\omega)$ be a compact symplectic $2n$--manifold, equipped
with a compatible
almost-complex structure $J$, and let ${\cal L}\to X$ be a hermitian line
bundle with a connection having curvature $-i\omega$.
Then, for any $\epsilon>0$, the following properties hold up to a suitable
rescaling, replacing $\omega$ by $k\omega$ and ${\cal L}$ by
${\cal L}^{\otimes k}$ for some large integer $k$ (such that $k^{-1/2}\ll
\epsilon$) \cite{Do1,Do3}:
\medskip

(1)\qua Near every point of $X$ there exist ``approximately holomorphic''
co-ordin\-ate charts, ie, local complex co-ordinates in
which the almost-complex structure $J$ differs from the standard complex
structure by at most a fixed multiple of $\epsilon$, and the symplectic 
form $\omega$ is uniformly bounded.
\medskip

(2)\qua  For every point $p\in X$ there exist $n+1$ ``localised'' sections
$\sigma_0,\dots,\sigma_n$ of ${\cal L}\to X$, and a real-valued
function $F$ which decreases exponentially fast away from $p$, such that

\begin{itemize}
\item $\sigma_i$ and its covariant derivatives are bounded by $F$;

\item $\overline{\partial} \sigma_i$ and its covariant derivatives are
bounded by $\epsilon F$;

\item $|\sigma_0|$ is bounded below near $p$, and the ratios
$\sigma_i/\sigma_0$ define local complex co-ordinates near $p$.
\end{itemize}

\noindent These key ingredients
make it possible to construct a pair of sections $\sigma_0,\sigma_1$ of
$\cal L$  such that the map $\sigma_1/\sigma_0$ is a
symplectic Lefschetz pencil. More precisely, $\sigma_0$ and $\sigma_1$
are obtained as linear combinations of the above-mentioned ``localised''
sections, chosen in a manner which ensures that $\sigma_0$, $\sigma_1$,
and $\partial(\sigma_1/\sigma_0)$ satisfy suitable uniform transversality properties (ie,
whenever one of these quantities vanishes, its derivative is surjective
and satisfies a uniform {\it a priori} lower bound)~\cite{Do3}.

In the near-symplectic case, we can use the same methods away from the
zero set $\Gamma$ of the near-symplectic form, while over a small
neighbourhood of $\Gamma$ we can construct the pencil explicitly from a
local model. The difficulty comes from the intermediate region. To adapt
the machinery to this situation, we consider the non-compact symplectic
manifold $Z=X\setminus \Gamma$,
equipped with a suitably rescaled symplectic form $k\omega$, and two
compact subsets $K\subset Z_0\subset Z$ (the complements of two
concentric tubular neighbourhoods of $\Gamma$ in $X$). Our goal will be
to show that, for a suitable choice of almost-complex structure on
$Z$, the following properties
hold (see below for more precise statements) and imply
Theorem~\ref{thm:main}:%

(1)\qua  near every point of $Z_0$ there exist local approximately holomorphic
co-ordinate charts;

(2)\qua  near every point of $K$ there exist localised 
sections of $\cal L$, with support contained in $Z_0$, and with the same
properties as in the symplectic case;

(3)\qua  there are sections $\sigma_0,\,\sigma_1$ of ${\cal L}\to Z$ such that
$\sigma_1/\sigma_0$ is a Lefschetz pencil outside of $K$, and satisfying
appropriate uniform bounds over $Z_0$ and transversality estimates over the
intermediate region $Z_0\setminus K$.
\medskip

In the rest of this section, we give precise formulations of these
three properties (``hypotheses''), and show how they can be used to
construct a Lefschetz pencil. We place ourselves in a more general
setting, although the reader may wish to keep in mind the above motivation.

   Let $(Z,\omega)$ be a symplectic $2n$--manifold, not necessarily
   compact and let $J$ be a  compatible almost-complex structure
   on $Z$. Suppose we have a hermitian line bundle $\mathcal{L}\rightarrow
   Z$ with a connection having curvature $-i \omega$. We also suppose
   that we have given compact subsets $Z_{0}$ and $K$ of $Z$, such that
   $Z_{0}$ contains a neighbourhood of $K$. 
    We wish to formulate three ``hypotheses'' bearing on
   various data in this situation, involving
   certain numerical parameters. One collection of parameters will
   be denoted $C_{1}, C_{2}, \dots$ which we abbreviate to a single
   symbol $\uC$. These give bounds on the geometry of the set-up:
   the precise number of parameters $C_{i}$ is unimportant, it would
   probably be possible to reduce them to a single constant $C$,
   but this would mean considerable loss of accuracy if one was
   actually interested in implementing the proof numerically. The
   important parameter is a small number $\epsilon$ which, roughly,
    measures
   the deviation from holomorphic geometry. In the third hypothesis
   we will introduce three parameters
    $\kappa_{1},\kappa_{2}, \kappa_{3}$ which we sometimes denote
    by $\ukappa$. These are a measure of transversality of certain data.    
   
   \begin{hyp} Hypothesis $H_{1}(\epsilon, \uC)$.
   
   For each point $p$ of $Z_{0}$ there is a  co-ordinate chart
   $\chi_{p}\co B^{2n}\rightarrow Z$ centred on $p$ such that 
   \begin{itemize}
   \item The pull-back
   $\chi_{p}^{*}(J)$ of the almost-complex structure on $Z$ is close
   to the standard structure $I$ on $B^{2n}\subset \bC^{n}$ in that
   $$ \Vert \chi_{p}^{*}(J)-I\Vert_{C^{r}}\leq C_{1} \epsilon. $$
   \item The pull-back of the symplectic form satisfies uniform
   bounds $$\Vert \chi_{p}^{*}(\omega)\Vert_{C^{r}}\leq C_{2}$$ and
   $ \chi_{p}^{*}(\omega)^{n}\geq
   C_{2}^{-1}$. \end{itemize}
   Here $r$ is a fixed integer, $r=3$ will do.
 \end{hyp}
          We call such a chart an ``approximately holomorphic chart'',
          where of course the notion depends on the parameters $\epsilon,
          C_{i}$.

\medskip          
          
\noindent {\bf Remark}\qua In essence, this hypothesis asserts that the manifold
          has bounded geometry and that the norm of the Nijenhuis tensor
          is $O(\epsilon)$.

\medskip          
          
          Before stating the next hypothesis we formulate a definition.
          Let $U\subset V\subset W$ be subsets of $Z$ (with $U$ open) and let $F$ be a
          positive function on~$Z$.
          \begin{defn}
          An $F$--localised, $\epsilon$--holomorphic system over $U$, relative
          to $V$ and $W$, consists of $n+1$ sections $\sigma_{0},\dots,\sigma_{n}$
          of $\mathcal{L}\rightarrow Z$ such that \begin{itemize}
          \item The support of any section  $\sigma_{i}$ is contained in
          the interior of $W$;
          \item $\vert \nabla^{p} \sigma_{i}\vert \leq F$ throughout $Z$,
          for $p\leq r$ and all $i$;
          \item $\vert \nabla^{p}\db \sigma_{i}\vert \leq \epsilon F$ in
          $V$ for $p\leq r$ and all $i$;
           \item $\vert \sigma_{0}\vert \geq 1$ in $U$~--~this means that
           we can define a map $f\co U\rightarrow \bC^{n}$ by the ratios
           $\sigma_{i}/\sigma_{0}$;
           \item The Jacobian of $f$ (defined using the volume form
            $\omega^{n}$
           on $Z$) is not less than $1$.
            \end{itemize}
            \end{defn}     
            
            Now we can state the following:
           \begin{hyp} 
            Hypothesis $H_{2}(\epsilon, \uC)$    
          
          There is a finite collection of approximately holomorphic charts
          $\chi_{i}$, $i=1,\dots, M$ mapping to balls $B_{i}$ contained in $Z_{0}$ such that
          \begin{itemize}
          \item For a fixed $\lambda= \frac{C_{3}}{1+ C_{3}}$, 
          the balls $\lambda B_{i}= \chi_{i}(\lambda B^{2n})$ cover $K$.
          We define $K^{+}$ to be the union of the balls $B_{i}$.
          \item There are positive functions $F_{i}$ on $Z$ and for each
          $i$ an $F_{i}$--localised, $\epsilon$--holomorphic system over
          $B_{i}$, relative to $K^{+}, Z_{0}$. For each point $q$ in
          the support of any section making up this system there is an approximately holomorphic
          chart centred on $q$ with image contained in $Z_{0}$.
          \item For each point $p$ of $Z$,
          $$  \sum_{i} F_{i}(p) \leq C_{4}. $$
          \item For all $D>1$ we can divide the set $\{1,\dots, M\}$ into
          $N=N(D)$ disjoint subsets $I_{1},\dots,I_{N}$ where
          $$  N(D)\leq C_{5} D^{C_{6}}, $$
          and if $p$ is contained in a ball $B_{i}$ for $i\in I_{\alpha}$
          then
          $$  \sum_{j\in I_{\alpha},j\neq i} F_{j}(p) \leq C_{7} e^{-D}. $$   \end{itemize}     
          \end{hyp}

\noindent {\bf Remark}\qua In essence, this hypothesis states that associated
          to each point there are
           approximately holomorphic sections of the line bundle which
           on the one hand decay rapidly away from the point, and on the
           other hand give an approximately holomorphic projective embedding of a neighbourhood of
           the point.
          
          \medskip
              
          The third hypothesis bears on a pair of sections $\sigma_{0},\sigma_{1}$
          which should be thought of as giving a model for a pencil outside
          $Z_{0}$. 

          Recall that, given a $\mathbb{CP}^1$--valued map $F$ defined over an open
          subset of $Z$ and a constant $\kappa>0$, we say that {\em $\partial F$ is $\kappa$--transverse to $0$}
          if at any point where $\vert \partial F\vert <\kappa$ the covariant
          derivative $\nabla \partial F$ is invertible and the inverse
          has norm less than $\kappa^{-1}$.

          \begin{hyp}
          Hypothesis $H_{3}(\epsilon, \kappa_{1},\kappa_{2},\kappa_{3},\uC)$.
          
          There are sections $\sigma_{0},\sigma_{1}$ of $\mathcal{L}\rightarrow Z$
          such that
          \begin{itemize}
          \item $F=\sigma_{1}/\sigma_{0}$ is a topological Lefschetz pencil over
          $Z\setminus K$, with symplectic fibres.
          \item $\vert \nabla^{p} \sigma_{i}\vert \leq C_{8}$ in $Z_{0}$,
          for $p\leq r$.
          \item $\vert \nabla^{p} \db \sigma_{i}\vert \leq C_{9}\epsilon
          $ in $K^{+}$.
          \item $\vert \sigma_{0}\vert^{2} + \vert \sigma_{1}\vert^{2}\geq
          C_{10}^{-1}$ in $Z_{0}\setminus K$; thus $F=\sigma_{1}/\sigma_{0}$
          defines a map from $Z_{0}\setminus K$ to the Riemann sphere
         $S^{2}$. 
          \item The complex-linear component $\partial F$ of the derivative
          of $F$ is $\kappa_{1}$--transverse to $0$ throughout $Z_{0}\setminus
          K$.
          \item $\vert \db F \vert \leq \max(  \epsilon \kappa_{2}, \vert \partial
          F\vert - \kappa_{3})$ throughout $Z_{0}\setminus K$
           \end{itemize}
          \end{hyp}          
          
          With all this in place we can state our  general theorem.
          \begin{thm} \label{thm:glob}
          There is a universal function
           $\epsilon_{0}(\ukappa,\uC)$ with
          the following property.
          If we have data satisfying hypotheses $H_{1}(\epsilon, \uC),
          H_{2}(\epsilon,\uC), H_{3}(\epsilon,\ukappa,
          \uC)$ and if 
           $\epsilon\leq \epsilon_{0}(\ukappa,\uC)$
          then there is a topological Lefschetz pencil on $(Z, \omega)$ with
          symplectic fibres, equal
          to $\sigma_{1}/\sigma_{0}$ outside $Z_{0}$.
          \end{thm}

          We will not say much about the proof of Theorem \ref{thm:glob}, which would
          essentially repeat the whole of the paper \cite{Do3} (see also
          \cite{Do1}, \cite{Au1}, \cite{Au2}).
          While there are no new ideas involved in the proof, the theorem
          extends the previous results in two different directions. On
          the one hand the theorem is a ``relative'' version of the
          previous results, extending a Lefschetz pencil which is already
          prescribed over a subset of the manifold. On the other hand, the dependence on parameters is made more explicit:
          in the earlier results the parameter $\epsilon$ is essentially
          $k^{-1/2}$ where one works with a {\it fixed} almost complex
          structure but scales the symplectic form by a factor $k$. The
          new result allows us to vary the almost complex structure at
          the same time as $k$, which will be one of the main ideas in our
          construction.

           We outline the proof of Theorem \ref{thm:glob}.  
           Introduce a parameter $c\in (0,1)$ and
          consider modifying the sections $\sigma_{0},
          \sigma_{1}$ to
          $$ \tilde{\sigma}_{0}= \sigma_{0} + \sum a_{j} s_{j}\ , 
          \tilde{\sigma}_{1} = \sigma_{1} + \sum b_{j} s_{j},
        $$ where the $s_{j}$ run over all the sections comprising the systems
        provided by Hypothesis $H_{2}$ and the coefficients $a_{j}, b_{j}$
        are complex numbers to be chosen, with the constraint that
         $$\vert a_{j}\vert, \vert
        b_{j}\vert \leq c. $$ The arguments of \cite{Do3} show that for any fixed
        $c$ and for small enough $\epsilon$ we can choose the coefficients
        such that $\tilde{F}= \tilde{\sigma}_{1}/\tilde{\sigma}_{0}$ is
        close to being a symplectic Lefschetz pencil over $K$, in that we can
        find a set of disjoint balls of radii $O(\epsilon)$ and obtain a Lefschetz
        pencil over $K$ by modifying $\tilde{F}$ inside these balls (in
        order to obtain the desired local model at the critical points).
        Since the sections
        $s_{j}$ are supported in $Z_{0} $ the map $\tilde{F}$ agrees with
        the model pencil outside $Z_{0}$. The new issue has to do with
        the intermediate region $Z_{0}\setminus K$, where we argue as follows.
       
        Suppose that a map  $\tilde{F}$ obtained by the
         procedure above satisfies
        \begin{itemize}
        \item $\partial \tilde{F}$ is $\tkappa_{1}$--transverse to $0$,
        \item
        $ \vert \db \tilde{F}\vert \leq \max (\nu,
         \vert \partial
        \tilde{F}\vert - \tkappa_{3}) , $
\end{itemize} for some $\nu, \tkappa_{1}, \tkappa_{3}>0$. By construction
we will also have bounds
    $$ \vert \nabla^{p} \tilde{F} \vert \leq C, $$
    for $p\leq 3$ and some fixed $C$. We claim that there is a $\nu_{0}$
    depending only on $C, \tkappa_{1}, \tkappa_{3}$ such that 
    if $\nu\leq \nu_{0}$ the map $\tilde{F}$ can be modified over a number
    of small disjoint balls to yield a symplectic Lefschetz fibration.
    
    By construction, the map $\tF$ agrees with the model $F$ outside the support
    of the $s_{j}$ and by Hypothesis $H_{2}$ we have a good co-ordinate
    chart centred on any point $q$ in the union of these supports. If
    $\vert \db \tF\vert < \vert \partial \tF\vert$ at $q$ then $\tF$ is a
    fibration with symplectic fibres near $q$. If on the other hand $\vert
    \db \tF \vert \geq \vert \partial \tF \vert$ then  we must have
    $  \vert \partial \tF\vert\leq \nu $ at $q$. It follows from the transversality
    estimate on $\partial \tF$ that if $\nu$ is sufficiently small  compared with
    $\tkappa_{1}$ then $q$ is close to a zero  of $\partial \tF$: more
    precisely we can find such a zero $p$ at a distance $O(\nu/\tkappa_{1})$
    from $q$. Adjusting constants slightly, we can suppose that there is
    a good co-ordinate chart centred at this point $p$ and contained in
    $Z_{0}$. 
        
        Now we clearly have $\vert \db \tF\vert\leq \nu$ at $p$. We
     claim that the derivative $\vert \nabla \db \tF\vert$ is $O(\nu^{1/2})$
     at $p$. To see this,  suppose that $\vert \nabla \db \tF(p)\vert=A$.
     Then for any small $r$,
     we can find a point $p'$ at distance $r$ from $p$ with
     $\vert \db\tF(p')\vert \geq Ar-\frac{C}{2}r^2$. If $r$ is small
     compared with $\tkappa_{3}/C$ we have $\vert \partial \tF\vert < \tkappa_{3}$
     at $p'$ so it follows that $\vert \db \tF\vert(p')\leq \nu$. Combining
     the inequalities gives $A\le \frac{\nu}{r}+\frac{Cr}{2}$. Taking $r$ of
     the order of $\nu^{1/2}$ we obtain the desired bound $A=O(\nu^{1/2})$. 
     Now considering the Taylor series
     of $\tF$ at $p$ just as in \cite{Do3}, Section 2, we see that $\tF$ can be
     modified in a ball of radius $\rho$ to obtain a new map which is a Lefschetz fibration
     over the ball provided we can find a  radius $\rho$ which satisfies
     $$    \nu^{1/2} \ll \rho \ll \tkappa_{1}/C.$$
     This will be possible if $\nu$ is small and we see that moreover the original
     point $q$ will lie inside the ball. So we conclude that, after making
     these modifications we obtain the desired fibration.
     
     With this discussion in place we now return to complete the proof.
     Recall that, under our hypotheses, we do not have any $\epsilon$ bound
     on $\db s_{j}$ outside $K^{+}$. What we do have is a bound
     $$  \vert \nabla^{r}( \tF - F)  \vert \leq B c $$
     for a suitable constant $B$. It follows that if $c$ is sufficiently
    small then $\partial \tF $
     is $\kappa_{1}/2$--transverse to $0$. Similarly
     $$ \vert \db \tF\vert \leq \vert \db F \vert + Bc \leq \max( Bc+ \epsilon
     \kappa_{2}, \vert \partial \tF \vert + 2Bc -\kappa_{3}) . $$
        We set $\tkappa_{1} = \kappa_{1}/2$ and choose $c$ so small that
        $2Bc\leq \kappa_{3}/2$. Then we can take $\tkappa_{3}=\kappa_{3}/2$.
        Thus we have a $\nu_{0}=\nu_{0}(\tkappa_{1}, \tkappa_{2})$, as
        above. Now we also choose $c$ so small that $ Bc\leq \nu_{0}/2$.
        Then if $\epsilon$ is so small that $\epsilon \kappa_{2}\leq \nu_{0}/2$
        we achieve the desired properties for our function $\tF$.

\section{Definition of the almost-complex structure}\label{sec:sec3}

\subsection{Set-up} 
 
         In this section we put our problem in the general framework considered
         in Section \ref{sec:ahtheory}. To simplify notation we will consider a case where
         the singular set has just one component and the model is $N_{+}$.
         (At the end of the proof, in \S \ref{ss:oddcase} below, we discuss the easy extensions
         to the general case.) Thus we suppose that $X$ is a compact Riemannian
         $4$--manifold containing an isometrically embedded
          copy  $N\subset X$ of the standard model $N_{+}$ and that $\omega$
           is a closed self-dual $2$--form on $X$ which
          is equal to the standard form $\Omega$ in $N_{+}$ and  which
          does not vanish
          outside $N_{+}$. We suppose that there is a unitary line bundle
          with connection $\mathcal{L}\rightarrow X$ having holonomy $-1$ around the
          zero set and with curvature $-i\omega$. For large odd integers
          $k$ we consider the line bundle $\mathcal{L}^{\otimes k}$ with curvature $-ik\omega$.
          Clearly the standard form $\Omega$ on $\bR^{4}$ scales with weight~$3$.
          Thus we can identify the pair $(N,k\omega)$  with the
          form induced by $\Omega$ on the quotient of $B^{3}(k^{1/3})\times \bR$ under the translations
          $t\mapsto t+2\pi\bZ k^{1/3}$, where $B^{3}(k^{1/3})$ is the ball
          in $\bR^{3}$ of radius $k^{1/3}$. We will denote this form again
          by $\Omega$. It is convenient to put
           $\epsilon = k^{-1/3}$; this is the essential parameter in the construction
           which will eventually be made very small. Throughout the proof
           our attention will be focussed on this region $N$ on which we
           take our standard co-ordinates $(x_{0}, x_{1}, x_{2}, t)$ (so
           $\vert \ux\vert\leq \epsilon^{-1}$). We
           recall that $\Omega$ is given by
              \begin{equation} \Omega = ( 2x_{0} dx_{0} - x_{1} dx_{1} - x_{2}dx_{2})
              \wedge dt + 2x_{0}dx_{1}\wedge dx_{2} - x_{1} dx_{2}\wedge dx_{0} - 
              x_{2} dx_{0}\wedge dx_{1}.\label{eq:om}\end{equation}
              So $$\Omega^{2} = (4x_{0}^{2} + x_{1}^{2}+ x_{2}^{2})\, dx_{0}
              \wedge dx_{1}\wedge dx_{2}\wedge dt. $$
              It will be convenient to write
              \begin{equation}
                p= (4x_{0}^{2}+x_{1}^{2} + x_{2}^{2})^{1/4}, \label{eq:p}\end{equation}
              so $\Omega^{2}$ is $p^{4}$ times the standard volume form.
  
           To match up with the set-up in Section \ref{sec:ahtheory}, we let $K\subset X\setminus
           \Gamma$ be the subset corresponding to $\vert \ux\vert \geq
           10$ and let $X_{0}$ be the subset corresponding to $\vert \ux\vert\geq
           1$.             
                                 
          The great benefit for us given by Honda's result \cite{Ho2}, reducing to
          this standard model, is that there are two obvious symmetries:
          translation in the $t$--direction and rotation in the $(x_{1}, x_{2})$ plane.
          We use the standard polar co-ordinates $(r,\theta)$ in
          the $(x_{1}, x_{2})$ plane and we define
          \begin{equation}   H= x_{0} r^{2}. \label{H}\end{equation}
          Then one readily checks that $H$ is the Hamiltonian for the rotation
          action and that
          $$\Omega = dQ\wedge dt + dH \wedge d \theta . $$
          Recall here that $Q$ is the quadratic form
          \begin{equation} Q(\ux)= 
          x_{0}^{2} - \frac{1}{2}(x_{1}^{2}+x_{2}^{2}). \label{eq:Q}\end{equation}
          In these $(Q,t,H, \theta)$ co-ordinates the Euclidean co-ordinate
           $x_{0}$ is defined implicitly as the root of the cubic equation
           \begin{equation}   x_{0}^{3} - Q x_{0} - \frac{H}{2}  = 0,
           \label{eq:cubic}\end{equation} 
           having the same sign as $H$.

         We want to define a suitable almost-complex structure $J$ on $X\setminus
         \Gamma$. This structure will depend on the parameter $\epsilon$.
          It is a standard fact that the compatible almost-complex
         structures on an oriented Riemannian $4$--manifold are parametrised
         by the unit self-dual $2$--forms, so we have one structure $J_{0}$ 
          corresponding to the form $\frac{\omega}{\vert \omega\vert}$,
          which is smooth away from $\Gamma$. In our co-ordinates on $N$
          this structure $J_{0}$ can be described as follows. We let
          $\un$ be the unit vector field on $\bR^{3}$
           $$    \un= p^{-2}( 2x_{0}, -x_{1}, -x_{2}). $$
           Then $J_{0}$ is characterised by the conditions that
           $$J_{0}(\un) = \frac{\partial}{\partial t}\ \ ,\ \  J_{0}(\frac{\partial}{\partial
           t}) = -\un, $$
           while on the orthogonal plane $\un^{\perp}$ in $\bR^{3}$, $J_{0}$
           is given by the standard rotation by $\pi/2$ (with orientation
           fixed by that of $\un$). Notice that $\un$ is the normalised
           gradient vector field of the quadratic   function $Q$ on $\bR^{3}$,
           so the planes $\un^{\perp}$ are tangent to the family of real
           quadric surfaces qiven by the level sets $\{Q(\ux)=q\}$ of $Q$.
            Thus these
           quadric surfaces are {\it complex curves} for the almost-complex
           structure $J_{0}$. More precisely, we have a 2--parameter family
           $\Sigma_{q,t}$ of Riemann surfaces in $N$.
           
           The almost-complex structure $J$ we want to use is a modification
           of $J_{0}$. We set
           \begin{equation}  J(\un) =p^{2}\psi^{-2} \dt  \ \    ,\ \
           J(\frac{\partial}{\partial t}) =   -p^{-2}\psi^{2}\un;
           \label{Jpsi}\end{equation}
           where $\psi=\psi_{\epsilon}(\ux)$ is a function which we will
           specify shortly. On the orthogonal plane $\un^{\perp}$ we define
           $J$ to be the same as $J_{0}$, thus the $\Sigma_{q,t}$ are still
           complex curves for the almost-complex structure $J$. We require that the function
           $\psi$ be equal to $p$ once $\vert \ux\vert\geq \epsilon^{-1}= k^{1/3}$ so we
           can extend $J$ over the whole of $X$ by the standard structure
           $J_{0}$. The form $k \omega$ and the almost-complex structure
           $J$ define a Riemannian metric $g=g_{\epsilon}$ on 
           $X\setminus \Gamma$ in the standard way: outside $N$ this is just the original
           metric scaled by a factor $k\frac{|\omega|}{\sqrt{2}}$. 
           
            In terms of the $(Q,t,H, \theta)$ co-ordinates, the almost
            complex structure $J$ in the $(\dQ,\dt)$ plane is given by 
            \begin{equation}
               J(\dQ)= \psi^{-2} \dt \ ,\ J(\dt)= -\psi^{2} \dQ. \label{eq:Jform}
\end{equation}
            Writing the almost-complex structure in the $\dH,\dtheta$--plane
            explicitly is equivalent to finding the conformal structure
             induced on the
            quadric surfaces~-- which is just the structure induced from
            the embedding in $\bR^{3}$. A short calculation, which we leave
            as an exercise for the reader, shows that the metric $g$ is
            given in these co-ordinates by
            \begin{equation} g= \psi^{-2} dQ^{2} + \psi^{2} dt^{2} + 
            p^{-2} r^{-2} dH^{2} + p^{2} r^{2} d\theta^{2}.  \label{eq:metric}\end{equation}
Thus \begin{equation}
J(\dH)= p^{-2} r^{-2} \dtheta\ \ ,\ \ J(\dtheta)= -p^{2} r^{2} \dH   . 
\label{eq:conf}\end{equation}
We will now specify the function $\psi$ and hence the almost complex structure.
It is convenient to make $\psi$ a function of $p$, depending also on the
parameter $\epsilon$. Notice that $p$ is essentially equivalent to the
square root of the Euclidean norm:
$$    \vert \ux\vert \leq p^{2} \leq 2\vert \ux \vert . $$
\begin{lem}\label{lem:defpsi}
There  are constants $c_{r}$ such that for all sufficiently small 
$\epsilon$ we can find a smooth, positive, non-decreasing, function $\psi(p)$ on the interval
$[1,\epsilon^{-1/2}]$ with following properties:
\begin{itemize}
\item $\psi(p)=\epsilon$ if $p\leq \frac{1}{2} \epsilon^{-1/2}$;
\item $\psi(p)= p$ if $p\geq \frac{9}{10} \epsilon^{-1/2}$;
\item $\psi(p)\leq c_{0} p$;
\item $ \psi(p) \leq c_{0} \epsilon p^{4}$;
\item $\vert \frac{\psi^{(r)}}{\psi}\vert \leq c_{r}  \epsilon
p^{2r}$ (where $\psi^{(r)}$ denotes the r-th derivative of $\psi$.)
\end{itemize}
\end{lem}

\begin{figure}[ht!]
\setlength{\unitlength}{4mm}
\centering
\begin{picture}(12,10.8)(-1,-1)
\put(0,0){\vector(0,1){9.8}}
\put(0,0){\vector(1,0){11}}
\put(0,0.5){\line(1,0){4}}
\put(8.5,7.5){\line(1,1){2}}
\put(6,5){\line(1,0){1}}
\qbezier[120](4,0.5)(5,0.5)(5.3,2.75)
\qbezier[80](5.3,2.75)(5.5,5)(6,5)
\qbezier[60](7,5)(7.5,5)(8,6.5)
\qbezier[50](8,6.5)(8.25,7.25)(8.5,7.5)
\put(4,-0.1){\line(0,1){0.2}}
\put(4,-0.2){\makebox(0,0)[ct]{\small $\frac{1}{2}\epsilon^{-\frac12}$}}
\put(7,-0.1){\line(0,1){0.2}}
\put(7,-0.2){\makebox(0,0)[ct]{\small $\frac{3}{4}\epsilon^{-\frac12}$}}
\put(10,-0.1){\line(0,1){0.2}}
\put(10,-0.2){\makebox(0,0)[ct]{\small $\epsilon^{-\frac12}$}}
\qbezier[3](4,0)(4,0.25)(4,0.5)
\qbezier[30](7,0)(7,2.5)(7,5)
\qbezier[54](10,0)(10,4.5)(10,9)
\put(-0.1,0.5){\line(1,0){0.2}}
\put(-0.2,0.5){\makebox(0,0)[rc]{\small $\epsilon$}}
\put(-0.1,5){\line(1,0){0.2}}
\put(-0.2,5){\makebox(0,0)[rc]{\small $\frac{1}{2}\epsilon^{-\frac12}$}}
\put(-0.1,9){\line(1,0){0.2}}
\put(-0.2,9){\makebox(0,0)[rc]{\small $\epsilon^{-\frac12}$}}
\qbezier[36](0,5)(3,5)(6,5)
\qbezier[60](0,9)(5,9)(10,9)
\end{picture}
\caption{The function $\psi(p)$}
\label{fig:graphpsi}
\end{figure}
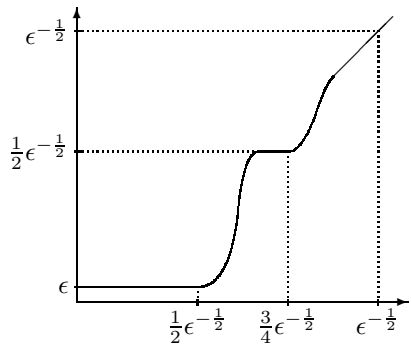

To prove the Lemma we give an explicit construction. Choose a smooth function on $[0,1]$ equal to
$0$ for small values and to $1$ for values near $1$. Using this in an
obvious way, we define for any $T>1$ a function $\alpha_{T}$, equal to
$1$ on the interval $[1,T]$ and supported in $(0,T+1)$. Likewise we choose
a smooth function $g(t)$, equal to $t$ for $t\geq \frac{9}{10}$ and to $\frac{1}{2}$
for $t\leq \frac{3}{4}$. For fixed $T$, let $f$ be the solution of the differential equation
$$    \frac{df}{dt}= \alpha_{T} f $$
with $f(t)=1$ for $t\leq 0$. Thus $f$ takes a constant value $L(T)$ say
for large values of $t$ (that is, for $t\geq T+1$). Clearly $L$ is approximately
$e^{T}$ for large values of $T$. Given a small $\epsilon$ we choose
$T$ so that $L=\frac{1}{2} \epsilon^{-3/2}$. Thus this $T=T(\epsilon)$
is much less than $\epsilon^{-1/2}$ for small $\epsilon$: we can assume that
  $T<\frac{1}{4} \epsilon^{-1/2}$.  Now define
$$  \psi_{0}(p) = \epsilon f(p-\frac{1}{2} \epsilon^{-1/2}). $$
Thus $\psi_{0}(p) = \epsilon$ for $p\leq \frac{1}{2}\epsilon^{-1/2}$ and
$\psi_{0}(p)= \frac{1}{2} \epsilon^{-1/2}$ for $p\geq \frac{3}{4}\epsilon^{-1/2}$.
Next define
$$\psi_{1}(p) = \epsilon^{-1/2} g (\epsilon^{1/2} p).$$
Thus $\psi_{1}(p) = p$ for $p\geq \frac{9}{10}\epsilon^{-1/2}$ and $\psi_{1}(p)$ takes
the same constant value $\frac{1}{2} \epsilon^{-1/2}$ as does $\psi_{0}$
for $p$ near $p_{0}=\frac{3}{4}\epsilon^{-1/2}$. So finally we define $\psi$
to be equal to $\psi_{0}$ for $p\leq p_{0}$ and to $\psi_{1}$ for $p\geq
p_{0}$ (see Figure \ref{fig:graphpsi}).

It is straightforward to check that this function satisfies the requirements
of the Lemma.

We now fix the almost complex structure to be the one defined by any function
$\psi$ which satisfies the requirements of Lemma \ref{lem:defpsi}, for example the function
constructed above.
\begin{prop}\label{prop:hyp1}
      There are constants $C_{1},C_{2}$ such that the symplectic
       manifold $(X\setminus \Gamma, k\omega)$ with the prescribed subsets
       $K\subset\subset X_{0}$ and  almost complex structure $J$ depending
       on $\epsilon= k^{-1/3}$ satisfies
       Hypothesis $H_{1}(\epsilon, \uC)$ for all large enough~$k$.
       \end{prop} 
 
          This is proved in Subsection \ref{ss:proofprop2} below.
          The essential idea of the proof  is the following.
          Away from $\Gamma$ what we have is just the familiar \lq\lq flattening''
          of the manifold by rescaling. The region $N$ is foliated by the
           Riemann surfaces $\Sigma_{q,t}$ and the
          almost-complex structure gives vector fields $\dt$ and $\dQ$ transverse
          to these. If the flow defined by these vector fields preserved
          the conformal structure of the Riemann surfaces we would have
          an integrable  structure and we could introduce genuine
           local holomorphic
          co-ordinates. The flow by $\dt$ obviously preserves the conformal
          structure, so the whole difficulty comes from the distortion
          in the conformal structure appearing in the flow of $\dQ$. However,
          the almost-complex structure and resulting metric $g$
           have been arranged so that the small
           parameter $\epsilon$ makes the length of $\dQ$ very large so,
           measured 
           with respect to this metric, the conformal distortion is very
 small and we can find approximately holomorphic co-ordinates.
 
 \subsection{Holomorphic co-ordinates} 
 
          While it is not really essential for the proof of Proposition
          \ref{prop:hyp1}, we will now find explicit holomorphic 
          co-ordinates -- ie, holomorphic
          functions -- on the Riemann
          surfaces $\Sigma_{q,t}$. These functions will also be crucial
          to the work in the later parts of the proof. The existence of
          the circle symmetry means that we are able to construct these
          by  elementary methods.
          
           Consider first the surface in $\bR^{3}$ defined by the equation
          $Q(\ux)=-1$,
         and take $x_{0}$ and $\theta$ as co-ordinates.
         We seek a holomorhic function $f$ on this surface of the form 
         $f=u(x_{0}) e^{i\theta}$. If we differentiate  Equation~(\ref{eq:cubic})
         we find that, with $Q$ fixed,
         \begin{equation}  \frac{\partial x_{0}}{\partial H}= p^{-4}. 
         \label{eq:dxdH}\end{equation}
         By Equation~(\ref{eq:conf}), the Cauchy-Riemann equations for $f$ on the surface are
         $$ \frac{\partial f}{\partial H} + i\,p^{-2} r^{-2}\, \frac{\partial f}{\partial
         \theta}=0, $$ so
   we see that $u(x_{0})$ must satisfy the
         equation
         \begin{equation} \frac{du}{dx_{0}} = \frac{p^{2}}{r^{2}}u=
         \frac{\sqrt{3x_{0}^{2} +1}}{\sqrt{2}(x_{0}^{2}+1)}
         u. \label{eq:CR1}\end{equation}
         We choose $u$ to be the solution of this equation with $u(0)=1$.
         Thus
         \begin{equation}    u(x_{0}) = \exp\left( \int_{0}^{x_{0}} \frac{\sqrt{3x^{2}+1}}
         {\sqrt{2}(x^{2} +1)} \ dx \right). \label{eq:int}\end{equation}
         We can evaluate this integral explicitly in terms of elementary
         functions,  but the formula that results is too cumbersome to be
         much use to us.
         Notice that $u(-x_{0}) = u(x_{0})^{-1}$.
         Clearly $u$ has the asymptotic behaviour
         \begin{equation}  u(x_{0}) \sim A x_{0}^{\nu} \label{eq:asymp}\end{equation}
         as $x_{0}\rightarrow +\infty$, where $\nu= \sqrt{3/2}$ and 
         \begin{equation}   A= \exp\left(\int_{0}^{\infty} \frac{ 
         \sqrt{3x^{2}+1}-\sqrt{3}\,x}{\sqrt{2}(x^{2}+1)}
         \, dx\right)={(2\sqrt{3})}^{\sqrt{3/2}}(\sqrt{3}-\sqrt{2}).
         \label{eq:constant}\end{equation}
         We  now define the function $F^{+}$ on the set $\{\ux: Q(\ux)<0\}$
         by
            \begin{equation}   F^{+}(\ux) = a^{\nu}
             u(\frac{x_{0}}{a})
            e^{i\theta}, \label{eq:F+}\end{equation}
            where $a=\sqrt{-Q}$. 
            The function $F^{+}$ is holomorphic on each quadric surface
            $Q(\ux)= -a^{2}$ for $a>0$, since scaling by $a^{-1}$ maps
            these conformally to the quadric $Q(\ux)=-1$.
            The asymptotic behaviour (\ref{eq:asymp}) implies that as 
            $\ux$ tends to the null cone  with
            $x_{0}$ fixed and positive $F^{+}(\ux)$ tends to $A x_{0}^{\nu}e^{i\theta}$,
            while if $x_{0}$ is fixed and negative $F^{+}(\ux)$ tends to
            zero on the null cone. We take these limiting values as the
            definition of $F^{+}$ on the null-cone. 
            Symmetrically, we define a function
            $F^{-}$ on $\{\ux:Q(\ux)<0\}$ by
            $$   F^{-}(\ux) = a^{\nu} u(-\frac{x_{0}}{a})
            e^{-i\theta}, $$
            so $F^{-}$ is also holomorphic on each surface, and $F^{+}F^{-}
            =a^{2\nu}=(- Q(\ux))^{\nu}$. The function $F^{-}$ now tends to zero on the
            part of the null cone where $x_{0}>0$.

            We follow a similar procedure on the set where $Q(\ux)>0$.
            On the sheet of the surface $\{Q(\ux)=1\}$ on which $x_{0}$ is
            positive we have a holomorphic function of the form $v(x_{0})e^{i\theta}$
            where, for $x_{0}>1$, the function $v$ satisfies
            $$ \frac{dv}{dx_{0}} = \frac{\sqrt{3x_{0}^{2}-1}}
            {\sqrt{2}(x_{0}^{2}-1)}
            v. $$
            This defines $v$ (with $v(1)=0$) up to a multiplicative constant,
            and we fix the constant by requiring that $v(x_{0})\sim
            A x_{0}^{\nu}$, where $A$ is given by Equation~(\ref{eq:constant}) above.
             Then we define $F^{+}$
            on $\{\ux: x_{0}>0,\ Q(\ux)>0\}$ by
            $$   F^{+}(\ux) = b^{\nu}v(\frac{x_{0}}{b})
            e^{i\theta},$$
           where $b=\sqrt{Q}$.  Symmetrically, we define $F^{-}(\ux)$
             to be $\overline{F^{+}(-\ux)}$ on the set $\{\ux: {x_{0}<0},\ Q(x_{0})>0\}$.
                                    
            To summarise, define open sets
            $$G^{+} = \{ \ux\in {\bR}^{3}: x_{0}>0 \ {\rm if}\  Q(\ux)\geq
            0\}, $$
            $$G^{-}= \{ \ux\in \bR^{3}: x_{0}<0 \ {\rm if}\  Q(\ux)\geq 0\}.
            $$
              Then we have:                       
            \begin{prop}
            The functions $F^{\pm}$ are smooth on $G^{\pm}$ and holomorphic
            on each connected component of the quadric surfaces $Q(\ux)=q$
            in $G^{\pm}$. \end{prop}
             The proof of this is a straightforward calculus argument involving
             the analytic continuation of the function $u(x_{0})$ to imaginary
             values of $x_{0}$.

 \subsection{Proof of Proposition \ref{prop:hyp1}}\label{ss:proofprop2}

    Let $\ur$ be a point in $\bR^{3}$ with $\vert \ur\vert =1$ and let
    $\Sigma$ be the quadric surface passing through $\ur$. We choose  a map
    $$L\co  D \times (-\textstyle{\frac{1}{4},\frac{1}{4}})\rightarrow \bR^{3},$$
    where $D$ is the unit disc in $\bC$,
    with the following properties.
    \begin{itemize}
    \item $L(0,0)= \ur$ and $z\mapsto L(z,0)$ gives a conformal parametrisation
    of a neighbourhood of $\ur$ in $\Sigma$.
    \item $H(L(z,q))$ and $ \theta(L(z,q))$ are independent of $q$
    \item $Q(L(z,q))= Q(\ur) + q$.
    \end{itemize}
    To construct this map we first choose a conformal parametrisation $L(z,0)$
    and then extend by integrating the vector field $\dQ$. This can all
    be done explicitly, using the conformal parametrisation by $F^{+}$
    above, but we do not need the detailed formulae; the crucial point for
       the proof of Proposition \ref{prop:hyp1} is the behaviour of the data under  
       scaling.
       The complex structure on the quadric surfaces pulls back to a leaf-wise
       structure on $D\times (-\frac{1}{4},\frac{1}{4})$ which is described by a matrix-valued
       function
       $     J(z,q)$. By construction $J(z,0)$ is the standard matrix $J_{0}$
       so
       $$    J(z,q)= J_{0} + q K(z,q)$$ say, with $K$ smooth.
       The pull-back by $L$ of the $2$--form $*(dQ\wedge dt)$ can be written as
       $$   A(z,q)\,i\,dz \wedge d\overline{z}, $$
       for some positive function $A$, with $A(z,q)\geq A_{0}>0$.
    As $\ur$ varies in the unit sphere we get a family of such maps and
    it is clear that, by compactness of the sphere,
     we can choose these so that $K$ and $A$ satisfy uniform 
     $C^{\infty}$--estimates on their derivatives, and $A_{0}$ is fixed
     independent of $\ur$.
     Having said this we will not complicate our notation by keeping the
     $\ur$--dependence explicitly.

     Now consider the point $\uR= \lambda \ur$ for some $\lambda\geq 1$.  Let
     $\psi_{0}$ be the value of the function $\psi$ at this point. We define
     a map $M(z,q,\tau)$ into $\bR^{4}$
     $$  M(z,q,\tau) = \Bigl(\lambda L\biggl( \frac{z}{\lambda^{3/2}}, 
     \frac{\psi_{0}}{\lambda^{2}} q\biggr), \psi_{0}^{-1} \tau\Bigr). $$
     The fourth condition of Lemma \ref{lem:defpsi}  implies that
      $\psi_{0}/\lambda^{2}=O(\epsilon)$, so we
     can suppose that $M$ is defined on $D\times I \times I$ for some fixed
     interval $I$. Then
     $$  M^{*}(\Omega)=dq\wedge d\tau + A\biggl(\frac{z}{\lambda^{3/2}}, 
     \frac{\psi_{0}}{\lambda^{2}} q\biggr)\, i\, dz\wedge d\overline{z} . $$
      Clearly, then, $M^{*}(\Omega)$ satisfies uniform $C^{\infty}$ bounds
      and with volume form bounded below by $A_{0}$ as the point $\uR=\lambda
      \ur$ ranges over the set $\{\vert \uR \vert\geq 1\}$. To prove
       Proposition \ref{prop:hyp1} we need to show that the almost-complex
      structure differs from the standard one in these co-ordinates
       by $O(\epsilon)$, with all derivatives. This almost-complex structure
       is given  by a matrix valued function which is the direct sum 
      \begin{equation}   J(\frac{z}{\lambda^{3/2}}) \oplus \left(
       \begin{array}{cc}
         0& -\Psi^2/\psi_{0}^2\\ \psi_{0}^2/\Psi^2 &0\end{array}
         \right) 
         \end{equation}
       where $\Psi$ is the composite $\psi \circ p \circ M$.
       
       Now the first term is
       $$ J(\frac{z}{\lambda^{3/2}})= J_{0} + \frac{\psi_{0} q}{\lambda^{2}}K(\frac{z}{\lambda^{3/2}},
       \frac{\psi_{0} q}{\lambda^{2}}). $$
       This satisfies the required bound since $\psi_{0} 
       \lambda^{-2}=O(\epsilon)$. Thus the real work involves the second
       term: we want to show that all  derivatives of 
       $1- \Psi/\psi_{0}$ are
         $O(\epsilon)$.  
        
        Return again to the function $L(z,q)$. Write
         $$    p(L(z,q))= G(z,q).$$
         Using homogeneity, our function $\Psi$ is given in the co-ordinates
         $M(z,q,\tau)$ by
         $$  \Psi(z,q,\tau)= \psi( \lambda^{1/2} G(\frac{z}{\lambda^{3/2}},
          \frac{\psi_{0} q}{\lambda^{2}})) . $$          
          We are left then with the elementary task of showing that the
          hypotheses in Lemma \ref{lem:defpsi} bound the derivatives of this composite
          function. For simplicity we will just work at the origin of the
          co-ordinates. We claim that 
          $$\lambda^{1/2} G(\frac{z}{\lambda^{3/2}},
           \frac{\psi_{0}}{\lambda^{2}}
          q)= \lambda^{1/2} G(0,0)+ \lambda^{-1} B(z,q), $$
          where $B$ is a smooth function, depending on the parameters
          $\lambda,\psi_{0}$ but all of whose derivatives are bounded.
          For if we write the Taylor series of $G$ in the schematic form
          $G(z,q)= \sum a_{IJ} z^{I} q^{J}, $ then
          $$  B(z,q)= \sum_{(I,J)\neq (0,0)} a_{IJ} \lambda^{3/2} 
          \lambda^{-3I/2 } \lambda^{-2J} \psi_{0}^{J} \ z^{I} q^{J}.$$    
          Now the assertion follows from the fact that 
          $\psi_{0}\leq C \lambda^{1/2}$. Thus our function is
          $$ 1-\Psi/\psi_{0} = 1- \psi(p_{0})^{-1} \psi(p_{0} + \lambda^{-1} B(z,w)), $$
          The fact that all derivatives of this are $O(\epsilon)$ follows
          from the condition
          $$ \psi^{(r)}\leq \epsilon c_{r} p^{2r} \psi, $$
          in Lemma \ref{lem:defpsi}.

       It is now straightforward to complete the proof of Proposition \ref{prop:hyp1}.
       We use the maps $M$ as above, together with their obvious translates
       in the $t$ variable, to get co-ordinate charts over a neighboorhood
       of $N\cap X_{0}$. Over the remainder of $X_{0}$ we can use the familiar
       rescaled osculating co-ordinates, just as in the case of compact
       symplectic manifolds.

       \section{Construction of approximately holomorphic sections}\label{sec:construct}

    We now start to work towards the verification of Hypothesis 2, involving
    sections of the line bundle ${\cal L}^{\otimes k}$ over $X$. The crucial constructions
    and arguments will take up this Section \ref{sec:construct} and the following Section
    \ref{sec:estimates}. As one would expect, the essential issues involve the local model
    around the zero set. Thus in Sections \ref{sec:construct} and \ref{sec:estimates} 
    we will work with a line
    bundle ${\cal L}$ over $\bR^{4}$ with a connection  of curvature
    $-i\Omega$. We use the almost-complex structure $J$, defined in the
    previous section, over the complement
    in $\bR^{4}$ of the $t$--axis. In Section~\ref{sec:completion} we will adapt our constructions
    to the $4$--manifold $X$. We write the line bundle ${\cal L}$ over $\bR^{4}$
    as the tensor product $${\cal L}={\cal L }_{1}\otimes {\cal L}_{2}$$
    where ${\cal L}_{1}$ has curvature $-i\,dH\wedge d\theta$ and ${\cal L}_{2}$
    has curvature $-i\, dQ\wedge dt$.

    We will omit some of the steps required to give a complete verification
    of Hypothesis 2. The proofs that we do give seem to us quite long enough,
    having in mind that the whole discussion is largely a matter of elementary
    calculus and geometry in $\bR^{3}$, and the techniques we develop can
    easily be extended to cover the parts we do not go through in detail.

\subsection{Holomorphic sections over the quadric surfaces}\label{ss:secl1}
    
    In this section we will work with the Hermitian line bundle  
${\cal L}_{1}$. We can  ignore the $t$--variable and consider ${\cal L}_{1}$ as a
line bundle over $\bR^{3}$. Our goal is to find 
 sections of ${\cal L}_{1}$ over suitable open sets in $\bR^{3}$
  which are {\it holomorphic}
 along the quadric surfaces and with appropriate localisation and smoothness
 properties. Exploiting the fact that the
 rotations in the $x_{1}, x_{2}$ plane act as symmetries of the whole set-up,
 we can find the desired sections by elementary methods. 
 
  Fix a trivialisation of ${\cal L}_{1}$ in which the connection form is
  $-iHd\theta$. We define  the section $\sigma$ of ${\cal L}_{1}$, in this trivialisation,
 to be
 \begin{equation}   \sigma = \exp( - \frac{ p^{6}}{18}) \label{eq:sigma}\end{equation}
 
 \begin{lem}
 The section $\sigma$ is holomorphic along each of the quadric surfaces
 in $\bR^{3}\setminus \{0\}$. \end{lem}

        With the connection form $-iH d\theta$, the Cauchy-Riemann equation
        for a holomorphic section $\sigma$ of ${\cal L}_{2}$ is:
        $$ p^2r^2 \frac{\partial \sigma}{\partial H} + i\biggl(
        \frac{\partial\sigma}{\partial\theta} - iH \sigma\biggr)=0.$$
        We have
        \begin{equation}  p^{4}= 6 x_{0}^{2} - 2Q
        \label{eq:implicitp}\end{equation}
        so,
        on a surface with $Q(\ux)$ constant,
        $$ 4p^{3} \frac{\partial p}{\partial H} = 12 x_{0} \frac{\partial
        x_{0}}{\partial H}= \frac{12x_{0}}{p^{4}}, $$
        using Equation~(\ref{eq:dxdH}). Thus
        \begin{equation}
         \frac{\partial p}{\partial H}= \frac{3 x_{0}}{ p^{7}},
         \label{eq:partialH}\end{equation}
        and the Cauchy-Riemann equation for a section with no $\theta$
        dependence is
        $$ \frac{3 x_{0}r^2}{ p^{5}}  \frac{\partial \sigma}{\partial p}=
         - H \sigma. $$
         But, since $H=x_{0}r^{2}$, this is just
         $$\frac{\partial \sigma}{\partial p}
= -\frac{p^{5}}{3} \sigma, $$
with solution $\sigma = \exp(-p^{6}/18)$.

         The section $\sigma$ can obviously be regarded as being localised
         at the origin in $\bR^{3}$, with exponential decay as we move
         away from the origin. We obtain more sections~--~holomorphic along
         the quadric surfaces~--~by multiplying $\sigma$ by suitable functions.
         The basic model to have in  mind here is that in ordinary flat
         space, say $\bC$. The Gaussian  $\exp(-\frac{ \vert z \vert^{2}}
         {4})$ represents a holomorphic section $s_{0}$ of the Hermitian line bundle
         with curvature $-i dx\wedge dy$ in a trivialisation in which the connection
         matrix is $ -\frac{i}{2}(xdy- ydx)$.  Given a point $a\in \bC$
         let $f_{a}$ be the holomorphic function
         $$   f_{a}(z)= \exp(\frac{\overline{a}z}{2}- \frac{\vert a\vert^{2}}{4}).$$
         Then $f_{a}s_{0}$ is a holomorphic section with norm
         $\exp\left(-\frac{\vert z-a\vert^{2}}{4}\right)$, concentrated around the
         point $a$ in $\bC$.

              To implement this idea in our setting, consider a section
              $\hat{\tau}=\exp(f)\,\sigma$ on one of the quadric surfaces,
              where $f=\mu+i\nu$ is a holomorphic function on the surface. We want to locate
              the points where $\vert \hat{\tau}\vert$ is stationary. In our
              trivialisation, these are points where the $H$ and $\theta$
              derivatives of $\mu + \log |\sigma|$ vanish. Since $\sigma$ is
              independent of $\theta$ and  $\frac{\partial
              \sigma}{\partial H}= -p^{-2}r^{-2} H \sigma$, the conditions
              are:
              $$  \frac{\partial \mu}{\partial H}= p^{-2} r^{-2} H,\ \ \frac{\partial
              \mu}{\partial \theta}=0. $$
             But the Cauchy-Riemann equations are
             $$  \frac{\partial \mu}{\partial H}= p^{-2} r^{-2} \frac{\partial
             \nu}{\partial \theta},\ \ \frac{\partial \nu}{\partial H}= - p^{2}
             r^{2} \frac{\partial \mu}{\partial \theta}, $$
             so the conditions just become:
             \begin{equation}  \frac{\partial f}{\partial \theta} =  i H\label{eq:hamil}\end{equation}

               Now,  given fixed $H_{0}, \theta_{0}$ we want to construct a section
              $\tau=\tau_{H_{0},\theta_{0}}$ of the line bundle ${\cal
              L}_{1}$ over a suitable open set in $\bR^{3}$ which,  on each quadric surface $Q(\ux)=q$, is holomorphic and
              which can be regarded as concentrated at the point in the
              surface with co-ordinates
              $H=H_{0}$, $\theta=\theta_{0}$. For simplicity we suppose $H_{0}\neq
              0$. The construction is simpler
              in the region where $Q(\ux)> 0$, and we begin with that case.
              We first assume $H_0>0$, in which case we consider
              the component where $x_{0}>0$. Here  we define
              \begin{equation}
                 \hat{\tau}=\hat{\tau}_{H_{0},\theta_{0}} = 
              \exp(\frac{H_{0}}{F^{+}(H_{0},\theta_{0},Q)} \ F^{+})\,\sigma. \label{eq:hattau}\end{equation}
              That is, we take the function $f$ above to be $A F^{+}$ where
              $A$ is, on each surface $Q(\ux)=q$,  the constant $H_{0}/F^{+}(H_{0},\theta_{0},q)$.
              Now $\frac{\partial F^{+}}{\partial \theta}= i F^{+}$ so
              $  \frac{\partial f}{\partial \theta}= iA F^{+}$ which, by
              construction, is
              equal to $iH$  when $H=H_{0}$, $\theta=\theta_{0}$. So the modulus of this section $\hat{\tau}$
              has a critical point at $(H_{0}, \theta_{0})$, which
              we will see is a maximum (cf \S \ref{ss:tau-est}). Now we normalise by defining
              $$  \tau_{H_{0},\theta_{0}} = \lambda \hat{\tau}_{H_{0},\theta_{0}},
              $$ where $\lambda= \vert \hat{\tau}(H_{0},\theta_{0}) \vert^{-1}$.
              Thus the value of $\vert \tau\vert $ at the point with co-ordinates
              $(H_{0},\theta_{0})$ on each quadric surface is $1$.
              
              If $H_0<0$, we work symmetrically on the region where $Q(\ux)> 0$ 
              but $x_{0}<0$ with the function $F^{-}$, setting
              $$  \hat{\tau} = \exp( - \frac{H_{0}}{F^{-}(H_{0},\theta_{0},Q)}
              F^{-})\,\sigma. $$
              The complication comes from the region $Q(\ux)<0$ where we
              need to use a combination of the functions $F^{\pm}$, smoothly
              interpolating between the two cases already defined.
              
              Consider the quadric surface $\{ Q(\ux)=-1\}$ on which we
              have functions $H$, $p$, and $u=u(x_{0})$ (defined
              by Equation (\ref{eq:int})). Any of $u, H, x_{0}$
              can (along with $\theta$) be used as a co-ordinate on the surface.
              For example we can regard $x_{0}$ as a function $x_{0}(H)$.
               Equation~(\ref{eq:asymp})
              implies that the positive function on this quadric
              $$ D=  (p^{2}-x_{0}) r^{2} u $$
              tends to infinity as $x_{0}\rightarrow \pm \infty$. Thus
              $D$ has a strictly positive minimum value, $\eta$ say. (The
              significance of this number will appear in the proof of
              Lemma \ref{lem:lemma6} in \S \ref{ss:tau-est}.)
              Now given small $\delta>0$ choose an even function $g$ on
              $\bR$
               with 
              \begin{itemize}
              \item $g'(h)\geq 0$ for $h\geq 0$
              \item $g(h) \geq \vert h\vert$
              \item
              $g(h)=\delta/2 $ for $\vert h\vert \leq \delta/4$ and $g(h)=\vert
              h \vert$ for $\vert h \vert \geq \delta$.
              
              \end{itemize}
              
               It is clear that
              if $\delta$ is sufficiently small we will have
              \begin{equation}  g(h)- h \leq \frac{\eta}{U(h)}
                \label{eq:smallbeta} \end{equation}
              for all $h>0$, where $U(h)= u(x_{0}(h))$. We fix such a $\delta$ and hence,
              once and for all, a function $g$.
              Define $\varphi(h)= \frac{1}{2}(h+ g(h))$ so 
              $$ \varphi(h) - \varphi(-h) = h $$
                   $$    \varphi(h) + \varphi(-h) = g(h),$$
                       and $\varphi(h)$ vanishes if $h<-\delta$.
       Now, on the set where $Q(\ux)<0$ write $Q(\ux) = -a^{2}$
              and define a section $\hat{\tau}_{H_{0}, \theta_{0}}$ of ${\cal
              L}_{1}$ by
              \begin{equation} \hat{\tau}_{H_{0}, \theta_{0}} = 
                \exp( \frac{\alpha}{F^{+}(H_{0},
              \theta_{0}, -a^{2})} F^{+} + \frac{\beta}{F^{-}(H_{0}, \theta_{0},
              -a^{2})} F^{-})\, \sigma, \label{eq:hattau2}\end{equation}
$$ \alpha = a^{3} \varphi(\frac{H_{0}}{a^{3}}),\leqno{\rm where} $$      
      $$         \beta = a^{3} \varphi(-\frac{H_{0}}{a^{3}}).  $$
               $$ \alpha-\beta= H_{0} ,\ \ \alpha+\beta= a^{3}g(\frac{H_{0}}{a^{3}}).\leqno{\rm Thus}$$        On each quadric surface $Q(\ux)=-a^{2}$ the section 
              $\hat{\tau}=\hat{\tau}_{H_{0},\theta_{0}}$ is holomorphic,
 since $\alpha,\beta$ and $F^{\pm}(H_{0}, \theta_{0}, -a^{2})$ are all
 constant on the surface. We claim that, on each surface, $\vert \hat{\tau}\vert$
  is stationary at the point where $H=H_{0}$ and $\theta=\theta_{0}$.
  Indeed $\hat{\tau}= e^{f} \sigma$ where $ f= A F^{+} + B F^{-}$ and $A,B$
  are constants on the surface. So
  $$  \frac{\partial f}{\partial \theta} = iA F^{+}- i BF^{-}$$
  which is equal to $i(\alpha-\beta)$ at the given point. Then the claim
  follows from the fact that $\alpha-\beta = H_{0}$. 
  Once again, we define $\tau_{H_{0},\theta_{0}}$ by normalising so that
  the modulus is $1$ at the critical point.
  
  To sum up, if $H_{0}> 0$ we have defined sections $\tau_{H_{0},\theta_{0}}$
  separately over the two regions $\{ Q(\ux)<0\}$ and $\{Q(\ux)>0,\ x_{0}>0\}$.
  However it follows from the construction that  these sections have the
  same limit over the positive part of the null cone, and define a smooth
  section over the region $G^{+}\subset \bR^{3}$. This is because the coefficient
  $B$ of $F^{-}$ vanishes near the positive part of the null cone. Likewise
  if $H_{0}< 0$ we get a section $\tau_{H_{0},\theta_{0}}$ defined over
  $G^{-}$. We obtain the following:  
  \begin{prop}
  For any $H_{0}\neq 0$, $\theta_{0}$ the section $\tau_{H_{0},\theta_{0}}$ defined
   above is a smooth section of ${\cal L}_{1}$ over $G^{+}$ or $G^{-}$.
   The section is holomorphic along each connected component of the quadric
   surfaces in its domain of definition and has modulus $1$ at the points
   with co-ordinates $(H_{0},\theta_{0})$.
   \end{prop}
   
   Note that some of the steps in the construction work equally well when
   $H_{0}=0$ but there are some  difficulties. From one point of view
   this is because we are really attempting to  define a family of sections indexed
   by the set of integral curves of the vector field $\frac{\partial}{\partial
   Q}$ on $\bR^{3}\setminus \{0\}$ and this set, in its natural topology,
   is not Hausdorff. To avoid these essentially irrelevant complications we do
   not define sections $\tau_{H_{0}, \theta_{0}}$ when $H_{0}=0$.

\subsection{Sections of ${\cal L}_{2}$ and cut-off functions}\label{ss:secl2}

In this subsection we first define suitable sections of the line bundle
${\cal L}_{2}$ over $\bR^{4}$. Recall that this has curvature $-i dQ\wedge dt$.
let $(\ux', t')$ be a point of $\bR^{4}$ where $\ux'$ has $(Q,H,\theta)$ 
co-ordinates $(Q_{0}, H_{0},
\theta_{0})$. Let $\psi_{0}$ be the  value of the function $\psi$ at $\ux'$.
We can choose a trivialisation of the
bundle such that the connection form is 
$$ -\frac{i}{2}( (Q-Q_{0}) dt - (t-t')
dQ).$$
In this trivialisation, we define a section by  
\begin{equation}  \hat{\rho}_{\ux', t'} = \exp\left( - \frac{ \psi_{0}^{2}(t-t')^{2} +
 \psi_{0}^{-2} (Q-Q_{0})^{2}}{4}
\right).\label{eq:rhodef}\end{equation}
(The trivialisation is ambiguous up to an overall phase, so this definition
is not strictly precise, but we can ignore this here.)  Notice that in a region where $\psi$
is constant the section will be a holomorphic section of ${\cal L}_{2}$;
we postpone until Section \ref{sec:estimates} the estimates for $\db \hat{\rho}$ in general.
Obviously $\vert \hat{\rho}\vert $ achieves its maximum value $1$ at points where
$Q=Q_{0}$, $t=t'$.

The section $\hat{\rho}_{\ux',t'}$
 decays rapidly away from the surface $Q(\ux)=Q_{0}$. We will now introduce
 a cut-off function to construct a section which vanishes outside a neighbourhood
 of this surface. Let $\chi(q)$ be a fixed, standard, cut-off function
 equal to $1$ for $\vert q\vert\leq 1$ and vanishing when $\vert q\vert\geq
 2
 $. Let $b_{1}$ be a small positive constant, to be fixed later,
 and define a function $\chi_{Q_{0}}$ on $\bR^{3}$ by
 \begin{equation}     
 \chi_{Q_0}= \chi\left( \frac{\epsilon}{b_{1}} \frac{Q-Q_{0}}{\psi_{0}} \right)
 .\label{eq:Qcutoff}\end{equation}
 
 Then set
 \begin{equation}    \rho_{\ux',t'} = \chi_{Q_{0}} \hat{\rho}_{\ux',t'}
 \label{eq:rho}\end{equation}
  We now return to the sections $\tau_{H_{0},\theta_{0}}$ defined in the
 previous section. We want to modify these by suitable cut-off functions
 to overcome the difficulties with their domains of definition. This cut-off
 construction will depend on another small positive parameter $b_{2}$. Let $c_{0}$
 be the constant from Lemma \ref{lem:defpsi}, so $\psi(p)\leq c_{0} \epsilon p^{4}$ for
 $\vert p\vert \geq 1$. Recall
 that $p^{4}$ is the quadratic form $4x_{0}^{2}+ r^{2}$ on $\bR^{3}$. We  choose the constant
 $b_{2}$ so that $b_{2}c_{0}<\frac{1}{10}$, say. Then the quadratic form $Q+b_{2}c_{0}p^{4}$
 is indefinite. Define $G_{0}^{\pm}\subset \bR^{3}$ by
 $$   G_{0}^{\pm}= \{ \ux: Q(\ux)+b_{2} c_{0}p^{4}<0 \ \ {\rm or}\ \ 
  Q(\ux) + b_{2}c_{0}p^{4}
 \geq 0 \ \ {\rm and}\ \ \pm x_{0}>0\}. $$
 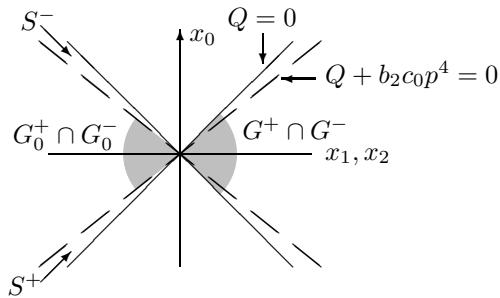
\begin{figure}[b]
 \setlength{\unitlength}{5mm}
 \centering
 \begin{picture}(8,7.5)(-4,-3.5)
 \psset{unit=\unitlength}
 \pswedge[fillstyle=solid,fillcolor=lightgray,linewidth=0,
          linecolor=lightgray](0,0){1.5}{-45}{45}
 \pswedge[fillstyle=solid,fillcolor=lightgray,linewidth=0,
          linecolor=lightgray](0,0){1.5}{141.5}{218.5}
 \put(-3.5,0){\line(1,0){7}}
 \put(0,-3){\vector(0,1){6.3}}
 \put(3.8,-0.15){\small $x_1,x_2$}
 \put(0.2,3){\small $x_0$}
 \put(-3,-3){\line(1,1){6}}
 \put(3,-3){\line(-1,1){6}}
 \multiput(3.75,-3)(-1.125,0.9){7}{\line(-5,4){0.75}}
 \multiput(-3.75,-3)(1.125,0.9){7}{\line(5,4){0.75}}
 \put(1.6,0.4){\small $G^+\cap G^-$}
 \put(-4.5,0.3){\small $G_0^+\cap G_0^-$}
 \put(1.2,3.4){\small $Q=0$}
 \put(2.2,3.2){\vector(0,-1){0.8}}
 \put(3.6,2){\vector(-1,0){0.9}}
 \put(3.8,1.9){\small $Q+b_2c_0p^4=0$}
 \put(-3.7,3.4){\vector(1,-1){0.8}}
 \put(-4.3,3.3){\small $S^-$}
 \put(-3.7,-3.4){\vector(1,1){0.8}}
 \put(-4.65,-3.75){\small $S^+$}
 \end{picture}
 \caption{The sets $G^\pm$ and $G_0^\pm$}
 \end{figure}
 Thus $G_{0}^{\pm}\subset G^{\pm}$ and $G_{0}^{+}\cup G_{0}^{-} = \bR^{3}
 \setminus \{0\}$. Let $N$ be the $1$--neighbourhood, in the metric $g$,
  of the plane-minus-disc $\{ (0,x_{1},x_{2})\in
 \bR^{3}: x_{1}^{2}+ x_{2}^{2} \geq 4\}$.
  It is easy to check (using the fact that $Q=-\frac{1}{2}p^4$ for $x_0=0$,
 the formula $\|\frac{\partial}{\partial Q}\|_g=
 \psi^{-1}$ and the estimates on $\psi$) that we can choose $b_{2}$ small
 enough (depending on the constants in Lemma \ref{lem:defpsi}) so that 
 \begin{equation} 
 G_{0}^{+} \cap G_{0}^{-} \supset N\label{eq:overlap}. \end{equation}
 We now fix a value of $b_{2}$ such that (\ref{eq:overlap}) holds. 
 Let $\lambda$ be  a standard cut-off function with 
 $\lambda(q)=1 $ for $q\leq -1$ and $\lambda(q)=0$ for $q>-1/2$. Define
a function  $L$ on the set $\{ \vert \ux \vert\geq 1\}$ in $\bR^{3}$ by
$$    L= \lambda( \frac{\epsilon}{b_{2}} \frac{Q}{\psi}). $$
Suppose a point $\ux$ lies in the support of $\nabla L$. Then we must have  
 $$    -\frac{b_2}{\epsilon}\psi < Q<0. $$ 
 Thus $ - b_2c_{0} p^{4}< Q<0$. So the support of $\nabla L$ is contained
 in the set
 $$ \{ \ux: \vert \ux\vert>1,\ Q<0,\ Q + b_2c_{0} p^{4}>0\}$$
 which is the disjoint union of two components,
 $$  S^{\pm}= (G^{\pm}\setminus G^{\pm}_{0})\cap\{\vert \ux\vert>1\}. $$
 It follows that there are smooth functions $\hat{L}^{+}, \hat{L}^{-}$ on $\{\vert
 \ux\vert>1\}$, supported in $G^{+}$, $G^{-}$ respectively  and equal to
 $1$ on $G^{+}_{0}$, $G^{-}_{0}$ respectively, such that $\hat{L}^{\pm}$ and $L$
 have the same restriction to $S^{\pm}$. Finally, define
 $$   L^{\pm} = \chi( \frac{2}{\vert \ux\vert})\  \hat{L}^{\pm}. $$ 
 Now suppose that the point $\ux' \in \bR^{3}$ with co-ordinates
  $Q_{0}, H_{0}, \theta_{0}$ has $\vert \ux'\vert >3. $
  Suppose that $H_{0}\neq
 0$. We
 define a section $\tau^{*}_{\ux'}$ of ${\cal L}_{1}$ as follows. 
  If $H_{0}>0$ the section $\tau_{H_{0},\theta_{0}}$ is smooth on
  $G^{+}$ and we set
  $$   \tau^{*}_{\ux'}= L^{+} \tau_{H_{0}, \theta_{0}}, $$
  extending in the obvious way by zero outside the support of $L^{+}$.
  Thus 
  $\tau^{*}_{\ux'}$ is equal to the section $\tau_{H_{0},\theta_{0}}$~--~holomorphic
  along the quadric surfaces~--~near $\ux'$, and the modulus 
  of $\tau^{*}_{\ux'}$ at the point $\ux'$ is $1$. In fact, because
  of (\ref{eq:overlap}), the $1$--ball $B_{\ux'}$ 
  centred at $\ux'$ in the metric $g$ is contained in $G_0^+$, and by estimating the norm
  of $d(|\ux|^2)$ one can verify that $B_{\ux'}\subset\{|\ux|>2\}$. Hence
  $\tau^{*}_{\ux'}$ is equal to $\tau_{H_{0},\theta_{0}}$ on the unit ball
  $B_{\ux'}$. 

  We proceed similarly if $H_0<0$ (using $L^-$ instead of $L^+$).
  Finally, we combine this with the other construction. For 
  $(\ux', t')$ as above, we set
  \begin{equation}  s_{\ux', t'}=   \tau^{*}_{\ux'} \otimes \rho_{\ux',t'}.
  \label{eq:sdef} \end{equation}
What we have now achieved is a collection of sections of the line bundle ${\cal
L}$ and in the next section we will derive the estimates which will ultimately allow us to
verify Hypothesis 2. That hypothesis requires rather more input data. Associated
to each point $(\ux',t')$ we need not just one section $s_{\ux',t'}$
 of ${\cal L}$ but a triple of
sections $(s, s', s'')$ say, so that $s'/s$ and $s''/s$ give local approximately
holomorphic co-ordinates. We will not go through this part of the construction
in detail, since it would not contain any new ideas. For example, one approach
is to define $s', s''$ by differentiating the section $s_{\ux', t'}$
with respect to the parameters $\ux', t'$.

\noindent {\bf Remark}\qua To illustrate this possible approach to the
construction of $s'$ and $s''$, consider the much simpler example of
a line bundle with curvature $-i\,dx\wedge dy$ over $\mathbb{C}$.
Given any $z'=x'+iy'\in\mathbb{C}$, this line bundle admits a holomorphic
section $s_{x',y'}$ with $|s_{x',y'}(z)|=\exp(-\frac{1}{4}|z-z'|^2)$; in
a trivialisation where the connection is $\nabla=d+\frac{i}{2}(y\,dx-x\,dy)$,
such a section is given eg, by $s_{x',y'}(z)=\exp(-\frac{1}{4}|z|^2+
\frac{1}{2}\overline{z'}\,z-\frac{1}{4}|z'|^2)$. Differentiating with
respect to $x'$ and $y'$, one easily checks that $(\frac{\partial}{\partial
x'}-\frac{i}{2}y')\,s_{x',y'}=\frac{1}{2}(z-z')\,s_{x',y'}$ and 
$(\frac{\partial}{\partial
y'}+\frac{i}{2}x')\,s_{x',y'}=-\frac{i}{2}(z-z')\,s_{x',y'}$. The ratio
of either one of these sections to $s_{x',y'}$ defines a local
holomorphic co-ordinate near $z'$.
In our case the sections $s_{\ux',t'}$ depend on four real parameters;
as in the example, the sections obtained by differentiating in directions
which belong to a same complex line are essentially redundant, and we
are left with two sections $s'$ and $s''$ whose ratio to $s_{\ux',t'}$
gives local approximately holomorphic co-ordinates near $(\ux',t')$.

\section{Estimates for approximately holomorphic sections}\label{sec:estimates}

\subsection{Estimates for $\tau$} \label{ss:tau-est}

In this subsection (\ref{ss:tau-est}) and the following (\ref{ss:dtau-est}) we develop estimates for
the sections constructed in \S \ref{ss:secl2}. Fix $H_{0}$ and $\theta_{0}$ and suppose
that $H_{0}>0$ (of course there will be symmetrical statements for the
case $H_{0}<0$). Then we have defined a section $\tau=\tau_{H_{0},\theta_{0}}$
of ${\cal L}_{1}$ over the open set $G^{+}\subset \bR^{3}$. We introduce
some notation. Let $\ux$ be a point in $G^{+}$, with co-ordinates $(Q,H,\theta)$.
Let $\ux'$ be the point in $G^{+}$ with co-ordinates $(Q,H_{0},\theta_{0})$.
Let $\ux''$ be the point with co-ordinates $(Q,H_0,\theta)$ if $\ux$ does not
lie on the positive $x_{0}$--axis, and otherwise set
$\ux''=\ux'$. We define two functions $S= S_{H_{0},\theta_{0}}$ and $L=L_{H_{0},\theta_{0}}$
on $G^{+}$. The value $S(\ux)$ is the distance in the metric $g$ from
$\ux$ to $\ux''$, measured along the quadric surface through $\ux$.  The value $L(\ux)$ is $1/2\pi$ times the length, in
the metric $g$, of the orbit of $\ux'$ under the rotation action. Now for
$\alpha>0$ we define a function $E_{\alpha}=E_{\alpha,H_{0},\theta_{0}}$
on $G^{+}$ by
\begin{equation}  E_{\alpha}(\ux)= \exp(-\alpha\left( S(\ux)^{2} + (\theta-\theta_{0})^{2}
L(\ux)^{2}\right)\,)\label{eq:decayfn}\end{equation}
 (Here we interpret $(\theta-\theta_{0})$ as taking values in $(-\pi,\pi]$;
 thus $L(\ux)(\theta-\theta_{0})$ is the distance in the metric $g$ from
 $\ux'$ to $\ux''$, measured along the circle orbit.)

 Now given $c>0$ let $\Omega^{+}_{c}$ be the set
 \begin{equation} \Omega^{+}_{c}= \{ \ux: \ux\in G^{+},\ \vert \ux\vert
 >1,\ Q(\ux)< -c\ \  {\rm if}\ \  x_{0}<0\}. \label{eq:omegadef} \end{equation}
 
 The  result we will prove in  this section
is the following:
\begin{prop}\label{prop:tau-est}
For any $c$ there are $C,\alpha$ (independent of $H_{0}, \theta_{0}$)
such that in $\Omega_{c}^{+}$,
$$  \vert \tau \vert \leq C E_{\alpha}. $$ \end{prop}

Recall that, given $\ux$ and $H_{0}, \theta_{0}$, we write $\ux'$ for the point
with co-ordinates $H_{0}, \theta_{0}$ on  the quadric surface through $\ux$.
In \S \ref{ss:dtau-est} below we will prove:
\begin{prop} \label{prop:dtau-est}
For any $c, r$ there are $C,\alpha$ such that at points
 $\ux\in \Omega_{c}^{+}$ for which $\vert \ux'\vert\geq 1$ and for all 
$p\leq r$:
\begin{itemize}
\item $\vert \nabla^{p} \tau \vert \leq C E_{\alpha}$,
\item $\vert \nabla^{p} \db \tau \vert \leq \epsilon\, C E_{\alpha}
$ at points where $\vert \ux\vert\geq 3$. \end{itemize}  \end{prop}
 
 Here, more precisely, $\db\tau$ is  defined by extending the section $\tau$
 to $G^{+}\times \bR$ but since there is no $t$ dependence we can formulate
 the result entirely within $\bR^{3}$.

 We begin the proof of Proposition \ref{prop:tau-est} by considering the restriction of
 $\tau$ to the sheet $\{ x_{0}>0\}$ of the quadric $Q(\ux)=1$.
 We may obviously suppose that $\theta_{0}=0$ and to begin with we consider
 the restriction to $\theta=0$. Thus we are considering the section
 $\tau$ over a single arc, homeomorphic to $[0,\infty)$. In our analysis we
 will use  two
 convenient co-ordinates on this arc. One co-ordinate is the function $v$,
 the modulus of the holomorphic function $F^{+}$. The other co-ordinate
 is the arc length $s$, measured from the intersection with the $x_{0}$--axis,
 in the metric $g$. We write $v_{0}, s_{0}$ for the co-ordinate values corresponding
 to $H=H_{0}$; ie, corresponding to the point $\ux'$. The co-ordinates
 $v$ and $s$ both run from $0$ to $\infty$ and the asymptotic relation
 between them is
      $$ s\sim C v^{\sqrt{3/2}}, $$ 
      as $s, v\rightarrow \infty$. In fact,
       in terms of the radial co-ordinate $r$, we have
       $$ s\sim C' r^{3/2}\ , \ v \sim C'' r^{\sqrt{3/2}}. $$
       The corresponding asymptotic relations hold for the mutual derivatives
       of these different co-ordinate functions.  
       
        Recall that our basic section is $\sigma=\exp(-p^{6}/18)$.
 We can write $p^{6}/18$ as a function of $v$, $f(v)$ say, 
 on this arc. Thus $f$ is an increasing function of $v$, asymptotic
 to a multiple of $v^{\lambda}$, where $\lambda=\sqrt{6}>2$. We introduce
 a piece of notation. For a function $g$ of a real variable $v\in [0,\infty)
 $ we write
 \begin{equation}  \Delta_{g}(v,v_{0}) = g(v)-g(v_{0}) - (v-v_{0})
 g'(v_{0}). \label{eq:Deltadef}\end{equation}
 The relevance of this, working with the co-ordinate $v$ over the arc,
 is that
  our definition of the section $\tau$ 
   is just
 $$\tau(v) = \exp\left(-\Delta_{f}(v,v_{0})\right).$$
 To see this, note that $\frac{d p}{d x_0}=3x_0p^{-3}$ 
 (by differentiating $p^4=6x_0^2-2Q$ with $Q$ fixed), and
 $\frac{d v}{d x_0}=p^2r^{-2}\,v$ (by definition of $v$). Therefore
 $f'(v)=H/v$, and $\Delta_f(v,v_0)=\frac{1}{18}(p^6-p_0^6)-\frac{H_0}{v_0}
 (v-v_0)=-\log |\tau|$.
 The first point to note is:
 \begin{lem}
 The function $f$ is a convex function of $v$, its second derivative
 is strictly positive.\end{lem}
 To see this, recall that $f'(v)=H/v$, so we have
 to show that $H/v$ is an increasing function of $v$, or equivalently
 of the variable $x_{0}$. Now, with $Q$ fixed,
 $$ \frac{dH}{dx_{0}}= 6x_{0}^{2}-2= p^{4}\ , \ \frac{dv}{dx_{0}} = \frac{p^{2}
 x_{0}}{H} v. $$
 Thus $$\frac{d}{dx_{0}}(H/v)= \frac{p^{4}}{v}- \frac{H}{v^{2}}
 \frac{p^{2}x_{0}}{H} v = \frac{p^{2}}{v} (p^{2}- x_{0}). $$
 This is positive since $p^{2} = \sqrt{4x_{0}^{2} + r^{2}} > x_{0}$.

 This Lemma shows that the modulus of the section $\tau$ does indeed attain
 a unique maximum at the point $v=v_{0}$. Next we need:
 \begin{lem}\label{lem:deltabound}
 Suppose $f$ is a function of $v\in [0,\infty)$ and $f''(v)\geq k(1+ v^{\lambda-2})$
 for some $k>0$, $\lambda>2$. Then there is a constant $c$ such that
 $$ \Delta_{f}(v,v_{0}) \geq c 
 \left((1+ v)^{\lambda/2} - (1+v_{0})^{\lambda/2}\right)^{2}.$$
 \end{lem} 
 
To prove this Lemma note that we can write
\begin{equation}\Delta_{f}(v_{1},v_{0}) = \int_{v_{0}}^{v_{1}} f''(v)(v_{1}-v)
\,dv. \label{eq:integralform}\end{equation}
Thus the hypothesis implies that $\Delta_{f}(v,v_{0})\geq \Delta_{g}(v,v_{0})$
where $g(v)= k(\frac{v^{2}}{2} + \frac{v^{\lambda}}{\lambda(\lambda-1)})$.
So, for a suitable constant $c$,
$$\Delta_{f}(v,v_{0}) \geq c\left( (v^{\lambda}- v_{0}^{\lambda})
-\lambda (v-v_{0}) v_{0}^{\lambda-1} + (v-v_{0})^{2}\right). $$
The convexity of the function $v^{\lambda}$ implies that the expression
$$v^{\lambda} - v_{0}^{\lambda} - \lambda(v-v_{0}) v_{0}^{\lambda-1}
$$ is non-negative. By considering the scaling behaviour
under simultaneous scaling of $v$ and $v_{0}$ (or by using the Taylor
formula), one sees that it is bounded below
by a positive multiple of  $(v-v_{0})^{2} (v^{\lambda-2} +
v_{0}^{\lambda-2})$. Thus
$$ \Delta_{f}(v,v_{0}) \geq c (v-v_{0})^{2} (1+ v^{\lambda-2}
+ v_0^{\lambda-2}) . $$
On the other hand it is clear that
\begin{eqnarray*}
\left( (1+v_{0})^{\lambda/2} - (1+ v)^{\lambda/2}\right)^{2}
&\leq& c'(v-v_{0})^{2} \left( (1+ v)^{\lambda/2-1} + (1+v_{0})^{\lambda/2-1}\right)^{2}
\\ &\leq& c''(v-v_{0})^{2} \left( 1+v^{\lambda-2}+v_{0}^{\lambda-2}\right)
\end{eqnarray*}
which implies the desired result.

In the case of $f(v)=\frac{1}{18}p^6$ the quantity $f''(v)=\frac{d}{dv}(
\frac{H}{v})=\frac{r^2}{v^2}(p^2-x_0)$ is bounded below by a positive
multiple of $(1+v^{\lambda-2})$, so by Lemma \ref{lem:deltabound} we
obtain a bound on $|\tau|$.
Now the functions $s$ and $\tilde{s}= (1+v)^{\lambda/2}$ have the same
asymptotic behaviour, so the derivative $\frac{ds}{d\tilde{s}}$ is bounded above and below
by positive constants. Thus we see that on this arc
$$\vert\tau(s)\vert \leq \exp(-c(s-s_{0})^{2}),$$
which is precisely the statement of Proposition 5 (on the arc).

Still working on the surface $Q(\ux)=1$, $x_{0}>0$, we now consider the dependence on the
 angular variable $\theta$. (Recall that we are assuming $\theta_{0}=0$.)
 We have 
 $$\log \vert \tau(s,\theta)\vert= \log \vert \tau(s,0)\vert - 
 \frac{v}{v_{0}} H_{0} (1-\cos\theta). $$
 Now $H_{0}$ is bounded below by a multiple of $s_{0}^{2}$. For $\theta\in
 [-\pi,\pi]$, the function
 $(1-\cos\theta)$ is bounded below by a multiple of $\theta^{2}$. Thus
 we have
 \begin{equation}   \vert \tau(s,\theta)\vert \leq \exp(-c\left( (s-s_{0})^{2} + 
 \frac{v}{v_{0}}s_{0}^{2} \theta^{2}\right) ). \label{eq:bound}\end{equation}
   Recall that
 we defined $L=L(\ux')$ to be $1/2\pi$ times the length of the circle orbit through
 $\ux'$. This is $p\,r$, evaluated at $\ux'$, which is bounded
 above and below by multiples of $s_{0}$. Thus, on this quadric surface, we have
 $$ E_{\alpha}(s,\theta) = \exp(-\alpha\left( (s-s_{0})^{2}+ C s_{0}^{2}\theta^{2}\right)).
 $$
 The difficulty comes from the term $\frac{v}{v_{0}}$ in  
 Equation~(\ref{eq:bound}). For this we use:
 \begin{lem}\label{lem:voverv0}
 There is a constant $C>0$ such that
 $$  (s-s_{0})^{2} +  s_{0}^{2}  \theta^{2} \leq C ( (s-s_{0})^{2}
 + \frac{v}{v_{0}} s_{0}^{2} \theta^2)$$
 for all $s,s_{0}\geq 0$ and $\theta\in [-\pi,\pi]$.
 \end{lem} 
 We consider the function $v/s$ as a function of $s$. This tends to a positive
 limit as $s\rightarrow 0$ and tends to zero as $s\rightarrow \infty$.
 Thus there is a constant $b$, independent of $s$ and $s_0$, such that 
 $$ \frac{v_{0}}{s_{0}}\leq b \frac{v}{s} $$
 whenever $v<v_{0}$. This means that whenever $v/v_{0}\leq 1/2b$ we have
  $ s\leq s_{0}/2$.  So either $v/v_{0}>1/2b$ in which case the desired
  inequality holds with $C=2b$, or $s_{0}^{2}\leq 4 (s-s_{0})^{2}$ in which
  case the inequality holds with $C=4\pi^{2}+1$.

 Lemma \ref{lem:voverv0} implies that $\tau$ is bounded by a suitable function $E_{\alpha}$
 on the quadric surface $Q(\ux)=1$, $x_{0}>0$. We can extend this bound to
 the entire cone $Q(\ux)>0$, $x_{0}>0$ in a very simple way, by homogeneity.
 We will use this principle repeatedly below, so we will spell it out
 clearly now.
 If we write, in the fixed trivialisation of the line bundle ${\cal L}_{1}$
 \begin{equation}   
 \tau=\tau_{H_{0}, \theta_{0}}= \exp(-A(\ux; H_{0},\theta_{0})) \label{eq:twovar}\end{equation}
 then the function $A$ satisfies
 \begin{equation} A(\lambda \ux; \lambda^{3} H_{0}, \theta_{0}) = \lambda^{3} A(\ux;
 H_{0}, \theta_{0}). \label{eq:homog}\end{equation}  
 The functions $\log E_{\alpha}$ satisfy exactly the same scaling behaviour
 $$\log E_{\alpha,\lambda^{3} H_{0},\theta_{0}}(\lambda \ux) = \lambda^{3}
 \log E_{\alpha, H_{0}, \theta_{0}} (\ux). $$
 Thus the bound $\vert \tau\vert \leq E_{\alpha}$ on the quadric surface,
 for all choices of the parameter $H_{0}$, immediately gives the same bound
 over the whole cone. 
 
 This scaling behaviour may be clearer if we change notation and regard
 $A$ and $E_{\alpha}$  as  functions of pairs of points 
 $\ux,\ux'$ on the same quadric surface.
 Then the scaling reads
 $$  A(\lambda \ux, \lambda \ux') = \lambda^{3} A(\ux,\ux')\ , \ \log E_{\alpha}(\lambda
 \ux,\lambda \ux') = \lambda^{3} \log E_{\alpha}(\ux,\ux'). $$
We now follow a similar argument for the region where $Q(\ux)<0$. By the
same scaling argument it suffices to work on the quadric $Q(\ux)=-1$.
Again, we begin with arc on this quadric where $\theta=0$. We have two
different co-ordinates on this arc. One is the arc length $s$ in the metric
$g$  which now runs
from $-\infty$ to $\infty$. The other is the function $u$, the modulus
of $F^{+}$, which runs over $(0,\infty)$. The function $u$ is asymptotic
to $|s|^{\pm \sqrt{2/3}}$ as $s\rightarrow \pm \infty$. The choice of the
parameter $H_{0}>0$ defines corresponding values $u_{0}>1$, $s_{0}>0$. 
Recall that over this arc our section $\tau$ is given by
$$   \tau= \exp(- f(u) +  
\alpha \frac{u}{u_{0}} + \beta \frac{u_{0}}{u}+c(u_0)) $$
where now $f$ is  $p^{6}/18$, expressed as a function of $u$ on the quadric,
$c(u_0)$ is a normalisation constant ensuring that $|\tau(u_0)|=1$,
and $\alpha$ and $\beta$ are defined by $u_{0}$ as in \S \ref{ss:secl1}.  
\begin{lem}  \label{lem:lemma6}
$\log \vert \tau(u) \vert 
$ has just one critical point, when $u=u_{0}$ and this point is a
global maximum.   \end{lem}

To see this, we have
$$  \frac{d}{du} (-\log \vert \tau\vert) = \frac{H}{u} - \frac{\alpha}{u_{0}}
+ \beta \frac{u_{0}}{u^{2}}. $$
We want to see that this vanishes only when $u=u_{0}$, where it vanishes
by construction (since $\alpha-\beta=H_{0}$). Thus it suffices to show
that the function $\frac{H}{u} +\beta \frac{u_{0}}{u^{2}}$ is an
increasing function of $u$, or equivalently of $H$. 
Now
$$  \frac{d}{dH}\left( \frac{H}{u} +
 \beta \frac{u_{0}}{u^{2}}\right)=
    \frac{1}{u} - \frac{H}{p^{2} r^{2} u} - 2\beta \frac{u_{0}}{p^{2}
    r^{2} u^{2}}, $$
    using the fact that $\frac{du}{dH}= \frac{u}{p^{2} r^{2}}$.   
 Rearranging terms, we need $$2\beta u_{0} < (p^{2}-x_0) r^{2}u.$$
 But this is precisely the condition we required in the choice of $\alpha,
 \beta$ defining~$\tau$ (see Equation (\ref{eq:smallbeta}), and recall
 that $\beta=\varphi(-H_0)=\frac{1}{2}(g(H_0)-H_0)\le \frac{\eta}{2u_0}$, where
 $\eta=\min\{(p^2-x_0)r^2u\}$), so the assertion follows. 
   This discussion also shows that $\log \vert \tau\vert$ is a concave function
   of $u$ along the arc $\theta=0$. Thus $u_{0}$ is a maximum along this
   arc. On the other hand the $\theta$--dependence is again proportional
   to $\cos\theta$ so clearly the maximum on each circle $u={\rm constant}$
   is attained when $\theta=0$.

We claim that, on the arc $\theta=0$, 
$$   \tau \leq \exp (-\alpha (s-s_{0})^{2}). $$ 
The proof  follows the same pattern as in the positive case above. Recall
that $\beta=0$ once $u_{0}$ is bigger than some $K>1$ say. 
When $u,u_{0}>K$ the argument is identical. There are then various other
cases to check, a task which we will largely leave to the reader. We just
discuss  
two representative sample cases. First, if $u_{0}=1$ then we have to show that
$$  f(u)-f(1) - c(u+u^{-1} -2) \geq c s^{2}. $$
This holds when $u$ is close to $1$ by the critical point analysis above.
When $u\rightarrow \infty$ the left hand side grows like $f(u)\sim
u^{\sqrt{6}}$ since
$\sqrt{6} >1$, and this is the same growth as $s^{2}$. Similarly when
$u\rightarrow 0$. For the second case, consider $u\rightarrow 0$ and
  $u_{0}\rightarrow \infty$. Then we have to show that
  $$   f(u)-f(u_{0}) - (u-u_{0}) f'(u_{0}) \geq c (s_{0}^{2} + s^{2}).
  $$ Now $f(u_{0}) + (u-u_{0}) f'(u_{0})$ grows like 
   $(1-\sqrt{6})u_{0}^{\sqrt{6}}$, which is large and negative,  while
    $f(u)$ grows like $u^{-\sqrt{6}}$, which is large and positive.
   Thus the left hand side is bounded below by a multiple of $u_{0}^{\sqrt{6}}+
   u^{-\sqrt{6}}$, or equivalently $s_{0}^{2} + s^{2}$.

 \subsection{Estimates for derivatives of  $\tau$}\label{ss:dtau-est}
 
 In this section we obtain estimates for the derivatives of a section
 $\tau=\tau_{H_{0},\theta_{0}}$. We suppose that $H_{0}>0$, so
 $\tau$ is defined over $G^{+}$. Recall that given $\ux\in G^{+}$ we write
 $\ux'$ for the point with co-ordinates $(Q,H_{0},\theta_{0})$ where $Q=Q(\ux)$.
  In the fixed trivialisation of ${\cal
 L}_{1}$ we write $\tau=\exp(-A)$ as above. We fix a positive integer $r$
 and $c>0$.  The result we prove is:
 \begin{prop}\label{prop:dtaudq}
 For any $\tilde\alpha>0$ there is a constant $C$ such that
 $$   \vert \ux\vert^{2r}  \Bigl| \Bigl({\frac{\partial}{\partial
Q}}\Bigr)^{r} A\Bigr|
 \leq C   E_{\tilde\alpha}^{-1}, $$
 in the set where $\vert\ux\vert, \vert \ux'\vert \geq 1$ and 
 $Q(\ux)<-c$  if $ x_{0}<0$.
 \end{prop}
 
 It is not hard to deduce Proposition \ref{prop:dtau-est} from this. The simplest case is
 the estimate on $\vert \db \tau\vert $. Since $\tau$ is holomorphic along the quadric
 surfaces we have
   $$  \vert \db \tau\vert = \psi \left\vert \frac{\partial \tau}{\partial Q}\right\vert
      = \psi \left\vert \frac{\partial A}{\partial Q} e^{-A} \right\vert \leq C \psi 
      \left\vert \frac{\partial A}{\partial Q}\right\vert E_{\alpha}, $$
      using Proposition \ref{prop:tau-est}. Now $\psi \leq C \epsilon p^{4} \leq C \epsilon
      \vert \ux \vert^{2}$ so Proposition \ref{prop:dtaudq} yields
      $$  \vert \db \tau \vert \leq \epsilon\,C E_{\tilde\alpha} $$
      (over the given set) for some $\tilde\alpha$ smaller than $\alpha$.
      The other estimates in Proposition \ref{prop:dtau-est} are obtained similarly. Using
      the fact that $\tau$ is holomorphic along the surfaces we can bound
      the partial derivatives in the $(H,\theta)$ directions in terms of
      $|\tau|$ (via either elliptic theory or the Cauchy integral formula). Thus we can estimate
      any partial derivative of $\tau$ by the derivatives in the $Q$ variable.
      We leave the details to the reader.

 We follow the same pattern as in the previous subsection, proving Proposition
 \ref{prop:dtaudq} first
 on the cone $\{Q(\ux)>0, x_{0}>0\}$. Again we can exploit homogeneity under
 scaling $\ux\mapsto \lambda \ux$,  $H_{0} \mapsto \lambda^{3} H_{0}$. Thus
 we begin by considering the restriction of $ \vert \ux \vert^{2r} \bigl(\frac{\partial}
 {\partial Q}\bigr)^{r} A$ to the surface $Q(\ux)=1$, $x_{0}>0$. The $\theta$--variable
 will play  essentially no role, so we suppose $\theta_{0}=0$ and restrict
 to the arc $\Gamma$ where $\theta=0$. We recall from the previous subsection
 that we have two useful co-ordinates
 along this arc, one the function $v$ and the other the arc length $s$.
 In what follows we will also have to bring in a third co-ordinate, the
 restriction of the function $H$. Recall also that the fixed parameter
 $H_{0}$ corresponds to values $v_{0}, s_{0}$~--~ie, the co-ordinates of
 the point $\ux'$ in the different parametrisations of the arc. For a suitable fixed
 $N$ we write
 $$   R(s,s_{0}) = \left( \frac{1+s}{1+s_{0}}\right)^{\!\!N} + \left( 
 \frac{1+s_{0}}{1+s}\right)^{\!\!N}.
 $$ With all these preliminaries out of the way, what we actually prove
 is:
 \begin{prop}\label{prop:dtaudq2}
 For any $r$ there are $N, C$ such that on the arc $\Gamma$
 $$  \vert \ux \vert^{2r} \Bigl\vert \Bigl(\frac{\partial}{\partial Q}\Bigr)^{r}A
 \Bigr\vert\leq C R(s,s_{0}) (s-s_{0})^{2}. $$
 \end{prop}
 To see that this implies Proposition \ref{prop:dtaudq} in the positive cone, we argue as follows.
 For any point $\ux$ in this cone we define $s$ to be the length of the
 obvious arc in the quadric surface through $\ux$, running from the $x_{0}$
 axis to $\ux$. Similarly we define $s_{0}$ to be the length of the arc
 in the same quadric surface to the point $\ux'$. Thus $S(\ux)= s-s_{0}$.
 The function $ \vert \ux \vert^{2r} \vert \bigl(\frac{\partial}{\partial Q}\bigr)^{r}
 A\vert$ is homogeneous of degree $3$ under rescaling, while $s, s_{0}$ are 
 homogeneous of degree $3/2$.
 Thus the estimate in Proposition \ref{prop:dtaudq2} scales to the general estimate
 \begin{equation} 
 \vert \ux \vert^{2r} \Bigl\vert \Bigl(\frac{\partial}{\partial Q}\Bigr)^{r}A
 \Bigr\vert \leq C R(Q^{-3/2}s, Q^{-3/2}s_0)\, (s-s_{0})^{2}.
 \label{eq:scaled} \end{equation}      
 We use the following:
 \begin{lem}
 For any $b,\,\beta>0$ there is a $C$ such that
     $$ R(Q^{-3/2}s, Q^{-3/2}s_0)\, (s-s_{0})^{2}\, e^{-\beta(s-s_{0})^{2}}\leq
      C$$
      provided that $Q\geq b$ or $s,\,s_{0}\geq b$. \end{lem}
      The proof is elementary and left to the reader. We can obviously choose $b$ so that
      $\vert \ux\vert\geq 1$, $\vert \ux'\vert\geq 1$ implies that $Q\geq b$ or
      $s, s_{0}\geq b$. Hence the Lemma and Equation~(\ref{eq:scaled}) imply
      Proposition \ref{prop:dtaudq} in the positive cone.
      
      We now turn to the heart of the matter: the proof of Proposition
      \ref{prop:dtaudq2}. The complication here is the interaction between the three co-ordinates
      $H, v, s$ on $\Gamma$.  For a function $f$ on $\Gamma$ and two points
      $\ux,\ux'$ on $\Gamma$ we write
      $$  \Delta(f; \ux, \ux')= \Delta_{f}(v,v_{0}) $$ where
      on the right hand side we understand that we use the co-ordinate
      $v$ to parametrise $\Gamma$ and $v,v_{0}$ are the co-ordinates
      of $\ux,\ux'$. 
      \begin{lem}\label{lem:deltaf}
      Suppose $f$ is a smooth function on $\Gamma$ and
       $f\sim H^{\mu}$, $\frac{df}{dH}\sim \mu H^{\mu-1}$, $\frac{ d^{2}
      f}{dH^{2}} \sim \mu (\mu-1) H^{\mu-2}$. Then for a suitable $N$ depending
      on $\mu$ we have:
      $$ \vert \Delta(f;\ux,\ux')\vert\leq C R(s,s_{0}) (s-s_{0})^{2} (1+
      H)^{\mu-1},$$
      where $H$ corresponds to the point $\ux$.
      \end{lem}
      To see this we express $f$ as a function of $v$, $f\sim C v^{\sqrt{6}\mu}$.
      We have $$\frac{d^{2}f}{dv^{2}} \sim  C v^{\sqrt{6}\mu-2}$$
      The integral formula Equation~(\ref{eq:integralform}) gives
      $$ \vert \Delta(f;\ux,\ux')\vert \leq C (v-v_{0})^{2} (1+ v_{*})^{\sqrt{6}\mu
      -2}, $$
      where $v_{*}$ is one of $v,v_{0}$ (which one depending on the sign
      of $\sqrt{6}\mu -2$ and which of $v, v_{0}$ is the larger). The function
      $s$ is asymptotic to a multiple of $v^{\sqrt{3/2}}$, hence
      $$   \vert s-s_{0}\vert \geq C \vert v-v_{0}\vert
       (1+ v_{**})^{\sqrt{3/2}-1}, $$
       where $v_{**}$ is the smaller of $v,v_{0}$. Then the result follows
       by elementary arguments. 
       (The point is that introducing the function $R$
       allows us to essentially interchange $v,v_{0}$ in our estimates.)

       Now consider the function $A=A(Q,H,H_{0})$. By construction this
       satisfies 
         $$
          A(Q,H_{0}, H_{0})=0 ;\ \ \frac{\partial A}{\partial H}\Big\vert_{H=H_{0}}=0.
          $$
          In other words, $A$ vanishes to second order along the \lq\lq
          diagonal'' $H=H_{0}$. Differentiating $r$ times with respect to
          $Q$, we see that $\bigl(\frac{\partial}{\partial Q}\bigr)^{r}A$
          also vanishes to second order along the diagonal. This means that
          on the arc $\Gamma$ it is equal to
          $$ \Delta( \Bigl( \frac{\partial}{\partial Q}\Bigr)^{r} A; \ux,\ux').
          $$
          Thus we see that, on $\Gamma$,
          $$  \Bigl( \frac{\partial}{\partial Q}\Bigr)^{r} A = B_{1}-B_{2},$$
          where
          $$  B_{1} = \Delta( \Bigl(\frac{\partial}{\partial Q}\Bigr)^{r}
          \frac{p^{6}}{18}; \ux,\ux') $$
          $$  B_{2} = \Delta( 
          \Bigl(\frac{\partial}{\partial Q}\Bigr)^{r} \frac{ H_{0}F^{+}(H,Q)}{F^{+}(H_{0},Q)};
          \ux,\ux').$$
          Now $f=\bigl(\frac{\partial}{\partial Q}\bigr)^{r} p^{6}$ is a
          homogeneous function of degree $3-2r$ on $\bR^{3}$. It follows
          that $f\sim C H^{\lambda}$ on $\Gamma$, where $\lambda= 1-\frac{2r}{3}$
          (since $H$ is homogeneous of degree $3$); similarly for
          the derivatives of $f$. Applying Lemma \ref{lem:deltaf}
          we see that
              $$\vert B_{1}\vert \leq C R(s,s_{0}) (s-s_{0})^{2} (1+H)^{-2r/3}$$
            Now on $\Gamma$, $\vert \ux\vert^{2} \leq C (1+H)^{2/3}$ so
            we obtain
            $$ \vert \ux\vert^{2r} \vert B_{1}\vert \leq C R(s,s_{0}) (s-s_{0})^{2},
            $$
            which is just the form of estimate we need. 
            
            The term $B_{2}$ is more complicated. Regard $v$ as a function
            of $H$~--~taking $Q=1$. Then we can write
            $$  \frac{H_{0} F^{+}(H,Q)}{F^{+}(H_{0}, Q)} = Q^{\sqrt{3/8}} v(H/Q^{3/2})
            \frac{H_{0}}{ Q^{\sqrt{3/8}} v(H_{0}/Q^{3/2})}. $$
            Set
            $$  f_{p} = \left(\frac{\partial }{\partial Q}\right)^{p} Q^{\sqrt{3/8}}
            v(H/Q^{3/2}), $$
            $$   g_{q} = \left( \frac{\partial}{\partial Q}\right)^{q}
            \frac{H_{0}}{ Q^{\sqrt{3/8}} v(H_{0}/Q^{3/2})}. $$
            Then $f_{p}$, $g_{q}$ are smooth functions on $\Gamma$ (ie, 
            we set $Q=1$ after performing the differentiation). We have
            $$  B_{2} = \sum_{p+q=r}  g_{q}(H_{0})\Delta(f_{p}; \ux, \ux'). $$    
          
          Now, regarded as a function of $x_{0}$, it is easy to see that
          $v$ has a series expansion for $x_{0}$ large:
          $$  v=  x_{0}^{\sqrt{3/2}}(a_{0} + a_{1} x_{0}^{-2} +\dots ).
          $$
          This means that $v(H)$ has an expansion
          $$  v(H)= H^{1/\sqrt{6}} ( b_{0} + b_{1} H^{-2/3} + \dots). $$
          $$  Q^{\sqrt{3/8}} v(H/Q^{3/2}) = H^{1/\sqrt{6}} ( b_{0} + b_{1}
          Q H^{-2/3} +\dots). \leqno{\rm Hence} 
          $$
          So we see that $f_{p}\sim b_{p} H^{1/\sqrt{6}- 2p/3}$. Applying
          Lemma \ref{lem:deltaf} we get
          $$ \vert \Delta(f_{p}; \ux,\ux') \vert \leq C R(s,s_{0}) 
          (s-s_{0})^{2} (1+ H)^{1/\sqrt{6} -1- 2p/3} . $$
            $$  \vert g_{q}\vert \leq C (1+ H_{0})^{1-1/\sqrt{6}-2q/3}. 
            \leqno{\rm Similarly}
            $$
            $$ (1+H)^{2r/3} \vert g_{q} \Delta(f_{p}, \ux,\ux') \vert \leq
            C R(s,s_{0}) (s-s_{0})^{2} \left(\frac{1+ H}{1+H_{0}}
            \right)^{\frac{1}{\sqrt{6}}-1+\frac{2q}{3}}. \leqno{\rm So}
            $$
            Now changing the value of $N$ suitably, the power of $(1+H)/(1+H_{0})$
            can be absorbed into $R(s,s_0)$ and we get
            $$  \vert \ux\vert^{2r} \vert g_{q} \Delta(f_{p}, \ux,\ux')\vert
            \leq C R(s,s_{0}) (s-s_{0})^{2} ,$$
            Hence $\vert \ux\vert^{2r} B_{2} $ is bounded by a multiple
            of $R(s,s_{0}) (s-s_{0})^{2}$ and we have finished the proof
            of Proposition \ref{prop:dtaudq2} over the positive cone.

         We omit the details of the extension of this argument to 
          the region $Q(\ux)<0$.  Let us just explain where the condition
          $Q(\ux)<-c$ enters, if $x_{0}\leq 0$. Using homogeneity we can
          throw the calculations onto the quadric $Q(\ux)=-1$. We consider
          the arc $\theta=0$ in this quadric on which we have arc length
          co-ordinates $s$ for the point $\ux$ and $s_{0}$ for the point
          $\ux'$. Alternatively, we can use the co-ordinates $H, H_{0}$.
          Then $H_{0}$ and $s_{0}$ are positive by hypothesis. The problem
          comes when $H$ and $s$ are large and {\it negative}.
          The function $u(H)$ for large positive $H$ has a series
          expansion       
          $$  u(H) = H^{1/\sqrt{6}}( b_{0} + b_{1} H^{-2/3} + \dots ),
          $$ just as before. For large negative $H$ on the other hand the
          series is
               $$u(H)= u(-H)^{-1} = (-H)^{-1/\sqrt{6}} (b_{0}^{-1} + \dots).
               $$
               This means that the ratio $H_{0} \frac{F^{+}(H,Q)}{F^{+}(H_{0},
               Q)}$, for $H_{0}\gg 0$ and $ H\ll 0$, is
               $$  \frac{ H_{0}^{1-1/\sqrt{6}} (-H)^{-1/\sqrt{6}} (-Q)^{\sqrt{3/2}}}
               {(b_{0} - b_{1} Q H_{0}^{-2/3} + \dots)(b_{0}- b_{1} Q
               (-H)^{-2/3} + \dots)}. $$
               The presence of the term $(-Q)^{\sqrt{3/2}}$ makes for the
               difference with the previous case.
               When we differentiate $r$ times this term contributes so
               we only get the bound:
               $$ \left(\frac{\partial}{\partial Q}\right)^{r} H_{0} 
               \frac{F^{+}(H,Q)}{F^{+}(H_{0}, Q)} \leq C  
               H_{0}^{1 -1/\sqrt{6}} (-H)^{-1/\sqrt{6}}. $$
               This means that  we get
               $$  \vert \ux\vert^{2r} \vert B_{2}\vert\leq C
               {(-s)}^{4r/3-\sqrt{2/3}} \, s_{0}^{2-\sqrt{2/3}}.$$
              Now scaling back  and
              using homogeneity the derivative bound becomes
              $$  {(-Q)}^{\sqrt{3/2} - r}\,{(-s)}^{4r/3-\sqrt{2/3}} \,
              s_{0}^{2-\sqrt{2/3}}. $$
              If $r\geq 2$ this blows up as $Q\rightarrow 0$ for fixed
              $s<0$, $s_{0}>0$. (As we know
              it must since the functions are only H\"older continuous
              along the null cone.) On the other
              hand if $Q<-c$ then we can proceed to obtain a subexponential
              bound much as before.  
              We leave it to the reader to check that the additional
              subtleties induced by the presence of $F^-$ in the
              definition of $\tau$ for $Q<0$ (Equation (\ref{eq:hattau2}))
              do not affect things in any significant manner.

\subsection{Estimates for $s$}\label{ss:s-est}

Given a point $(\ux',t')$ in $\bR^{4}$ with $\vert \ux'\vert>3$ we have
defined a section $s=s_{\ux',t'}$ of ${\cal L}$. For $\alpha>0$ define
 a function
 \begin{equation} F_{\alpha} = \exp(-\alpha\left( \psi_{0}^{-2}(Q-Q_{0})^{2}
 + \psi_{0}^{2}(t-t')^{2}\right)\,)\end{equation}
 $$\Psi(\ux,\ux')= \frac{\psi}{\psi_{0}} + \frac{\psi_{0}}{\psi},
 \leqno{\rm Also\ define} $$ 
    $$     \delta(\ux,\ux')= \Bigl\vert \frac{\psi}{\psi_{0}}- \frac{\psi_{0}}{\psi}\Bigr\vert.
    \leqno{\rm and}$$
  The goal of this subsection
is to prove the following:
\begin{prop} \label{prop:s-est}
   For any $r,\,c$ 
   there are $\alpha,\, C$ such that for $p\leq r$\begin{itemize}
   \item $\vert \nabla^{p} s\vert \leq C \Psi(\ux,\ux')^{p} E_{\alpha} F_{\alpha} $
    everywhere,
   \item $\vert \nabla^{p}(\db s)\vert \leq  C (\epsilon + \delta(\ux,\ux'))
   \Psi(\ux,\ux')^{p+1} E_{\alpha}F_{\alpha}$
   throughout $\{ (\ux, t): \vert \ux \vert \geq 3\}$.
    \end{itemize}
   
    \end{prop}

The proof of this will require a number of steps. For simplicity we will
just prove the estimate on $\vert \db s\vert$~--~the extension to higher
derivatives is straightforward (using the appropriate results from \S \ref{ss:dtau-est}).
Since $s= \tau^{*} \otimes \rho$ we have
\begin{equation}   \vert \db s \vert \leq \vert \db \tau^{*}\vert \vert \rho\vert + \vert \tau^{*}
\vert \vert \db \rho\vert. \label{eq:split}\end{equation} 
Throughout this subsection and the next we will make frequent use of the bounds on the
derivative of the function $\psi$. Note that we have
\begin{equation} \frac{\partial p}{\partial Q}= \frac{3x_{0}^{2}+ Q}{p^{7}}
= O(p^{-3}) \end{equation}
(by differentiating (\ref{eq:cubic}) and (\ref{eq:implicitp}) with $H$ fixed),
while Equation~(\ref{eq:partialH}) gives
\begin{equation} \frac{\partial p}{\partial H} = \frac{3x_0}{p^7} = O(p^{-5}).\end{equation}
Thus Lemma \ref{lem:defpsi} implies that \begin{equation}
  \psi^{-1} \Bigl|\frac{\partial \psi}{\partial Q}\Bigr|\leq C \epsilon p^{-1}
,\ \ \psi^{-1} \Bigl|\frac{\partial \psi}{\partial H}\Bigr|\leq C \epsilon p^{-3}. \label{eq:slowvar} \end{equation}

\begin{lem}\label{lem:dbtaustar}
For suitable $C,\alpha$ we have
$$ \vert \db \tau^{*} \vert \leq  C \epsilon E_{\alpha} $$
when $\vert \ux \vert \geq 3$. 
\end{lem}
Recall that when $\vert \ux\vert>3$ the section $\tau^{*}$ is equal to
$\hat{L}^{+} \tau$, so 
\begin{equation}
\vert \db(\tau^{*})\vert \leq  \vert \hat{L}^{+}\db \tau \vert + \vert \nabla \hat{L}^{+}
\vert \vert \tau \vert . \label{eq:obvious}\end{equation} 
There are two issues here. The first issue is that the estimates of
Proposition~\ref{prop:dtau-est} for $\db \tau$
only hold in a region $\Omega^{+}_{c}$. However $\hat{L}^{+}$ vanishes at
points where $x_{0}<0$ and $Q(\ux) > - \frac{b_2\psi}{2\epsilon}$. Since
$\psi\geq \epsilon$ we see that the support of $\hat{L}^{+}$ lies in $\Omega^{+}_{c}$
with $c= b_{2}/2$ and the estimates of Proposition \ref{prop:dtau-est} deal with the
 first term in Equation~(\ref{eq:obvious}). The second issue concerns the
 term involving $\nabla \hat{L}^{+}$. 
Thus it suffices to show that
$$  \vert \nabla (\lambda(\frac{\epsilon}{b_{2}} \frac{Q}{\psi})) \vert \leq C
\epsilon. $$
 The derivative of $\lambda (\frac{\epsilon}{b_{2}} \frac{Q}{\psi})$ vanishes
 if $\vert Q\vert > \psi b_{2}/\epsilon$.
 Since the function $\lambda$ has bounded derivative it suffices to show
 that 
 $$\vert \nabla \Bigl(\frac{\epsilon}{b_{2}} \frac{Q}{\psi}\Bigr) \vert\leq
 C \epsilon, $$ when $\vert Q\vert \leq \psi b_{2}/\epsilon$.   
The relevant components of $\nabla$ with respect to the standard
 orthonormal basis of tangent vectors for our metric $g$ are
 $\psi \frac{\partial}{\partial Q}$ and $pr \frac{\partial}{\partial H}$.
  Consider first the $Q$ derivative. We have
  $$   \Big\vert \psi \frac{\partial}{\partial Q} \Bigl(\frac{\epsilon}{b_{2}} \frac{Q}{\psi}\Bigr)
 \Big\vert  = \frac{\epsilon}{b_{2}} \Big\vert 1- Q \psi^{-1} \frac{\partial \psi}{\partial
 Q} \Big\vert \leq \frac{\epsilon}{b_{2}} + \Big\vert \frac{\partial \psi}{\partial
 Q}\Big\vert \leq \frac{\epsilon}{b_{2}} + C \epsilon p^{-1} {\psi}. $$
 Now $\psi \leq C p $ by the third item of Lemma \ref{lem:defpsi}, so we are done.
 
 For the $H$ derivative we have similarly:
 $$  \Big\vert p r \frac{\partial}{\partial H} \Bigl(\frac{\epsilon}{b_{2}}\frac{Q}{\psi}\Bigr)
 \Big\vert = \Big\vert pr \frac{\epsilon}{b_{2}} \frac{Q}{\psi^2} \frac{\partial \psi}{\partial
 H}\Big\vert \leq pr \psi^{-1} \Big\vert \frac{\partial \psi}{\partial H}\Big\vert
 \leq C \epsilon,$$ 
which completes the proof of Lemma \ref{lem:dbtaustar}.

We now turn attention to the section $\rho$. We begin with $\hat{\rho}$.
\begin{lem} \label{lem:dbrhohat}
For any $\alpha<1$ there is a constant $C$ such that
$$  \vert \db \hat{\rho}\vert \leq C  \delta(\ux,\ux') F_{\alpha}.
 $$\end{lem}

In our standard orthonormal frame, and the given trivialisation
of ${\cal L}_{2}$, 
$$\db\hat{\rho}= \left(
  \psi\left( \frac{\partial}{\partial Q} + \frac{i}{2}(t-t')\right) + i \psi^{-1}
\left( \frac{\partial}{\partial t} - \frac{i}{2} (Q-Q_{0})\right) \right)
\hat{\rho}$$
 $$ \hat{\rho} = \exp\left(-\frac{\psi_{0}^{2}(t-t')^{2} +
\psi_{0}^{-2}(Q-Q_{0})^{2}}{4}\right).
\leqno{\rm where}
$$
$$  \frac{1}{2}\left(\frac{\psi_{0}}{\psi}- \frac{\psi}{\psi_{0}}\right) \left( \frac{Q-Q_{0}}{\psi_{0}}
- i\psi_{0} (t-t')\right) \hat\rho. \leqno{\rm This\ is}
$$
The Lemma follows from the fact that for any $\alpha<1$ there is a $C$ such that
$$   A e^{-A^{2}} \leq C e^{-\alpha A^{2}}. $$
Next we have the following:
\begin{lem}\label{lem:dbrho}
    For any $\alpha<1$ there is a constant $C$, depending on $b_{1}$, such that
    $$  \vert \db \rho\vert \leq C (\epsilon + \delta(\ux,\ux')) 
    \Psi(\ux,\ux') F_{\alpha} , $$
    in the set where $\vert \ux\vert>3$.
    \end{lem}
    Given the preceding lemma, we just have to estimate the derivative
    of the cut-off function $\chi(\frac{\epsilon}{b_{1}} \frac{Q-Q_0}{\psi_{0}})$.
    This is bounded in modulus by $C\frac{\epsilon}{b_{1}} \frac{\psi}{\psi_{0}}$,
    which gives the desired result.

    The main result (Proposition \ref{prop:s-est}) in the case of $\vert \db s\vert$  follows from
    Equation~(\ref{eq:split}) and Lemmas \ref{lem:dbtaustar} and \ref{lem:dbrho}, since we clearly have
    $$  \vert \tau^{*} \vert
\leq \vert \tau \vert \leq C E_{\alpha},\ \  \vert \rho \vert \leq \vert \hat{\rho}\vert
= F_{1}. $$
Notice that if we estimate $\db \tau^{*}$
 over the region $\vert \ux\vert \leq 2$ we get a new term involving the
 cut-off function $\chi(2/\vert\ux\vert)$ and our estimate is only
  as good as that on the full covariant derivative $\nabla \tau^{*}$.
This is why we only consider the case $|\ux|\ge 3$ in the
second half of Proposition \ref{prop:s-est}.

\subsection{Estimates on sums}\label{ss:sum-est}

    For each point $(\ux',t')$ with $\vert \ux'\vert \geq 3$ we have now
    got a section $s_{\ux', t'}$ obeying  estimates expressed in terms of functions
    $E_{\alpha}, F_{\alpha}, \Psi(\ux,\ux'), \delta(\ux,\ux')$.
    Moreover, $s_{\ux',t'}$ is supported in a set $S(\ux')\times \bR$ where
    $S(\ux')$ is the set of points $\ux$ in $\bR^{3}$ which satisfy the
    conditions \begin{itemize}
    \item  $\vert Q(\ux)- Q(\ux')\vert \leq \frac{2b_{1}}{\epsilon}
    \psi(\ux')\ ,$ \item
     $Q(\ux)\leq -\frac{b_{2}}{2\epsilon}\psi(\ux)$ 
    if $x_0$ and $x'_0$ have different signs, \item
    $  \vert \ux\vert \geq 1. $ \end{itemize}
    
    Given $\ux$ with $\vert \ux\vert \geq 1$ let $N(\ux)\subset \bR^{3}$
    be the set
    $$   N(\ux)= \{\ux': \vert \ux'\vert \geq 3,\ \ux \in S(\ux')\}. $$    
      The modulus of the section $s_{\ux',t'}$  at the point
    $(\ux',t')$ is $1$ and it is quite clear from the constructions that
    the section is not small on a ball (in the metric $g$) of uniform size.
    Let us say $\vert s_{\ux',t'}\vert \geq C^{-1} $ on the ball of radius
    $1/10$ centred at $(\ux',t')$. Our goal in this subsection is to prove:
    \begin{prop} \label{prop:sum-est} For any $\alpha$
         we can find a countable collection of points\break ${(\ux'_{i}, t'_{i})}_{i\in I}$
         with $\vert \ux'_{i}\vert \geq 3$ having the following properties:
         \begin{itemize}
         \item The balls $B_{i}$ of radius $1/10$ centred at the $(\ux'_{i},t'_{i})$
         cover $\{ (\ux,t) : \vert \ux \vert \geq 4\}$.
         \item Let $E_{\alpha,i}, F_{\alpha, i}$ denote the functions associated
         with these points and write $\delta_{i}= \delta(\,\cdot\,, \ux'_{i})$ and
         $\Psi_{i}= \Psi(\,\cdot\,, \ux'_{i})$. Then for any $p$ there is a $C$
         such that 
         $$ \sum_{i,\,\,\ux'_{i}\in N(\ux)} \delta_{i} E_{\alpha,i} F_{\alpha, i} \Psi_{i}^{p} \leq C\epsilon,
          $$
          and $$ \sum_{i,\,\,\ux'_{i}\in N(\ux)} E_{\alpha, i} F_{\alpha, i} \Psi_{i}^{p} \leq
          C. $$
         \item There is a constant $K$ such that for all $D>1$
          we can divide the index set $I$ into at most $K D^4$ disjoint subsets
          $I_{\mu}$, such that if $(\ux,t)$ is
          contained in a ball $B_{i}$ for $i\in I_{\mu}$ then
          for any $p$ there is a $C$ such that
            $$  \sum_{j\in I_{\mu},\,\,j\neq i,\,\,\ux'_{j}\in N(\ux)} \delta_{j} E_{\alpha, j}(\ux,t) F_{\alpha,
            j}(\ux,t) \Psi_{j}^{p}\leq C \epsilon\, e^{-D}, $$
             $$ \sum_{j\in I_{\mu},\,\,j\neq i,\,\,\ux'_{j}\in N(\ux)} E_{\alpha, j}(\ux,t) F_{\alpha,
             j}(\ux,t) \Psi_{j}^{p} \leq C\, e^{-D}. $$
\end{itemize}  \end{prop}

Notice that this Proposition does not involve the sections we have constructed,
 only the
geometry of the metric $g$ and the functions $F_{\alpha}, E_{\alpha}, \Psi,
\delta$.

To begin the proof of Proposition \ref{prop:sum-est} we consider the restriction of
the metric $g$ to the $(x_{0}, x_{1})$--plane. We first choose a sequence
of points on the $x_{0}$--axis such that the $\frac{1}{20}$--balls about these points
cover the portion $\vert x_{0}\vert >3$ of this axis. It is easy to check
then that the corresponding $\frac{1}{10}$--balls cover the neighbourhood $T=\{pr< \delta\}$ for some small
$\delta$. We then choose a collection of points in the half-plane $x_{1}>0$
 and outside $T$ such that the $\frac{1}{10}$--discs (in the metric $g$)
  about these points cover the complement $U$ of $T$ and the Euclidean
  ball $\{x_{0}^{2}+ x_{1}^{2} \leq 9\}$ in the half-plane.
 We denote the centres obtained in this way by $P'_{j}$
 and the $\frac{1}{10}$--discs by $D_{j}$. It is fairly clear that we can do this in such a way that any intersection
of more than $n$ discs $D_{j}$ is empty, for some fixed $n$.
    
    We now move to $3$--space. We use the balls centred on the axis
     to cover the relevant portion of the $x_{0}$--axis in $3$--space
    in the obvious way. Recall that the length of the circle orbit
    under rotations is $2\pi pr$. It is straightforward to check that there
    is a constant $R$ such that for each point $P'_{j}$ which is not on
    the axis
    $$ \frac{ \max_{D_{j}\cap U}(pr)}{ \min_{D_{j}\cap U}(pr)} \leq R.$$
    As a consequence of this we can, for each such $P'_{j}$, choose an
    integer $m_{j}$ which is comparable to $pr$ for all points in $
    D_{j}\cap U$. Then we get a cover of $\bR^{3}$, minus the Euclidean ball
    of radius $3$, in the following way. We take the images of these points
    $P'_{j}$ under
    rotations through multiples of $ 2\pi/M m_{j}$ for suitable fixed $M$,
    and the balls of radius $\frac{1}{10}$ centred on these points. In this way
    we get a collection of $\frac{1}{10}$--balls $B_{k}$ with centres $\ux'(k)$ in $\bR^{3}$ such that
    \begin{itemize}
    \item The balls $B_{k}$ cover $\{\vert \ux\vert>3\}$,
    \item The centre of any ball $B_{k}$ either lies on the $x_{0}$--axis or is
    contained in the orbit of a $P'_{j}$  under a cyclic subgroup of the rotation
    group, where the order of the cyclic group is bounded by a fixed multiple
    of $pr$, evaluated at the centre.                 
\end{itemize}

     Next we move to $4$--space. Equation~(\ref{eq:slowvar}) above shows that
     \begin{equation}   \psi^{-1} \vert \nabla \psi\vert \leq C \epsilon. \end{equation}
     This means that, once $\epsilon$ is sufficiently small, we can suppose
     that
     $$  \frac{ \max_{B_{k}} (\psi)}{\min_{B_{k}} (\psi)}\leq \frac{11}{10}, $$
     say. We fix a constant $M'$ and for each centre $\ux'(k)$ we take a countable
     collection of points
     $$ \Bigl(\ux'(k), \frac{\nu}{ M'\, \psi(\ux'(k))}\Bigr), \ \  \nu \in \bZ. $$
     This finally gives us our collection of centres $(\ux'_{i},t'_{i})$
     in $\bR^{4}$.   For a suitable choice of the constants $M$ and $M'$ we can arrange that
         the $\frac{1}{10}$--balls about the $(\ux'_{i}, t'_{i})$ cover
          $\{ (\ux,t): \vert \ux \vert >3\}$.

          Now let $(\ux, t)$ be  a point with $\vert\ux\vert\geq 3$. We
          want to study the sum
          \begin{equation}
            B(\ux,t)= 
            \sum_{i} E_{i,\alpha}(\ux,t) F_{i,\alpha}(\ux,t) \Psi(\ux, \ux'_{i})^p\delta(\ux,\ux'_{i})
               \label{eq:sum}\end{equation}
              with the set of centres $(\ux_{i}', t_{i}')$ obtained above.
              The manner in which these centres were chosen allows us to easily sum over the $\theta$
              and $t$--variables.
              \begin{lem}\label{lem:arithmsum}
              Let $u_{i}$ be an arithmetic progression $u_{i}=A i + C$,
              $A>0$,
              labelled by $i\in \bZ$. Then there are universal constants
              $k_{0}, k_{1}$ such that for all $B>0$
              $$ \sum_{i\in\bZ} \exp(-\left(\frac{u_{i}}{B}\right)^{2}) 
              \leq k_{0} + k_{1}\frac{B}{A}. $$
              \end{lem}
              This is standard and elementary.  When we consider the contribution
              to the sum in Equation~(\ref{eq:sum}) from the centres which
              lie in the same orbit under the translation action
              we get terms precisely of the form considered in
              the Lemma (with $A=\frac{1}{M'\psi_0}$ and
              $B=(\alpha\psi_0^2)^{-1/2}$, where 
              $\psi_0=\psi(\ux'_i)$). Thus we can reduce to a 3-dimensional
              problem by summing over translation orbits (which yields at
              most a uniform constant factor).
              
              The rotation action can also be factored out in a similar
              way, but requires a more careful treatment. Let $$\lambda=
              pr=(4H^2+r^6)^{1/4}.$$ The centres in a same rotation orbit yield (finitely
              many) terms of the form considered in Lemma \ref{lem:arithmsum},
              but now $A=\frac{2\pi}{Mm_j}\sim \lambda(\ux'_i)^{-1}$, 
              while $B=(\alpha L(\ux)^2)^{-1/2}\sim \lambda(\ux')^{-1}$,
              where $\ux'$ is the point introduced in
              \S \ref{ss:tau-est}, lying on the same quadric as $\ux$ but
              with $H(\ux')=H(\ux'_i)$. Hence, denoting by $(Q,H)$ and
              $(Q_0,H_0)$ the co-ordinates of $\ux$ and $\ux'_i$ respectively,
              the factor $\Sigma$ resulting from summation over a rotation orbit
              satisfies \begin{equation}\label{eq:sigmabound}
              |\Sigma|\le \min\Bigl(C +C \frac{\lambda(Q_0,H_0)}{
              \lambda(Q,H_0)},\ C'\lambda(Q_0,H_0)\Bigr)\end{equation}
              (using Lemma~\ref{lem:arithmsum} and
              the fact that the number of centres in the orbit is of
              the order of $\lambda(Q_0,H_0)$). We now use:
              
              \begin{lem}\label{lem:sigmabound}
              There is a constant $C$ such that 
              $$|\Sigma|\leq C + C\psi_0^{-1}|Q-Q_0|,$$
              where $\psi_0=\psi(\ux'_i)$.
              \end{lem}
              
              There are several cases to consider. First assume
              that $Q_0\ge -|H_0|^{2/3}$. Then the co-ordinates of 
              $\ux'_i$ satisfy $|x_0|\ge cr$ for some $c\in (0,\frac{1}{2})$
              (the positive root of the equation $c^2+c^{2/3}=\frac12$). Hence
              $r\le C|H_0|^{1/3}$, and 
              $\lambda(Q_0,H_0)=(4H_0^2+r^6)^{1/4}\le C|H_0|^{1/2}$.
              On the other hand $\lambda(Q,H_0)\ge |2H_0|^{1/2}$, so
              we get a constant bound on $|\Sigma|$ using Equation (\ref{eq:sigmabound}).
              In the other case $Q_0\le -|H_0|^{2/3}$,
              the co-ordinates of $\ux'_i$ satisfy $|x_0|\le cr$, so
              $r\sim |Q_0|^{1/2}$ and $p\sim |Q_0|^{1/4}$, so 
              $\lambda\sim |Q_0|^{3/4}$. If $Q\le \frac{1}{2}Q_0$ then
              $$|\Sigma|\le C+C(Q_0/Q)^{3/4}$$ is bounded by a uniform
              constant. Otherwise, we have $|Q-Q_0|\ge \frac{1}{2}|Q_0|$, so
              $$|\Sigma|\le C'\lambda(Q_0,H_0)=C'pr\le C \psi_0^{-1} p^2 r\le
              C\psi_0^{-1}|Q_0|\le C \psi_0^{-1}|Q-Q_0|.$$
              This completes the proof of the Lemma. Since the
              factor $\psi_0^{-1}|Q-Q_0|$ can be absorbed
              into $F_\alpha$ up to an arbitrarily small modification of the
              constant $\alpha$, Lemma \ref{lem:sigmabound} allows us to
              sum over rotation orbits.

              Thus we can reduce to a 2--dimensional
              problem. For this we adapt our notation slightly. We regard
              $Q$ and $H$ as functions on $\bR^{2}$ in the obvious way
              and for $P$ in the half-space $x_{1}\geq 0$ in $ \bR^{2}$
               let $\uSigma(P)$ be the part of the
              corresponding quadric through $P$ which lies in the half-space.
              Thus $\uSigma(P)$ can be identified with the quotient of
              one of our quadrics in $\bR^{3}$ under the rotation action.
              We write $N(P)$ for the quotient of the corresponding set
              $N$ defined above. For each of the centres $P'_{j}$ we have
              chosen above we write $E_{j}(P),F_{j}(P), \delta(P,P'_{j}),
              \Psi(P,P'_{j})$ for the corresponding functions on the half-plane.
               
               To prove the second item of Proposition \ref{prop:sum-est}
               it suffices to prove:
               \begin{prop}\label{prop:sum-est2}
                  Let  $\{P'_{j}\}\in \bR^{2}$ be the set of centres constructed above.
                  Then there is a $C$ such that for any $P\in \bR^{2}$
                  $$\sum_{j: P'_{j}\in N(P)}  \delta(P,P'_{j})\,E_{j}(P)\,F_{j}(P) 
                   \,\Psi(P, P'_{j})^{p} \leq
                  C \epsilon $$
                  $$ \sum_{j:P'_{j} \in N(P)} E_{j}(P) \,F_{j}(P) 
                  \,\Psi(P,P'_{j})^{p} \leq C $$\end{prop}

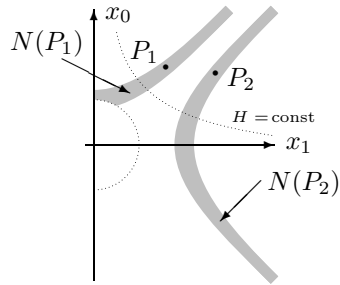
\begin{figure}[ht!]
\setlength{\unitlength}{6mm}
\centering
\begin{picture}(4,6)(0,-3)
\psset{unit=\unitlength}
\psbezier[linecolor=lightgray,linewidth=4pt](0,1.1)(1,1.1)(2,2)(3,3)
\psbezier[linecolor=lightgray,linewidth=4pt](0,0.9)(1,0.9)(2,2)(3,3)
\psbezier[linecolor=lightgray,linewidth=4pt](2.1,0)(2.1,1)(3,2)(4,3)
\psbezier[linecolor=lightgray,linewidth=4pt](1.9,0)(1.9,1)(3,2)(4,3)
\psbezier[linecolor=lightgray,linewidth=4pt](2.1,0)(2.1,-1)(3,-2)(4,-3)
\psbezier[linecolor=lightgray,linewidth=4pt](1.9,0)(1.9,-1)(3,-2)(4,-3)
\psarc[linewidth=0.4pt,fillstyle=solid,fillcolor=white,linestyle=dotted,
       dotsep=1pt](0,0){1}{-90}{90}
\psbezier[linewidth=0.4pt,linestyle=dotted,dotsep=1pt](0.5,2.5)(0.8,1)(2,0.5)(4,0.2)
\put(0,-3){\vector(0,1){6}}
\put(0.2,2.8){\small $x_0$}
\put(-0.2,0){\vector(1,0){4.2}}
\put(4.2,-0.1){\small $x_1$}
\put(1.6,1.72){\circle*{0.1}}
\put(0.75,1.9){\small $P_1$}
\put(2.7,1.6){\circle*{0.1}}
\put(2.9,1.3){\small $P_2$}
\put(-0.7,1.9){\vector(2,-1){1.5}}
\put(3.7,-0.8){\vector(-1,-1){0.9}}
\put(-1.9,2.1){\small $N(P_1)$}
\put(3.8,-1){\small $N(P_2)$}
\put(3,0.5){\tiny $H=$const}
\end{picture}
\caption{The set $N(P)$ 
(case 1: $Q(P)>-\frac{b_2}{2}$; case 2: $Q(P)<-\frac{b_2}{2}$)}\label{fig:NP}
\end{figure}

                The essential thing now is to understand the set $N(P)$.
                Notice first that if $P'\in N(P)$ and if the $x_0$
                co-ordinates $P_{0}$ and $P'_{0}$ have different signs
                then we have $Q(P) \leq -\frac{b_{2}}{2} \frac{\psi(P)}{\epsilon}$
                which implies that $Q(P) \leq -\frac{1}{2}b_{2}$. Thus we have the
                following ``quarter-space property'': if $Q(P)>
                -\frac{1}{2}b_{2}$ the sign of the co-ordinate $x_0$ on
                the whole of $N(P)$ is the same as that at $P$
                (see Figure \ref{fig:NP}).
                
                  Now, given $P$, let $\uSigma =\uSigma(P)$ be the
               the quotient of the quadric through $P$ as above.
                  We claim that $N(P)$ is contained in a 
                   ``thin neighbourhood'' of $\uSigma$. To state what
                   we need precisely: 
                   \begin{lem} If $b_1$ is sufficiently small then
                   for any point $P'$ in $N(P)$ 
                    the  corresponding level set $H^{-1}(H(P'))$ of $H$
                    meets
                   $\uSigma(P)$ in exactly one point $P''$,
                   and moreover if $\Gamma(P')$ is the connected arc of the
                   level set joining $P'$ to $P''$ then
                   $$  \frac{\max_{P^{*}\in \Gamma(P')} \vert P^{*}\vert }{\min_{P^{*}\in
                   \Gamma(P')}\vert
                   P^{*}\vert}\leq 11/10. $$
                  \end{lem}
                  This is fairly clear from a picture (see Figure 
                  \ref{fig:NP}), and can be verified
                  by routine calculations. 
                   Next we have the following:
                   \begin{lem}\label{lem:psivar}
                   If $b_{1}$ is sufficiently small then for any $P'$ in $N(P)$
                   we have $$\frac{\max_{P^{*}\in\Gamma(P')}\psi(P^{*})}
                   {\min_{P^{*}\in \Gamma(P')}
                   \psi(P^{*})} \leq 11/10 .$$  
                     \end{lem}
                     To prove this recall that by Equation (\ref{eq:slowvar})
                     and Lemma \ref{lem:defpsi} we have
                      $$ \Bigl|\frac{\partial \psi}{\partial
                     Q} \Bigr|\leq C \epsilon. $$
                    The variation of $Q$ over the connected arc $\Gamma(P')$
                    is at most $\frac{2b_{1}}{\epsilon} \psi(P')$. Integrating
                    over the arc we find  that for any point $P^{*}$ on $\Gamma(P')$,
                     $$ \vert \psi(P^{*})-\psi(P')\vert \leq C b_{1}\psi(P').
                     $$ Now we choose $b_{1}$ so small that $Cb_{1}\leq 1/50$ (say).

                     We now define a map $M$ from $N(P)$ to $\uSigma(P)\times
                     \bR$ by
                     $$   M(P')= (P'', Q(P')). $$
                     We define a metric $g_{0}$ on $\uSigma(P)\times \bR$ as follows.
                     In
                     the $\uSigma$ factor we take the metric induced by
                     $g$, and in the $\bR$ factor, with co-ordinate $Q$,
                     we take 
                     $$   \psi(P'')^{-2} dQ^{2}. $$
                     In other words, if we take $Q$ and $H$ as co-ordinates
                     we obtain the metric by \lq\lq freezing'' the coefficients
                     of $dQ^{2}$ and $dH^{2}$ at their values on $Q=Q(P)$.
                     Now by the quarter-space property above, the image of $M$
                      lies in a connected subset $\Sigma_{0}\times
                     \bR$ of $\uSigma\times \bR$, where 
                      $\Sigma_{0}$ lies in $\{P: \vert P\vert \geq c\}$
                      for some fixed $c>0$  depending on $b_{2}$.

                     \begin{lem}
                     If $b_{1}$ is sufficiently small then $M$ is an $11/10$ quasi-isometry
                     from the metric $g$ restricted to $N(P)$ to an open
                     subset in $\Sigma_{0}\times \bR$ with metric $g_{0}$.
                     \end{lem}
                     
                     For the $Q$ variable this follows from Lemma \ref{lem:psivar}. For the
                     $H$ variable we have to check that the variation
                     of $\log p r$ along the arc $\Gamma$ is small, which
                     follows from calculations similar to those above.

                     Now choose the arc length $s$ along $\uSigma(P)$ as
                     co-ordinate, taking the point $P$ as the origin $s=0$.
                     Thus we can regard the restriction of
                     $\psi$ to $\uSigma$ as a function $\psi(s)$. (This
                     notation is not really consistent with that used in
                     Section \ref{sec:sec3}, but we hope this will not cause confusion).
                     On $\Sigma_{0}$ we have
                     $$    \vert \frac{d\psi}{ds}\vert \leq C \epsilon
                     \psi$$
                     so if $s_{1}, s_{2}$ are the arc-length co-ordinates
                     of two points in $\Sigma_{0}$
                     \begin{equation}     \frac{\psi(s_{1})}{\psi(s_{2})} \leq e^{C
                     \epsilon \vert s_{1}-s_{2}\vert}. \label{eq:expbound}
                     \end{equation}
 We can now prove Proposition \ref{prop:sum-est2}. We just consider the first inequality,
                     the second being similar. The points $P'_{j}$ which contribute to
                     the sum  lie in $N(P)$ and we can map these by $M$
                     to get points $(s'_{j}, Q'_{j})$ in $\Sigma_0\times
                     \bR$. We use three facts:
                     \begin{itemize} \item The quasi-isometry property implies that
                     $ E(P,P'_{j})\leq \exp(-\alpha(s'_{j}-s)^{2})$,
                     for some $\alpha$.
                     \item  The function $\log \psi$  varies little
                     over the arcs $\Gamma(P')$, so we can replace
                     $\psi(P')$ by $\psi(P'')$ in estimating the sum.
                     \item The terms $\Psi(P, P'_{j})$ and $\psi(P)/\psi(P'_{j})$
                     appearing in the sum can be replaced by the exponential bound
                     Equation~(\ref{eq:expbound}) above.
                     \end{itemize}
                     
                     Putting all of this together,  it suffices to bound the sum
                     \begin{equation}
                      \sum \exp\left(-\alpha ((s'_{j})^{2}+\psi(s'_j)^{-2}
                       (Q'_{j}-Q)^2)\right) \exp (C\epsilon
                      \vert s'_{j}\vert) (\exp(C \epsilon\vert s'_{j}\vert)-1)
                      \label{eq:finalform} \end{equation}
                     Now it is easy to check that for any $\alpha'<\alpha$
                     we have an inequality
                     $$  ( e^{\epsilon A} - 1) e^{-\alpha A^{2}} \leq C
                     \epsilon e^{-\alpha' A^{2}}. $$
                     This means that, changing the value of $\alpha$ slightly,
                     it suffices to bound the sum  
                     \begin{equation}
                     \sum \exp(-\alpha((s'_{j})^{2} + \Bigl( \frac{Q'_{j}-Q}{\psi(s'_{j})}\Bigr)^{\!2}
                     )) \label{eq:simplified}\end{equation}
                     To do this we compare with the corresponding integral.
                     We consider the image $M(D_{j})$ of the $1/10$--disc centred on $P'_{j}$
                     under the map $M$ and let
                     $$  I_{j}(\beta)= \int_{M(D_{j})} e^{-\beta f} \frac{dQ'}{\psi(s')}
                     ds', $$
                     where
                     $$   f(s',Q')=
                     (s')^{2} + \Bigl(\frac{Q'-Q}{\psi(s')}\Bigr)^{\!2}
                     .$$
                      Now over $M(D_{j})$ the function $\psi(s')$ 
                      is essentially constant and the variations in
                      $Q'/\psi(s')$ and $s'$ are $O(1)$. It follows then
                      that there are constants $A$, $B$ such that
                      $$  \sup_{M(D_{j})} f \leq A f(s'_{j}, Q'_{j}) +
                      B . $$
                       This implies that
                       $$ e^{\beta B} I_{j}(\beta) \geq e^{-\beta A f(s'_{j},
                       Q'_{j})} \int_{M(D_{j})} \frac{dQ'}{\psi(s')} ds'.
                       $$
                       We take $\beta= \alpha/A$. Clearly
                       $$  \int_{M(D_{j})} \frac{dQ'}{\psi(s')} ds'\geq
                       c$$ for some fixed $c>0$. We see then that the sum
                       in Equation~(\ref{eq:simplified}) is bounded by a
                       multiple of 
                       $$ \sum_{j} I_{j}(\beta). $$
                       By construction of our open sets $D_{j}$, no more
                       than $n$ of the $M(D_{j})$ intersect, so
                       $$ \sum_{j} I_{j}(\beta) \leq n \int_{\bR^{2}}
                       e^{-\beta f} \frac{dQ'}{\psi(s')} ds'. $$
                       But this last integral can be evaluated explicitly
                       $$  \int_{\bR^{2}} e^{-\beta( s'^{2} + 
                       \left(\frac{Q-Q'}{\psi(s')}\right)^{2})} \frac{dQ'}{\psi(s')}
                       ds' = \frac{\pi}{\beta}. $$                     
                       This completes  the verification of the first two
                       items of Proposition \ref{prop:sum-est}. We omit the verification
                       of the third item which follows similar lines.

\section{Completion of proof}\label{sec:completion}
\subsection{Verification of Hypothesis 2}
\label{ss:hyp2}

In this subsection we will bring together the different strands of the
analysis in Sections \ref{sec:construct} and \ref{sec:estimates}
to complete the verification of Hypothesis 2.
The main issue we have to deal with is the fact that the model for
our neighbourhood $N$ of the zero set $\Gamma$ is a quotient of a tube in
$\bR^{4}$ under translations $t\mapsto t+ 2\pi\bZ\epsilon^{-1}$ whereas in
Sections \ref{sec:construct} and \ref{sec:estimates} we have worked
in $\bR^{4}$. To deal with this we go back to examine the definition of
the section $\hat{\rho}_{\ux',t'}$ in \S \ref{ss:secl2}.
To construct the line bundle corresponding to $\mathcal{L}_{2}$ on the quotient
space we proceed as follows. On $\bR^{4}$ we take a trivialisation of $\mathcal{L}_{2}$
in which the connection form is $-i (Q  + \frac{\epsilon}{2}) dt$. This $1$--form
is preserved by the translations so we get a line bundle with connection
over the quotient space in the obvious way. The factor $\frac{\epsilon}{2}$ means
that the holonomy is $-1$ around the zero set, as required. 
Now given $Q_{0}, H_{0}, t'$, the section $\hat{\rho}_{\ux',t'}$ we defined in
\S \ref{ss:secl2} is given, in this trivialisation, by
$$   \exp(-\frac{1}{4}\left(\psi_{0}^{2}(t-t')^{2} + \psi_{0}^{-2}(Q-Q_{0})^{2}\right))
\exp(iU)$$
where
$$  U=\frac{1}{2}(Q+Q_{0}+\epsilon)(t- t'). $$
 We now replace $t'$ by $t'_{\nu}=t'+ 2\pi\nu \epsilon^{-1}$
and form the sum 
\begin{equation}
    \Theta_{\ux',t'} = \sum_{\nu\in \bZ} \hat{\rho}_{\ux', t'_{\nu}},
\label{eq:deftheta}\end{equation}
working always in the fixed trivialisation of $\mathcal{L}_{2}$. Then
$\Theta_{\ux',t'}$ is a $2\pi \epsilon^{-1}$--periodic section.
Essentially these are the standard $\theta$--functions. 

The modulus of $\Theta_{\ux',t'}$ at the point $(\ux',t')$ is no longer
$1$. However it is very close to $1$, the difference is bounded by the
sum
   $$ 2 \sum_{\nu\geq 1} e^{-\pi^{2}\nu^{2}}, $$
   which is very small. More generally, the section $\Theta_{\ux,t'}$
   is very close to $\hat{\rho}_{\ux',t'}$ over a ball 
   (in the metric $g$) of radius $1/10$
   centred on $(\ux',t')$. This means that these sections have essentialy
   the same local behaviour as those considered before.    

The sections $\Theta_{\ux',t'}$ define sections of the corresponding line
bundle over the quotient space $N$ and we can repeat all the constructions
of Sections \ref{sec:construct} and \ref{sec:estimates} using these in
place of the $\hat{\rho}$. However it
easier to keep working in $\bR^{4}$. We can the reduce all the estimates
for this modified construction to those established before by the following
simple device. Recall that for any point $\ux'$, we have $\psi_{0}=\psi(\ux')\geq
\epsilon$. We can choose an integer $q$ such that 
$$q\leq \frac{\psi_{0}}{\epsilon}\leq
2q. $$
Now  we modify the  construction in \S \ref{ss:sum-est}, when we go from a covering in 3-space to
a covering in 4--space, slightly. We have centres $\ux'(k)$ in $\bR^{3}$ as
before and we take the sequence of centres
$$ (\ux'(k), \frac{\nu}{N q \epsilon}) \ \ \ \nu \in \bZ,$$
    where $N$ is some suitable fixed integer (independent of $\ux'(k)$, while
    $q$ depends on $\psi(\ux'(k))$).
The separation between these centres, in the metric $g$, is $\frac{\psi}{
Nq\epsilon}$ which lies between $N^{-1}$ and $2N^{-1}$: bounded above and below
independently of $\ux'(k)$. When we estimate the sum over these centres and
combine with the sum involved in the definition of $\Theta_{\ux',t'}$ we
get exactly the same form of sum considered in Lemma \ref{lem:arithmsum}. (Since we estimate
via the sum of moduli, the phase factors are irrelevant.)

The verification of Hypothesis 2 should now be clear. 

\begin{itemize}

\item For fixed $k$, and
hence $\epsilon$, we choose a covering of an appropriate annular region
around $\Gamma$ from the covering in $\bR^{4}$ constructed in
\S \ref{ss:sum-est}, adapted to the quotient
as above. Along with this covering we get a collection of approximately
holomorphic sections, multiplying the sections of Section~\ref{sec:construct}
by cut-off functions to extend over the $4$--manifold. There is just one very
small point to mention. In the covering constructed in 
\S \ref{ss:sum-est} some of the centres are taken to lie on
the $x_{0}$--axis, where the co-ordinate $H$ vanishes. On the other hand,
when we defined the sections $s_{\ux',t'}$ we ruled out this case. However
this is a completely artificial problem and we merely need to take sections
associated to points arbitrarily close to the axis. 

\item We extend this covering to the remainder of the $4$--manifold using
the familiar approximately holomorphic co-ordinates. Likewise for each ball
in the covering we have approximately holomorphic sections, defined just
as in the theory for compact symplectic manifolds. 

\item The localisation properties of the sections, expressed through the
convergence of the sums in the last two items of Hypothesis 2, follow from
the estimates in Section \ref{sec:estimates}. 
\end{itemize}

\subsection{The local model, verification of Hypothesis 3}
\label{ss:localmodel}

In this subsection we will  construct  
sections $\sigma_{0}, \sigma_{1}$ satisfying Hypothesis 3. The construction
is completely explicit but is reasonably complicated so we will perform
it in four stages.

\rk{Stage I}
     
     Consider the Riemann surface $\bC/2\pi i\bZ$ with the symplectic form
     $dx\wedge dy$, where $z=x+iy$ is the standard co-ordinate on $\bC$. Let
     $\mathbb{L}$ be the Hermitian holomorphic line bundle over $\bC/2\pi
     i\bZ$
     with a connection having curvature $-i\,dx\wedge dy$ and with holonomy $-1$
      around the
     circle $C$ corresponding to the imaginary axis. 
     
     \begin{lem}
     There are holomorphic sections $\theta_{0},\theta_{1}$ of $\,\mathbb{L}$
      such that
     \begin{itemize}
     \item The $\theta_{i}$ are bounded.
     \item The sections $\theta_{0},\theta_{1}$ have no common zeros and
     the map $$f^I=\theta_{1}/\theta_{0}\co \bC/2\pi i\bZ\rightarrow \bCP^{1}$$
     maps the circle $C$ bijectively to the circle $i\bR\cup \{\infty\}$
     in $\bCP^{1}$.
     \item The derivative $\partial f^I$ is $\lambda$--transverse to $0$ for
     some $\lambda>0$.\end{itemize}
     \end{lem}
     
These sections can be constructed as follows. Recall that the Weierstrass
$\wp$--function of the rectangular lattice $\Lambda=2\bZ\oplus 2\pi i \bZ$
is an even meromorphic function on the elliptic curve $\bC/\Lambda$ with a
double pole at the origin, representing it as a double cover of $\bCP^1$
ramified at $p_0=0$, $p_1=1$, $p_2=i\pi$ and $p_3=1+i\pi$. The meromorphic
function $\wp$ is the
quotient of two holomorphic sections of the line
bundle $\mathcal{O}(2p_0)$ over $\bC/\Lambda$. Since
$\wp(z)$ and $\wp(1-z)$ have the same ramification points, they
must differ by an automorphism of $\bCP^1$ (this also follows from the
fact that $\mathcal{O}(2p_0)$ and $\mathcal{O}(2p_1)$ are isomorphic).
Setting $a=\wp(1)$ and $b=\wp(\frac{1}{2})^2-
2\wp(\frac{1}{2})\wp(1)$, we have
$$\wp(1-z)=\frac{a\wp(z)+b}{\wp(z)-a}.$$
The line $\mathrm{Re}(z)=\frac{1}{2}$ is one of the two components of
the fixed point locus of the antiholomorphic involution
$z\mapsto 1-\overline{z}$ of $\bC/\Lambda$, and is mapped bijectively
by $\wp$ to the fixed point locus $\Theta$ of the involution
$$w\mapsto \frac{a\overline{w}+b}{\overline{w}-a}.$$
Choose a fractional linear transformation $\varphi\in\mathrm{Aut}(\bCP^1)$
mapping the circle $\Theta$ to the imaginary axis $i\bR\cup\{\infty\}$, and let
$$f^I(z)=\varphi(\wp(z+\frac{1}{2})).$$
Then $f^{I}$ is a doubly-periodic meromorphic function
which maps the imaginary axis to itself, without ramification. 
We can write $f^{I}$ as the quotient $f^{I}= \theta_{1}/\theta_{0}$
of two holomorphic sections of the line bundle $\mathcal{O}(2p')$ over
$\bC/\Lambda$,  where $p'=-\frac{1}{2}$. This degree 2 line bundle 
 can easily be seen to admit a holomorphic
connection with curvature $-i\,dx\wedge dy$ and holonomy $-1$ around the
   circle corresponding to the imaginary axis.

     Now recall that in our standard model around a component of $\Gamma$
     we write our  line bundle ${\cal L}$ as ${\cal L}_{1}\otimes {\cal
     L}_{2}$, where ${\cal L}_{1}$
     has curvature $-idH\wedge d\theta $ and ${\cal L}_{2}$ has curvature 
     $-idQ\wedge dt$. Writing
     $z=\epsilon^{-1} Q + i \epsilon t $, the pullback of ${\cal L}_{2}$
     descends to the quotient $\mathbb{C}/\Lambda$, where it can be identified
     with $\mathbb{L}$. Here we use the condition that the holonomy
     around each component of $\Gamma$ is $-1$. Thus we can regard 
     $\theta_{0}$ and $\theta_{1}$ as sections of ${\cal L}_{2}$.  Then define
     $$\sigma^{I}_{0}= \theta_{0} \otimes \sigma \ , \ \sigma^{I}_{1}=\theta_{1}\otimes
     \sigma,$$ 
     where $\sigma$ is the section of ${\cal L}_{1}$ constructed in Section
     \ref{sec:construct} above. 
     
     These sections $\sigma^{I}_{0},\sigma^{I}_{1}$ have some of the properties
     required by Hypothesis~3. Let $z_{r}\in \bC$ be the branch points
     of $f^I$. We can choose disjoint discs in $\bC$ of a fixed radius $\delta$ centred
     on the $z_{r}$. We also suppose that $\delta$ is chosen small enough
     that $\vert \mathrm{Re}(z_{r})\vert >2\delta$ for all $r$. Let $N_{r}$ be the tubular region in $\bR^{4}$ defined
     by the condition
     $\vert z-z_{r}\vert \leq \delta$. Then the sections have all the desired
     properties outside the region
     $$  (\bigcup N_{r}) \cap (X\setminus K), $$
     where we recall that $K$ is the set defined by $\vert \ux\vert\geq 10$.
     In the following stages we will modify the sections to achieve all
     the required properties. (In fact, except for the very last step,
     the modifications will only involve the \lq\lq numerator'' $\sigma_{1}^{I}$.)

\rk{Stage II}

In the second stage we improve the sections over the intersection of the
 tubular regions $N_{r}$ with the annulus $\Omega= \{ 2 <\vert \ux\vert<5\}$.
 We take a standard cut-off function $\beta$ supported in $[0,\delta)$ and equal
 to $1$ on $[0,\delta/2]$. Then define $\beta_{r}=\beta(\vert z-z_{r}\vert)$.
 Thus $\beta_{r}$ is supported in the tube $N_{r}$ and equal to $1$ on
 a half-sized tube. Recall that we have functions $F_{+},F_{-}$ which are
 holomorphic along the quadric surfaces $z=\rm{const}$. In Section
 \ref{sec:sec3} these were only defined over the subsets $G^\pm$, but we
 now extend them by zero over the complement of $G^\pm$. For a small parameter
 $\alpha$, to be chosen later,  we set:
 \begin{equation}
  \sigma^{II}_{0} = \sigma^{I}_{0},\qquad \sigma^{II}_{1}= \sigma^{I}_{1}
 + \alpha \sum_{r} \beta_{r} (F_{+}+ F_{-}) \sigma^{I}_{0}.
 \label{eq:sigmaII}\end{equation}
 Thus $  f^{II}= \sigma^{II}_{1}/\sigma^{II}_{0}$ is
 $$  f^{II} = f^I + \alpha \sum_{r} \beta_{r} (F_{+}+ F_{-}). $$
 
 \begin{lem}\label{lem:stage2}
 For sufficiently small $\alpha, \epsilon$ there  
 are $\kappa_{1},\kappa_{2},\kappa_{3}>0$ such
 that, over $\Omega$, 
 \begin{itemize}
 \item $\partial f^{II}$ is $\kappa_{1}$--transverse to $0$;
 \item $\vert \db f^{II}\vert\leq \mathrm{max} (\epsilon \kappa_{2}, \vert \partial
 f^{II}\vert - \kappa_{3}).$ \end{itemize}
 \end{lem}

 \noindent
 There are positive constants, independent of $\epsilon$, so that over
 $\Omega$
 \begin{itemize}
 \item $\vert \nabla \beta_{r}\vert \leq k_{1}$
 \item $\vert F_{+}+F_{-}\vert \leq k_{2}$
 \item $\vert \nabla(F_{+}+F_{-})\vert\leq k_{3}$
 \item $\vert \nabla_{z}(F_{+}+F_{-})\vert\leq k_{4}\epsilon$.\end{itemize}
 Here we write $\nabla_{z}$ for the component of the derivative in the
 $z$ direction. The existence of these bounds is fairly clear, there is
 just one point we want to spell out here. The function $F_{+}$
 is not smooth along the part of the null cone where $x_{0}<0$, but behaves
 like $(-Q)^{\nu}$ when $Q<0$ and vanishes when $Q\geq 0$, where $\nu$ is
 $\sqrt{3/2}$. Since $\nu>1$ we have a uniform bound on the first derivative, but
 one might worry about the higher derivatives. 
 In terms of $x=\mathrm{Re}(z)= Q/\epsilon$, $F_{+}$ behaves like $\epsilon^{\nu}(-x)^{\nu}$;
 so on the set where $x<-\delta$ all derivatives with respect to $z$
 are bounded by  multiples of $\epsilon^{\nu}$. Since our formulae only
 involve $F_{+}$ over the tubes $N_{r}$, on which $\vert x \vert >\delta$,
 we do not encounter any problems from the singularities of $F_{\pm}$.
  
 Let $M_{r}\subset N_{r}$ be the interior tube on which $\beta_{r}=1$.
 Then there is a $K_{1}>0$ such that $\vert \partial f^{I}\vert \geq K_{1}$
 outside the $M_{r}$ (but inside $\Omega$). So on this set
 $$  \vert \partial f^{II}\vert \geq \vert \partial f^{I} \vert- \alpha\,
 \Bigl\vert \sum \nabla \beta_{r}(F_{+}+F_{-}) + \beta_{r} \nabla(F_{+}+F_{-})\Bigr\vert.$$
 At any given point there is at most one term contributing to the sum (since
 the $N_{r}$ are disjoint) so we have
 $$  \vert \partial f^{II} \vert \geq K_{1}- \alpha(k_{1}k_{2} + k_{3}).
 $$
 Thus if we choose $\alpha< K_{1}/(10(k_{1}k_{2}+k_{3}))$ we have
 $ \vert \partial f^{II}\vert \geq 9K_{1}/10$ outside the $M_{r}$. On the
 other hand, outside the $M_{r}$, we have
 $$ \vert \db f^{II}\vert \leq \alpha(k_{1}k_{2}+ k_{3})\leq K_{1}/10.
 $$ 
 Now consider the situation inside a tube $M_{r}$ where
 $$  f^{II}= f^{I} + \alpha(F_{+}+F_{-}). $$
 Then $$\vert \db f^{II}\vert =\alpha \vert \db (F_{+}+F_{-})\vert\leq \alpha
 k_{4}\epsilon, $$
 since $F_{+}+F_{-}$ is holomorphic along the quadric surfaces and only
 the $z$ derivative contributes. Now on each quadric surface the holomorphic
 function $F_{+}+F_{-}$ is either unramified (for $Q>0$) or has two
 ramification points (where $x_0=0$ and
 $\theta\in\{0,\pi\}$, for $Q<0$). Let $p_{r}^\pm$ be
 the ramification points on the surface corresponding to $z_{r}$ and $B^\pm_{r}$
 be the $\delta$--balls about $p_{r}^\pm$. It is clear
 then that there is a $K_{2}>0$ such that in the intersection of $\Omega$
 and $M_{r}\setminus (B^+_{r}\cup B^-_{r})$, and once $\epsilon$ is sufficiently small,
 we have 
 $$   \vert \partial_{w}( F_{+}+F_{-}) \vert\geq K_{2}, $$
 where $\partial_{w}$ denotes the derivative along the quadric surfaces.
 Thus, on this set,
   $$  \vert \partial f^{II}\vert \geq \alpha K_{2}. $$
  On the other hand it is also clear that if $B^\pm_{r}$ meets the annulus
  $\Omega$ we have a  bound on the inverse of the
   Hessian of $f^{II}$ over $B^\pm_{r}$:
   $$\vert (\nabla \partial f^{II})^{-1} \vert \leq K_{3} \alpha^{-1}.$$   
   In sum then, $\partial f^{II}$ is $\kappa_{1}$--transverse to $0$ over
   $\Omega$ with
            $$\kappa_{1}= \min({\textstyle \frac{9}{10}}K_{1},
              \,\alpha K_{2}, \,\alpha K_{3}^{-1}),
            $$
            while
            $$\vert \db f^{II}\vert \leq \max(\kappa_{2} \epsilon,\, \vert
            \partial f^{II}\vert - \kappa_{3}) $$
            with $\kappa_{2}=\alpha k_{4}$, $\kappa_{3}= 8K_{1}/10$.

\rk{Stage III}
 
 The formulae (\ref{eq:sigmaII}) define $\sigma^{II}_{i}$ over all of $\bR^{4}$ but they
 do not satisfy the requirements of Hypothesis 3. One problem is that the
 section $\sigma^{II}_{1}$  is
  not $\epsilon$--holomorphic over $\{\vert x\vert \geq 10\}$
  because when we differentiate we pick up a term from $\nabla \beta_{r}$
  which is multiplied by the small parameter $\alpha$ but is not controlled
  by $\epsilon$. We now get over this problem.
  
  First we address  the fact that the functions $F_{+}, F_{-}$ are not
  smooth along the null cone. This is similar to the construction in \S \ref{ss:secl2}. 
  We define a function $\gamma_{+}$ in the
  region $\vert \ux\vert>0.5$ in the following way. We let $\gamma_{+}(\ux)=1$
  if $x_{0}>0$ and $\gamma_{+}(\ux)= \gamma_{\epsilon}(Q(\ux))$ if $x_{0}\leq
  0$, where $\gamma_{\epsilon}$ is a standard cut-off function, with
  $\gamma_{\epsilon}(t)= 1$ if $t\leq -\delta \epsilon$ and 
  $\gamma_{\epsilon}(t)=  0$ if $t\geq -\frac{1}{2}\delta\epsilon$. Once $\epsilon$ is sufficiently
  small, the function $\gamma_{+}$ is smooth in $\{\vert \ux\vert>0.5\}$. Now
  we put $\tilde{F}_{+}= \gamma_{+} F_{+}$. Then $\tilde{F}_{+}$ is a smooth
  function over $\{\vert \ux\vert>0.5\}$, holomorphic along the quadric surfaces.
   We define $\tilde{F}_{-}$ in a symmetrical fashion. Notice that 
  $\tilde{F}_\pm=F_\pm$ over $\bigcup N_r$.
   
    Let $\chi=\chi(\vert \ux\vert)$ be a
  standard cut-off function, equal to $1$ when $\vert \ux\vert \leq 5$ and
  zero when $\vert \ux\vert\geq 10$. Now we set
  $ \sigma^{III}_{0}= \sigma^{II}_{0}= \sigma^{I}_{0}$
  and
  $$ \sigma^{III}_{1}= 
  \chi \sigma^{II}_{1} + (1-\chi)(\sigma^I_1+\alpha(\tilde{F}_{+}+\tilde{F}_{-})
  \sigma^I_0).
  $$    
  These sections are well-defined everywhere, even though the $\tilde{F}_{\pm}$
  are not, because the factor $(1-\chi)$ vanishes when $\vert \ux\vert\leq
  0.5$.
  \begin{lem}
  There  are  constants $C,\kappa_{1},\kappa_{2},\kappa_{3}$ such that for small enough $\alpha$ and
  $\epsilon$ we have \begin{itemize}
  \item
  $ \vert \db \sigma_{i}^{III}\vert \leq C\epsilon $ in $\{\vert \ux\vert \geq
  10\}$
  \item if $f^{III}=\sigma^{III}_{1}/\sigma^{III}_{0}$ then
  over $\{2 \leq \vert \ux\vert\leq 10\}$, $\partial f^{III}$ is $\kappa_{1}$--transverse
  to $0$ and $\vert \db f^{III}\vert \leq \max(\epsilon \kappa_{2}, \vert \partial
  f^{III}\vert - \kappa_{3})$\end{itemize}
  \end{lem} Consider the second item of the Lemma.                    
    The proof of the previous Lemma applies equally well to any fixed annulus,
    with suitable adjustment of constants. Thus here we have to deal with
    extra terms introduced by, on the one hand,
     the passage from $F_{\pm}$ to $\tilde{F}_{\pm}$
     and on the other hand the introduction of the cut-off function $\chi$.
     The first issue is essentially covered by the
     discussion at the beginning of the proof of Lemma \ref{lem:stage2},
    which applies equally well to $\tilde{F}_{\pm}$.
    So we will simply ignore the distinction between $\tilde{F}_{\pm}
    $ and $F_{\pm}$, and consider the function
    $$  f^{I}+ \alpha\,(\chi\, {\textstyle \sum\beta_{r}} + (1-\chi)) (F_{+}+F_{-}). $$
    When we differentiate this we get a new term 
    $$\alpha\, \nabla\chi\, ({\textstyle \sum\beta_{r}}-1) (F_{+}+ F_{-})$$
    which is supported outside the $M_{r}$. The size of $\nabla \chi$ is
    bounded (independently of $\epsilon$): 
    $ \vert \nabla \chi\vert \leq k_{5}$ say. Then the size of the
    new term is bounded by $\alpha k_{2} k_{5}$. Thus the estimates inside
    $M_{r}$ are completely unchanged and outside $M_{r}$ we have
$$ \vert \partial f^{III} \vert \geq {\textstyle \frac{9}{10}}K_{1} - 
\alpha k_{2}k_{5}, \ \ \vert
\db f^{III}\vert \leq {\textstyle \frac{1}{10}}K_{1} + \alpha k_{2} k_{5}. $$
This establishes the second item of the lemma, once $\alpha$ is sufficiently
small and the constants $\kappa_{i}$ are adjusted suitably.

The first item of the lemma follows from the fact that on $\{\vert \ux\vert
\geq 10\}$ we have simply
   $$f^{III}= f^{I}+ \alpha (\tilde{F}_{+} + \tilde{F}_{-}) $$
   and we can apply the bounds on the derivatives of $\tilde{F}_{\pm}$,
   together with the rapid exponential decay of $\sigma$.

\rk{Stage IV}

In this final stage, we modify the construction to ensure that we get a
topological Lefschetz fibration over the inner region. For each point
$z$ in one of the discs $\vert z-z_{r}\vert <\delta$ we have a corresponding quadric
surface $\Sigma(z)$, say. We can use our standard co-ordinates $H,\theta$ to identify these
surfaces for different values of $z$, so we have diffeomorphisms $\tau_{z}:
\Sigma(z)\rightarrow \Sigma(z_{r})$. Let $\rho$ be a standard cut-off function
with $\rho(\ux)=0$ if $\vert \ux\vert\leq 1$ and $\rho(\ux)=1$ if $\vert \ux\vert\geq
2$. On the surface $\Sigma(z)$ we define
      $$  F_{\pm,r} = \rho F_{\pm} + (1-\rho) F_{\pm}\circ \tau_{z}. $$
      This defines new functions $F_{\pm,r}$ on the tube $N_{r}$ which
      are equal to $F_{\pm}$ when $\vert \ux\vert\geq 2$. Now define
      $\sigma^{IV}_{1}$ to be equal to
       $ \sigma^{III}_{1}$ in $\vert \ux \vert \geq 2$ and to be given
       by the modified formula
       $$  \sigma^{IV}_{1} = \sigma^{I}_{1} + \alpha \sum \beta_{r} (F_{+,r}
       + F_{-,r}) \sigma_{0}^I$$
       in the inner region $\vert \ux\vert \leq 2$.
 Again, we keep the same \lq\lq denominator'' $\sigma_{0}^{IV}=\sigma_{0}^{III}$.
       
       \begin{lem} When $\alpha$ is sufficiently small the ratio
       $ f^{IV}= \sigma^{IV}_{1}/\sigma^{IV}_{0}$ is a topological Lefschetz
       fibration over $\{\vert \ux\vert \leq 1\}$, with symplectic fibres. \end{lem}
       
       Notice that the statement of this lemma does not involve any almost
       complex structure or quantative estimates. Clearly the only issue
       involves the behaviour over the tubes $N_{r}$ and to prove the Lemma
       we consider an auxiliary almost-complex structure on the tubes~--~just
       the integrable product structure
       given by the identification with $D_{r}\times \Sigma(z_{r})$. Thus
       taking $w$ as a complex co-ordinate on $\Sigma(z_{r})$ our function
       has the simple form on $N_{r}$,
       $$  f^{IV}(z,w)= f^{I}(z) + \alpha\,\beta_{r}(z)\, g(w), $$
       where $g$ is the holomorphic function $F^{+}+F^{-}$ on $\Sigma(z_{r})$.
       This function $f^{IV}$ is holomorphic on the interior tube $M_{r}$ with nondegenerate
       critical points. So it suffices to show that
          $$ \vert \db f^{IV}\vert < \vert \partial f^{IV}\vert $$
          on $N_{r}\setminus M_{r}$, where now $\db, \partial$ refer to
          the product complex structure. Then on this region we still  have
          $$  \vert \partial f^{I}\vert \geq   K_{1} ,\  \  \vert \nabla \beta_{r}\vert
          \leq k_{1},\ \ \vert F_{+}+ F_{-} \vert \leq k_{1}. $$
          Now
          $$ \vert \db f^{IV} \vert  = \alpha \vert \nabla \beta_r\vert
          \,\vert g \vert \leq \alpha k_{1} k_{2}, $$
          while $$\vert \partial f^{IV} \vert \geq \vert \partial_{z} f^{IV}
          \vert = \vert \partial_{z} f^{I} + \alpha \,\nabla \beta_r\, g \vert
          \geq K_{1}-\alpha k_{1} k_{2}. $$
          Thus the result follows once $\alpha < K_{1}/2k_{1}k_{2}$. (The
          point of this proof is that we do not need to control the derivatives
          of $F_{\pm}$ in the inner region where $\vert \ux \vert<1$.)

       This essentially completes our construction. There is just one last
       issue; that we want to have sections defined over the whole manifold
       $X$ while up to now we have been working in the local model.
       So we define
       $$ \sigma_{i} = \phi \,\sigma^{IV}_{i} $$
       where $\phi$ is a cut-off function equal to $1$ for $\vert \ux \vert
       \leq c\,\epsilon^{-1}$ and to zero when $\vert \ux\vert \geq 2c\,\epsilon^{-1}$
       (for some fixed $c>0$). Thus these sections
       $\sigma_{i}$ can be extended by $0$ over the whole of $X$.
      Our final result is:
              \begin{prop} There are constants $\kappa_{1}, \kappa_{2},
              \kappa_{3}, \uC$ such that for a suitable value of $\alpha$,
              and once $\epsilon$ is  sufficiently small, the sections $\sigma_{0}, \sigma_{1}$ satisfy Hypothesis
              $H_{3}(\epsilon, \kappa_{1}, \kappa_{2}, \kappa_{3}, \uC)$.
              \end{prop}

        The proof of this proposition has been largely covered in the preceding
        lemmas. There is one point left over from Stage IV: we need to
        check that the map $f^{IV}$ satisfies the required transversality
         estimates over the annulus $\{1\leq \vert \ux \vert\leq 2\}$. Here
         the discussion follows the same lines as in Stage II, except that
         we  replace the functions $F_{\pm}$ by the linear combinations
         $$ F_{\pm,r} = \rho F_{\pm} + (1-\rho) F_{\pm}\circ \tau_{z} =
         F_{\pm} + (1-\rho) (F_\pm\circ\tau_z-F_\pm).$$
         But $F_{\pm,r}-F_{\pm}$ is $O(\epsilon)$ (along with its derivatives).
         So the extra term introduced here causes no problem.

        Finally, we check that the terms introduced by the cut-off function
        $\phi$ are much smaller than $\epsilon$ due to the rapid decay of
        $\sigma$ away from the origin; this completes the argument.
        
\subsection{The odd case}\label{ss:oddcase}

In all our discussion so far we have focussed on the case when the zero
set $\Gamma$ has just one component and the local model is the \lq\lq
even'' version $N_{+}$. We now consider the modifications required
for the general case. It is quite obvious that the existence of several
components makes no difference to the argument, all we have to discuss
is the case of the \lq\lq odd'' model $N_{-}$. In this case the map
$\sigma_{-}$ interchanges the two components of the positive cone and maps
$(H,\theta)$ to $(-H, -\theta)$ but preserves the co-ordinate $Q$. 

We begin with the last part of the construction, the local model
 in Section \ref{ss:localmodel} above. Since the sections $\sigma_{0}^{I}, \sigma_{1}^{I}$
 only depend on the $Q,t$ variables the first step goes through unchanged.
 In the later stages we use the fact that the involution interchanges the
 functions $F_{+}$ and $F_{-}$, and so preserves their sum. The upshot is
 that the whole construction in \S \ref{ss:localmodel} goes over immediately to the odd
 case. 
 
 The slightly more substantial discussion involves the construction of the localised
 sections in the odd case. Working in $\bR^{4}$, in \S \ref{ss:hyp2} we have  
 defined $2\pi\epsilon^{-1}$--periodic
 sections $\Theta_{\ux',t'}$. We write these as 
 $$  \Theta_{\ux',t'}= \Theta^{+}_{\ux',t'} + \Theta^{-}_{\ux',t'},$$
 taking the even and odd terms respectively in the sum (\ref{eq:deftheta}).
 Thus $\Theta^{\pm}_{\ux',t'}$
 are $4\pi\epsilon^{-1}$--periodic and the translation $t\mapsto t+2\pi\epsilon^{-1}$
 interchanges the two sections. Now we define
 $$  s_{\ux',t'} = \Theta^{+}_{\ux',t'} \otimes 
  \tau^{*}_{\ux'}+ \Theta^{-}_{\ux',t'}\otimes \tau^{*}_{\sigma_{-}(\ux')}. $$
 These sections are invariant under the map $\overline{\sigma}_-$ on $\bR^{4}$
 so descend to $N_{-}$. 

\section{The converse result} \label{sec:converse}

\subsection{Proof of Theorem \ref{thm:converse}}

 The proof of Theorem \ref{thm:converse} is very similar to that of
 Gompf's result for symplectic Lefschetz fibrations and pencils \cite{Go2},
 which in turn relies on a classical argument of Thurston \cite{Th}.

 Let $X$ be a compact oriented 4--manifold, and let $f\co X\setminus A\to S^2$
 be a singular Lefschetz pencil with singular set $\Gamma$. Let $B$ be the
 finite set of isolated critical points of $f$ in $X\setminus \Gamma$, near
 which $f$ is modelled on $(z_1,z_2)\mapsto z_1^2+z_2^2$. We assume
 that there exists a cohomology class $h\in H^2(X)$ such that $h(\Sigma)>0$
 for every component $\Sigma$ of a fibre of $f$ (if every component $\Sigma$
 contains a base point of the pencil, then we can choose $h$ to be
 Poincar\'e dual to the homology class of the fibre).
 \medskip

 {\bf Step 1}\qua We start by constructing a closed 2--form $\omega_0$ over a regular
 neighbourhood $U$
 of $A\cup B\cup \Gamma$, non-degenerate outside of $\Gamma$ and positive
 on the fibres of $f$, in the following manner. Near $A\cup B$, we take
 $\omega_0$ to be the standard K\"ahler form of $\bC^2$ in some local
 oriented co-ordinates in which $f$ is given by the standard models
 $(z_1,z_2)\mapsto z_1/z_2$ and $(z_1,z_2)\mapsto z_1^2+z_2^2$. Near
 a point $p\in\Gamma$, we have oriented local co-ordinates in which $f$
 is modelled on $(x_0,x_1,x_2,t)\mapsto x_0^2- \frac{1}{2}(x_1^2+x_2^2)+it$.
 Then we let  $$\omega_p=d\Bigl(\chi(|t|)\,x_0(x_1\,dx_2-x_2\,dx_1)\Bigr),$$
 where $\chi$ is a suitable smooth cut-off function, and we
 extend $\omega_p$ into a closed 2--form defined over a tubular neighbourhood
 of the component of $\Gamma$ containing $p$, supported near $p$.
 The 2--form $\omega_p$ vanishes on $\Gamma$, and its restriction to the
 fibres of $f$ is non-negative, and positive near $p$ (outside of $\Gamma$).
 By choosing a suitable finite subset $\{p_i\}$ of $\Gamma$ and
 setting $\omega_0=\sum_{i} \omega_{p_i}+f^*(\omega_{S^2})$, we obtain a
 closed 2--form defined over a neighbourhood of $\Gamma$, positive on the
 fibres, vanishing on $\Gamma$ and non-degenerate outside of $\Gamma$.
 \medskip

 {\bf Step 2}\qua Our next task is to construct local closed 2--forms over
 neighbourhoods of the fibres of $f$, compatible with our local model
 $\omega_0$ near $A\cup B\cup \Gamma$, and restricting positively to the
 vertical tangent spaces; we will then glue these into
 a globally defined 2--form. For this purpose, we choose a closed 2--form
 $\eta\in\Omega^2(X)$, with $[\eta]=h$.
 Since $U$ retracts onto a union of points and circles, $H^2(U)=0$, and
 there is a 1--form $\beta$ such that $\omega_0-\eta=d\beta$ over $U$.
 Extending $\beta$ to an arbitrary 1--form on $M$ with support in a
 neighbourhood of $U$, and replacing $\eta$ by $\eta+d\beta$, we can assume
 that $\eta_{|U}=\omega_0$. 

 Given any point $q\in S^2$, we can find a regular neighbourhood $V_q$ of the
 fibre $F_q=f^{-1}(q)\cup A$, and neighbourhoods $U''\subset U'\subset U$ of
 $A\cup B\cup \Gamma$, with the following properties:
 \begin{itemize}
 \item $V_q\cap U'$ retracts onto $F_q\cap (A\cup B\cup\Gamma)$;
 \item $V_q\setminus (V_q\cap U'')$ is diffeomorphic to a product 
 $D^2\times (F_q\setminus (F_q\cap U''))$;
 \item there exists
 a smooth map $\pi\co V_q\to V_q$ with image in $F_q\cup (V_q\cap U')$, equal to identity
 over $F_q\cup (V_q\cap U'')$.
 \end{itemize}
 The first and second properties can easily be ensured by shrinking $V_q$ so
 that all critical points of $f$ over $V_q$ lie close to the singular locus
 of $F_q$; the map $\pi$ can then be built by interpolating between the
 identity map over $V_q\cap U'$ and the projection map from $V_q\setminus
 (V_q\cap U'')$ to $F_q\setminus (F_q\cap U'')$ given by the product
structure.
 
 Since by assumption $[\eta]=h$
 evaluates positively over each component of $F_q$, shrinking $U'$ if
 necessary we can equip $F_q$
 with a (near) symplectic form $\sigma_q$ which coincides with $\eta$
 over $F_q\cap U'$, is symplectic
 over the smooth part of $F_q$, and such that $[\sigma_q-\eta_{|F_q}]=0$
 in $H^2(F_q, F_q\cap U')$ (ie, 
 $\int_\Sigma \sigma_q=h(\Sigma)$ for every component $\Sigma$ of $F_q$).
 Using the projection $\pi$ to pull back the 2--forms $\eta$ on $V_q\cap U'$
 and $\sigma_q$ on $F_q$, we obtain a 2--form $\tilde\eta_q$ on $V_q$
 with the following properties:
 \begin{itemize}
 \item $\tilde\eta_q$ is closed, and $[\tilde\eta_q]=h_{|V_q}$;
 \item $\tilde\eta_q$ coincides with $\eta$ over $V_q\cap U''$;
 \item $[\tilde\eta_q-\eta]=0$ in $H^2(V_q, V_q\cap U'')\simeq 
 H^2(F_q,F_q\cap U'')$;
 \item (shrinking $V_q$ if necessary)
 the restriction of  $\tilde\eta_q$ to $\mathrm{Ker}(df)$ is
 positive at every regular point of $f$ in $V_q$.
 \end{itemize}
 By the third property, there is a 1--form $\beta_q$ on $V_q$,
 vanishing identically over $V_q\cap U''$, such that
 $\tilde\eta_q=\eta+d\beta_q$.
 \medskip
 
 {\bf Step 3}\qua For each $q\in S^2$, the above construction yields a 2--form
 $\tilde\eta_q$ defined over a neighbourhood $V_q$ of the fibre $F_q$. By
 compactness, each $V_q$ contains the preimage of a neighbourhood
 $D_q$ of $q$ in $S^2$, and there is a finite set $Q\subset S^2$ such that the
 open subsets $(D_q)_{q\in Q}$ cover $S^2$. Consider a smooth partition of
 unity $\sum_{q\in Q} \rho_q=1$ with $\rho_q$ supported inside $D_q$,
 and define
  \begin{equation}\label{eq:tildeeta}
  \tilde\eta=\eta + d\Bigl(\,\textstyle\sum\limits_{q\in Q} 
 (\rho_q\circ f)\,\beta_q\Bigr).\end{equation}
 The closed 2--form $\tilde\eta$ coincides with $\eta$ over the intersection
 $\tilde{U}$ of the neighbourhoods $U''$ considered
 above for all $q\in Q$, and hence is well-defined over all $X$ even though
 (\ref{eq:tildeeta}) only makes sense outside of $A$.
 Moreover, the restriction of $\tilde\eta$ to a fibre $F_p$ of $f$ is
 $$\tilde\eta_{|F_p}=\textstyle \sum\limits_{q\in Q} \rho_q(f(p))\,
 (\eta+d\beta_q)_{|F_p}=\sum\limits_{q\in Q} \rho_q(f(p))\,\tilde\eta_{q|F_p},$$
 ie, a convex combination of positive forms; hence $\tilde\eta$ induces
 a symplectic structure on each fibre of $f$ (outside of the critical points).
 Hence, as in Thurston's original argument, for large enough $\lambda>0$ the
 2--form $$\omega_\lambda=\tilde\eta+\lambda\,f^*\omega_{S^2}$$ is closed and non-degenerate
 over $X\setminus (A\cup\Gamma)$, and restricts positively to the fibres of
 $f$; moreover $\omega_\lambda$ vanishes transversely along $\Gamma$, as expected.
 However, $\omega_\lambda$ does not extend smoothly over the base locus $A$, and
 we need to apply a trick due to Gompf \cite{Go2} in order to complete the
 construction.
 \medskip
 
 {\bf Step 4}\qua Near a base point of $f$, consider local co-ordinates in
 which $f$ is the projectivisation map from $\bC^2\setminus\{0\}$ to
 $\mathbb{CP}^1$, and denote by $r$ the radial co-ordinate and by
 $\alpha$ the pullback to $\bC^2\setminus\{0\}=\bR^+\times S^3$ of the
 standard contact form of $S^3$. Then we have 
 $$\omega_\lambda=\lambda\,f^*\omega_{S^2}+\omega_0=(\lambda+r^2)\,f^*\omega_{S^2}
 +\frac{1}{2}\,d(r^2)\wedge\alpha.$$
 Setting $R^2=\lambda+r^2$, we have
 $\omega_\lambda=R^2\,f^*\omega_{S^2}+\frac{1}{2}\,d(R^2)\wedge \alpha$.
 Hence, the radially symmetric map  $\varphi(z)=(\lambda+|z|^2)^{1/2}\,z/|z|$
 defines a symplectic embedding of $(\bC^2\setminus\{0\},\omega_\lambda)$ into
 $(\bC^2,\omega_0)$, whose image is the complement of a ball of radius
 $\lambda^{1/2}$.
 Therefore, by replacing the ball of radius $\epsilon$ around each point
 of $A$ in $(X,\omega_\lambda)$ by a standard ball of radius
 $(\lambda+\epsilon^2)^{1/2}$ in $(\bC^2,\omega_0)$ we can obtain a globally
 defined near-symplectic structure $\omega$. More precisely, $\omega$
 is naturally defined on the 4--manifold $Y$ obtained from $X$ by this
 cut-and-paste process; however $Y$ can easily be identified with $X$
 via a diffeomorphism which equals identity outside of an arbitrarily
 small neighbourhood of $A$. 
 
 Another viewpoint is to observe
 that $\omega_\lambda$ extends smoothly to the manifold $\hat{X}$ obtained
 by blowing up $X$ at the base points; gluing in standard balls in place
 of the exceptional divisors amounts to a symplectic blowdown of $(\hat{X},
 \omega_\lambda)$, and yields a well-defined near-symplectic form on $X$.
 In any case, one easily checks that the various requirements
 satisfied by $\omega_\lambda$ (vanishing along $\Gamma$, and positivity over
 the fibres of $f$) still hold for the modified form $\omega$; this
 completes the main part of the argument.

 The cohomology class of the constructed form $\omega$ is
 $h+\lambda\,f^*[\omega_{S^2}]$ (identifying implicitly $H^2(X)$ with
 $H^2(X\setminus A)$). If we assume that every component of
 every fibre contains a base point we can take $h$ to be Poincar\'e dual
 to the class of the fibre. In that case $f^*[\omega_{S^2}]=h$ (up to a
 scalar factor), so after scaling by $\frac{1}{1+\lambda}$ we can ensure that
 $[\omega]=h$ is Poincar\'e dual to the fibre. (However, since we have
 no control over the relative class $[\omega]\in H^2(X,\Gamma)$,
 deformations of near-symplectic forms in the class $h$ are not always
 generated by isotopies of $X$).

 Before we can state more precisely our uniqueness result for the
 deformation class of $\omega$, we consider again the positivity property
 for the restriction of $\omega$ to the fibres of $f$, and its implications
 for the local structure near a point of $\Gamma$.
 Recall that the first-order variation of $\omega$ at a point $x$ of
 $\Gamma$ yields canonically a linear map $\nabla\!_x\omega\co N\Gamma_x\to
 \Lambda^2T^*X_x$. Restrict locally $f$ to a normal slice $D$
 to $\Gamma$ through $x$ obtained as the preimage of a transverse arc to
 $f(\Gamma)$ through $f(x)$. Then the 2--jet of $f_{|D}$ at $x$ defines a
 non-degenerate quadratic form $Q$ on $TD_x\simeq N\Gamma_x$; and, if one
 approaches $x$ in the direction of a non-zero vector $v\in TD_x$, the plane
 field $\mathrm{Ker}\,df$ converges to $v^\bot=
 \mathrm{Ker}\,Q(v,\cdot) \subset TD_x$. Hence, the
 positivity condition on the restriction of the near-symplectic form
 $\omega$ to the fibres of $f$ implies that, for every
 $v\in TD_x\setminus\{0\}$, the 2--form $\nabla\!_x\omega(v)$
 evaluates non-negatively on the 2--dimensional subspace $v^\bot\subset TD_x$
 (since it is the limit of the tangent spaces to the fibres when approaching
 $x$ in the direction of $v$). However, in our case it is easy to check that
 the above construction of $\omega$ guarantees that $\omega_{|\mathrm{Ker}\,df}$ is
 bounded from below by a constant multiple of the distance to
 $\Gamma$ (ie, a constant multiple of the norm of $\omega$). Equivalently,
 the tangent spaces to the fibres near $x$ do not tend to degenerate
 to isotropic subspaces as one approaches $x$. This implies that
 the restriction of $\nabla\!_x\omega(v)$ to the 2--plane $v^\bot$ is in
 fact {\it positive} for every $v\neq 0$. Now we have the following:
 
 \begin{lem} Let $\omega_0,\ \omega_1$ be two near-symplectic forms
 with the same zero set $\Gamma$ and for which the smooth parts of the
 fibres of $f$ are symplectic. Assume moreover that for all $x\in \Gamma$,
 $v\in N\Gamma_x\setminus\{0\}$, $j\in \{0,1\}$, the restriction of 
 $\nabla\!_x\omega_j(v)$ to the limiting tangent plane $v^\bot$ is positive.
 Then $\omega_0$ and $\omega_1$ are deformation equivalent through
 near-symplectic forms with the same properties.
 \end{lem}

 \proof
 Start with the convex combinations $\omega_s=(1-s)\omega_0+s\,\omega_1$.
 For all $s\in [0,1]$, $\omega_s$ is a closed 2--form which vanishes on
 $\Gamma$ and evaluates positively on the fibres of $f$ outside of the
 critical points, but it may be degenerate at some points of $X\setminus
 \Gamma$. We can avoid this problem by deforming $\omega_0$ and $\omega_1$
 to make them standard over a small neighbourhood of $A$, choosing a large
 enough constant $\lambda>0$, and considering the 2--forms $\tilde\omega_s$
 obtained from $\omega_s+\lambda\,f^*\omega_{S^2}$ by inserting standard
 balls near the base points as described above.
 
 The 2--forms $\tilde\omega_s$ are closed and positive on fibres, they vanish
 on $\Gamma$, and if $\lambda$ is large enough they are non-degenerate
 outside of $\Gamma$ (away from $A\cup B$ this follows from Thurston's
 classical argument; and at a point $x\in A\cup B$ this follows from
 positivity on the fibres, which implies that $\tilde\omega_s$
 tames a naturally defined complex structure on $T_x X$). Moreover,
 $\tilde\omega_0$ and $\omega_0$ are deformation equivalent through the
 family of near-symplectic forms obtained by blowing down $\omega_0+t\,
 f^*\omega_{S^2}$ for $t\in [0,\lambda]$; and similarly for $\tilde\omega_1$
 and $\omega_1$. Hence, all that remains to be checked is the non-degeneracy
 of $\nabla\tilde\omega_s$ along $\Gamma$ for all $s\in [0,1]$.
 
 By assumption, for all $x\in \Gamma$ and $v\in N\Gamma_x\setminus\{0\}$,
 the 2--forms $\nabla\!_x\omega_j(v)$ ($j=0,1$), and consequently
 $\nabla\!_x\tilde\omega_j(v)$ too, evaluate positively on the limiting
 vertical tangent space $v^\bot$. Since this positivity condition is
 preserved by convex combinations, we conclude that
 $\nabla\!_x\tilde\omega_s(v)$ evaluates positively on $v^\bot$. Moreover
 this implies that
 $\nabla\!_x\tilde\omega_s(v)\neq 0$ for all $s\in [0,1]$, $x\in \Gamma$, 
 $v\in N\Gamma_x\setminus\{0\}$, which proves that $\tilde\omega_s$
 vanishes transversely along $\Gamma$ and hence is a near-symplectic
 form.

\subsection{Proof of Proposition \ref{prop:ns-sdharm}}

Consider $\bR^4$ with its standard Euclidean structure and orientation,
inducing a splitting $\Lambda^2\bR^4=\Lambda^2_{+,0}\oplus\Lambda^2_{-,0}$.
The wedge-product restricts to a given 3-dimensional subspace $P\subset
\Lambda^2\bR^4$ as a definite positive bilinear form if and only if $P$
can be written as the graph $P=\{\alpha+L(\alpha),\
\alpha\in\Lambda^2_{+,0}\}$ of a linear map $L\co \Lambda^2_{+,0}\to
\Lambda^2_{-,0}$ with operator norm less than $1$.
Therefore, positive definite subspaces
form a ``convex'' subset of the Grassmannian of 3-planes in $\Lambda^2\bR^4$.
Moreover, given an element $\beta\in\Lambda^2\bR^4$ with $\beta\wedge \beta>0$,
the space of all positive definite 3-planes containing $\beta$ is again
convex (and hence contractible).
In another guise, the set of positive definite subspaces can be identified
with the set of conformal classes of Euclidean metrics on $\bR^{4}$, ie,
for each such subspace $P$ there is a unique metric, up to scale, which
realises $P$ as its space of self-dual forms.

Given a near-symplectic form $\omega$ on $X$, our goal is to build a Riemannian
metric with respect to which $\omega$ is self-dual; for this purpose, we
first build a smooth rank 3 subbundle $P$ of $\Lambda^2T^*X$,
positive definite with respect to the wedge-product, and
such that $\omega$ is a section of $P$. The smoothness assumption
implies that, at every point $x\in\Gamma=\omega^{-1}(0)$, $P_x$ must coincide
with the image of the intrinsically defined derivative $\nabla\omega_x\co T_xX
\to \Lambda^2 T^*X_x$. We can extend the construction of $P$ first to a
neighbourhood of $\Gamma$, and then to all of $X$, using the convexity
property mentioned in the previous paragraph to patch together local
constructions by means of a partition of unity.

By the discussion above,
there is a unique conformal class $[g]$ which realises the subbundle $P$
as the bundle of self-dual forms. For any
metric $g$ in this conformal class, the 2--form $\omega$ is self-dual, and
then closedness implies harmonicity. This completes the proof of the first
statement in the Proposition.

We now consider the claim that
if $X$ is compact and $b_{2}^{+}(X)\geq 1$ then for generic
Riemannian metrics on $X$ one can obtain near-symplectic structures
from self-dual harmonic forms.  This is proved by considering the space
$\mathcal{C}$ of pairs $(g,a)$, where $g$ is a $C^{k,\alpha}$
Riemannian metric on $X$ and $a\in H^2_{+,g}$ is a cohomology class
such that $a^2=1$ and admitting a self-dual representative. The universal
bundle $\Lambda^+$ over $X\times \mathcal{C}$, whose fibre at $(x,g,a)$ is
$\Lambda^2_{+,g} T^*X_x$, admits a universal section $\Omega$ whose
restriction to $X\times \{(g,a)\}$ is the unique harmonic self-dual 2--form in
the given cohomology class. It can be shown that $\Omega$ is transverse
to the zero section of $\Lambda^+$ (see for example \cite{LB}, Section 3).
The statement follows by observing that the regular values
of the projection of $\Omega^{-1}(0)$ to $\mathcal{C}$ form
a dense subset of the second Baire category in $\mathcal{C}$.
Detailed proofs have already appeared in the literature, and the reader is
referred to \cite{Ho1} (Theorem 1.1) or \cite{LB} (Proposition~1).

The only remaining statement to prove is that $[\omega]\in H^2(X,\bR)$ can
be chosen to be the reduction of a rational class. However, this follows
readily from the observation that the set of all $(g,a)\in \mathcal{C}$
for which the self-dual harmonic form in the class $a$ has transverse zeros
is an open subset of $\mathcal{C}$, and therefore necessarily contains
points such that $a$ is proportional to a rational cohomology class.

\section{Topological considerations and examples} \label{sec:topology}

\subsection{Monodromy}
Consider a near-symplectic 4--manifold $(X,\omega)$ with
$\omega^{-1}(0)=\Gamma$, and a singular Lefschetz pencil
$f\co X\setminus A\to S^2$ such that each component of $\Gamma$ maps
bijectively to the equator as in
Theorem \ref{thm:main}. Up to a small perturbation we can assume that
$f$ is injective on the set $B$ of isolated critical points, and that
$f(B)\cap f(\Gamma)=\emptyset$. After blowing up the base points,
we obtain a new manifold $\hat{X}$, and $f$ extends to a well-defined map
$\hat{f}\co \hat{X}\to S^2$.

Let $V$ be a tubular neighbourhood of the equator in $S^2$, disjoint from
$f(B)$, and denote by $D_\pm$ the two components of
$S^2\setminus V$. Then we can decompose $\hat{X}$ into three pieces:
$X_+=\hat{f}^{-1}(D_+)$, $W=\hat{f}^{-1}(V)$, and $X_-=\hat{f}^{-1}(D_-)$.
The zero locus $\Gamma$ of the near-symplectic form is entirely contained
in $W$. The manifolds $X_\pm$ are symplectic, and the restriction of
$\hat{f}$ to $X_\pm$ yields two symplectic Lefschetz fibrations
$f_\pm\co X_\pm\to D_\pm$, with fibres $\Sigma_{\pm}$. 

Consider the quadratic local model $(\ux,t)\mapsto (Q(\ux),t)$ describing the
behaviour of $f$ near $\Gamma$: the fibres are locally given by hyperboloids
in $\bR^3$, two-sheeted for $Q>0$ and one-sheeted for $Q<0$, with a conical
singularity for $Q=0$. Hence, the fibres for $Q<0$ are obtained from those
for $Q>0$ by attaching a handle, which decreases the Euler characteristic
by $2$. Since the diffeomorphisms used to paste this local model into $f$
are oriented in the same manner for all components of $\Gamma$, the induced
normal orientations of the equator are consistent, and we can choose $D_+$
(resp.\ $D_-$) to correspond to positive (resp.\ negative) values of $Q$ in
the local models near all components of $\Gamma$. With this convention,
$\chi(\Sigma_-)=\chi(\Sigma_+)-2m$, where $m$ is the number of components
of $\Gamma$ (if we assume that $\Sigma_+$ is connected of genus $g$, then
the genus of $\Sigma_-$ is $g+m$).

Since the restriction of $\hat{f}$ to $W$ has no critical points outside of
$\Gamma$, the 4--manifold $W$ is a fibre bundle over $S^1$, whose fibre $Y$
(the preimage of a small arc transverse to the equator) defines a cobordism
between $\Sigma_+$ and $\Sigma_-$, consisting of a series of handle
attachment operations (one for each component of $\Gamma$). Hence
$W$ relates the boundaries of $X_+$ and $X_-$ to each other via
a sequence of fibrewise handle additions.

More precisely, identify $\overline{V}$ with $S^1\times [-\delta,\delta]$,
and consider for each $\theta\in S^1$ the two boundary fibres
$\Sigma_{\pm,\theta}=\hat{f}^{-1}(\theta,\pm \delta)$. Then
$\Sigma_{-,\theta}$ is obtained from $\Sigma_{+,\theta}$ by deleting
$2m$ small discs and identifying $m$ pairs of boundary components.
Conversely $\Sigma_{+,\theta}$ is obtained from $\Sigma_{-,\theta}$ by
cutting it open along $m$ disjoint simple closed curves, and capping the
boundary components with discs. 

Letting $\theta$ vary, the
union of these discs forms the tubular neighbourhood $U_L$ of a {\it link}
$L\subset \partial X_+$. The link $L$ intersects each fibre of $\partial
X_+$ in $2m$ points (ie, it is in fact a braid with $2m$ strands in
$\partial X_+$); these points are naturally partitioned into $m$
pairs, according to the manner in which the boundary components of
$\Sigma_{+,\theta}\setminus (\Sigma_{+,\theta}\cap U_L)$ are glued to each
other in order to obtain $\Sigma_{-,\theta}$. Since each pair of points
canonically corresponds to a component of $\Gamma$, the components of $L$
are naturally labelled (``coloured'') by components of $\Gamma$ (or,
less canonically, by integers $1,\dots,m$). 

Moreover, $L$ also carries naturally a {\it relative framing}, which
keeps track of the manner in which the boundary components of $\partial
X_+\setminus U_L$ with the same colour are identified. More precisely,
the relative framing is the choice of a smooth involution $\rho\co \partial
U_L\to\partial U_L$, preserving the fibration structure above $S^1$,
the colouring and the orientation, but exchanging the two components with
the same colour in each fibre, up to isotopy. Given two relative framings
$\rho,\rho'$, for each of the $m$ colours the restrictions of $\rho$ and
$\rho'$ to the corresponding components of $\partial U_L$ differ by an
element of $\pi_1\mathrm{Diff}(S^1)\simeq \bZ$. Hence, the set of
relative framings is a $\bZ^m$--torsor.

The monodromy of $\partial X_+$, the $2m$--strand braid
$L\subset \partial X_+$, the colouring $c\co L\to \{1,\dots,m\}$ and the
relative framing $\rho$ determine completely the topology of the fibred
cobordism $W$.

Recall that the symplectic
Lefschetz fibrations $f_\pm\co X_\pm\to D_\pm$ are determined by their
monodromies, which take values in the relative mapping class groups
$\mathrm{Map}(\Sigma_\pm,A)$, ie, the set of isotopy classes of
orientation-preserving diffeomorphisms of $\Sigma_{\pm}$ which coincide
with identity over a small neighbourhood of the base locus $A$. If we
assume that $\Sigma_\pm$ are connected of genus $g_\pm$ and the number of
base points is $n$, then $\mathrm{Map}(\Sigma_\pm,A)$ is nothing but the
mapping class group $\mathrm{Map}_{g_\pm,n}$ of a genus $g_\pm$ 
surface with $n$ boundary components. The
monodromy around each isolated singular fibre is a positive Dehn twist
along a simple closed curve (the corresponding vanishing cycle), and the
product of these Dehn twists is equal to the monodromy $\psi_\pm$ of the
boundary fibration $\partial X_\pm\to S^1$.

The coloured braid $L$ and the relative framing $\rho$ determine
a lift of $\psi_+$ from $\mathrm{Map}(\Sigma_+,A)$
to $\mathrm{Map}(\Sigma_-,A)$,
which we denote by $\hat{\psi}_+$. More precisely, starting from the mapping
torus $\partial X_+$ of $\psi_+$, by deleting a tubular neighbourhood of the
braid $L$ one obtains a new fibre bundle over $S^1$, whose fibre has genus $g_+$ and
$2m$ boundary components (if $L$ is trivial, this lifts $\psi_+$ from
$\mathrm{Map}_{g_+,n}$ to $\mathrm{Map}_{g_+,n+2m}$).
The colouring and the relative framing then specify a manner in which the
$2m$ boundary components are glued to each other, to obtain a bundle over
$S^1$ with closed fibres of genus $g_++m=g_-$, and whose monodromy is by
definition $\hat{\psi}_+\in \mathrm{Map}(\Sigma_-,A)$. Because this
3-manifold coincides with $\partial X_-$ up to a change of orientation, 
$\hat{\psi}_+\cdot \psi_-$ belongs to the kernel of the natural morphism
$\mathrm{Map}(\Sigma_-,A)\to \mathrm{Map}(\Sigma_-)$. However, because each
exceptional section of $\hat{f}$ obtained by blowing up $A$ has a normal
bundle of degree $-1$, the product $\hat{\psi}_+\cdot\psi_-$ is not
$\mathrm{Id}$, but rather the boundary twist $\delta_A\in
\mathrm{Map}(\Sigma_-,A)$, ie, the product of the Dehn twists along small
loops encircling the various points of $A$.

If we assume that the identity components in $\mathrm{Diff}(\Sigma_\pm,A)$
are simply connected (eg, if $\Sigma_\pm$ both have genus at least $2$),
then the manner in which the boundaries of $X_\pm$ and $W$ are glued to each
other is determined uniquely up to isotopy. The above data (the
monodromies of $X_\pm$, and the coloured link $L$ with its relative framing)
then determine completely the topology of $f$. Otherwise, the possible
gluings of $\partial X_\pm$ to the boundary of $W$ are parametrised by
elements of $\pi_1\mathrm{Diff}(\Sigma_{\pm},A)$.
\medskip

{\bf Example}\qua To make the above discussion more concrete, we briefly
consider the case where $X_+$ has no singular fibres ($X_+\simeq
\Sigma_+\times D^2$) and $\Gamma$ is connected. Then $L$ intersects each
fibre of $\partial X_+$ in two points, and $\Sigma_-$ is obtained
from $\Sigma_+$ by cutting it open at these two points and attaching a
handle in the manner prescribed by the relative framing of $L$.
The core of this handle is a simple closed loop $\gamma\subset \Sigma_-$,
which can be thought of as the ``vanishing cycle'' associated to the equator.

The link $L$ is an arbitrary element of the braid group $B_2(\Sigma_+)$,
ie, the fundamental group of the complement of the diagonal in the second
symmetric product of $\Sigma_+$. Depending on whether the monodromy
preserves or exchanges the two points of $\Sigma_+\cap L$ (ie, whether
$L$ has one or two components), the $S^1$--bundle over $S^1$ formed by the
``vanishing cycle'' inside $\partial X_-$ can be either a torus or a
Klein bottle. These two cases correspond respectively to the two local
models $N_+$ and $N_-$  described in the
Introduction for the behaviour of $\omega$ in a neighbourhood of $\Gamma$.

We finish with a simple remark illustrating the importance of
the relative framing of $L$.
Even when the braid $L$ is trivial, the boundary of $X_-$ need not be
diffeomorphic to $S^1\times \Sigma_-$: in general, the monodromy of
$\partial X_-$ can be an arbitrary power of the Dehn twist along the
vanishing cycle $\gamma\subset \Sigma_-$.

\subsection{Examples}

\begin{figure}[ht!]
\centering
\setlength{\unitlength}{0.5cm}
\begin{picture}(10,5)(-4.5,-2.5)
\qbezier[60](2.5,0)(2.5,1)(1.75,1.75)
\qbezier[60](0,2.5)(1,2.5)(1.75,1.75)
\qbezier[60](-2.5,0)(-2.5,1)(-1.75,1.75)
\qbezier[60](0,2.5)(-1,2.5)(-1.75,1.75)
\qbezier[60](2.5,0)(2.5,-1)(1.75,-1.75)
\qbezier[60](0,-2.5)(1,-2.5)(1.75,-1.75)
\qbezier[60](-2.5,0)(-2.5,-1)(-1.75,-1.75)
\qbezier[60](0,-2.5)(-1,-2.5)(-1.75,-1.75)
\qbezier[110](-2.5,-0.03)(-2,-0.4)(0,-0.4)
\qbezier[110](2.5,-0.03)(2,-0.4)(0,-0.4)
\qbezier[20](2.5,0.03)(2,0.4)(0,0.4)
\qbezier[20](-2.5,0.03)(-2,0.4)(0,0.4)
\put(-4.5,0){\makebox(0,0)[lc]{$W$}}
\put(-4.5,-2){\makebox(0,0)[lc]{$X_-$}}
\put(-4.5,2){\makebox(0,0)[lc]{$X_+$}}
\qbezier[10](-4.5,1)(-3.8,1)(-3.1,1)
\qbezier[10](-4.5,-1)(-3.8,-1)(-3.1,-1)
\put(-3.2,2){\vector(0,1){0.5}}
\put(-3.2,2){\vector(0,-1){0.9}}
\put(-3.2,-2){\vector(0,1){0.9}}
\put(-3.2,-2){\vector(0,-1){0.5}}
\put(-3.2,0){\vector(0,1){0.9}}
\put(-3.2,0){\vector(0,-1){0.9}}
\put(3.5,2){\vector(-1,0){1}}
\put(3.5,-2){\vector(-1,0){1}}
\put(3.8,0){\vector(-1,0){1}}
\put(4.5,2){\circle{0.8}}
\qbezier[25](4.1,1.98)(4.5,1.8)(4.87,1.98)
\qbezier[10](4.1,2)(4.5,2.18)(4.87,2)
\put(4.7,-2){\circle{0.8}}
\put(4.7,-2){\circle{1.7}}
\qbezier[15](3.85,-2)(4.07,-2.15)(4.3,-2)
\qbezier[15](5.5,-2)(5.32,-2.15)(5.1,-2)
\qbezier[7](3.85,-2)(4.07,-1.85)(4.3,-2)
\qbezier[7](5.5,-2)(5.32,-1.85)(5.1,-2)
\qbezier[20](5.1,0)(5.1,0.3)(4.8,0.3)
\qbezier[20](4.8,0.3)(4.6,0.3)(4.4,0)
\qbezier[25](4.4,0)(4.2,0.3)(4.5,0.6)
\qbezier[30](4.5,0.6)(4.9,0.9)(5.3,0.6)
\qbezier[25](5.3,0.6)(5.6,0.3)(5.6,0)
\qbezier[20](5.1,0)(5.1,-0.3)(4.8,-0.3)
\qbezier[20](4.8,-0.3)(4.6,-0.3)(4.4,0)
\qbezier[25](4.4,0)(4.2,-0.3)(4.5,-0.6)
\qbezier[30](4.5,-0.6)(4.9,-0.9)(5.3,-0.6)
\qbezier[25](5.3,-0.6)(5.6,-0.3)(5.6,0)
\qbezier[15](5.1,0)(5.35,-0.15)(5.6,0)
\qbezier[7](5.1,0)(5.35,0.15)(5.6,0)
\end{picture}
\caption{A genus 0/1 singular fibration}\label{fig:example1}
\end{figure}
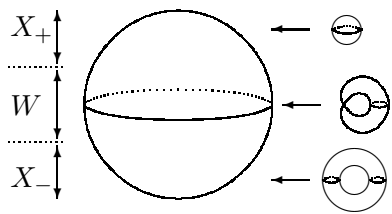

{\bf Example 1}\qua The simplest non-trivial examples of singular Lefschetz
fibrations $f\co X\to S^2$ are those where $\Gamma$ is connected, with a
neighbourhood modelled on $N_+$, there are
no isolated singular fibres, and the fibres are connected of genus $0$ over
$D_+$ and genus $1$ over $D_-$ (see Figure \ref{fig:example1}).

The total space of the fibration is a smooth 4--manifold $X$ obtained by gluing
together the three open pieces $X_-\simeq T^2\times D^2$ lying over the
southern hemisphere $D_-$, $W$ lying over
a neighbourhood of the equator, and $X_+\simeq S^2\times D^2$ lying over the
northern hemisphere $D_+$. The manifold $W$ is a product of $S^1$ with the
standard cobordism from the torus $T^2$ to sphere $S^2$, which is
diffeomorphic to a solid torus with a small ball removed. Hence, $W\simeq
S^1\times (S^1\times D^2\setminus B^3)$.

Because the diffeomorphism groups of $S^2$ and $T^2$ are not simply
connected, there are various possible choices for the identification
diffeomorphisms $\phi_\pm$ between the boundaries $S^1\times S^2$
(resp.\ $S^1\times T^2$) of $X_\pm$ and $W$.
Since $\phi_\pm$ must be compatible with the fibration structure over $S^1$,
they are described by families of diffeomorphisms of the boundary fibres,
ie, elements of $\pi_1\mathrm{Diff}(S^2)\simeq \bZ/2$ and
$\pi_1\mathrm{Diff}(T^2)\simeq \bZ^2$ (compare with the case of ordinary
sphere or torus bundles over $S^2$). 

Let us consider eg, the ``untwisted'' fibration $f\co X\to S^2$,
corresponding to trivial
choices for both gluings. This fibration admits a section with trivial
normal bundle (considering a point lying away from the ``vanishing cycle''
in each $T^2$ fibre, and the corresponding point in each $S^2$ fibre), and
its fundamental group is $\bZ$ (generated by a loop transverse to
the vanishing cycle in a $T^2$ fibre). Its total space is diffeomorphic
to the connected sum $(S^1\times S^3) \# (S^2\times S^2)$. Indeed, using the
decomposition of $S^3$ into two solid tori, it is easy to see that $X_-\cup W$
is diffeomorphic to the complement of an embedded loop $\gamma$ in
$S^1\times S^3$ (the $S^1$ factor corresponds to the direction transverse
to the vanishing cycles in the $T^2$ fibres). In the untwisted case, the
loop $\gamma$ projects to a single point in the $S^1$ factor, and represents
an unknot in $S^3$; in particular, it can be contracted into an arbitrarily
small ball in $S^1\times S^3$, and the attachment of the handle $X_+$ can be
viewed as a connected sum operation performed on $S^1\times S^3$. Observing
that $S^4$ splits into $(S^1\times B^3)\cup (D^2\times S^2)$ and hence that
the corresponding handle attachment operation turns $S^4\setminus S^1$ into a
$S^2$--bundle over $S^2$ (in this case $S^2\times S^2$), we conclude that $X$
is as claimed.

If we still glue $X_-$ via the trivial element in $\pi_1\mathrm{Diff}(T^2)$
but glue $X_+$ using the non-trivial element in $\pi_1\mathrm{Diff}(S^2)$,
then we obtain $(S^1\times S^3)\# \mathbb{CP}^2\#
\overline{\mathbb{CP}}{}^2$ instead. However, if eg, we twist the
fibration by a loop of diffeomorphisms of $T^2$
corresponding to a unit translation in the direction transverse to the
vanishing cycle, we lose the existence of a section, and the total space
becomes simply connected. In fact, the new total space $X'$ is
diffeomorphic to $S^4$. Indeed, $X_-\cup W$ is still the complement
of a closed loop in $S^1\times S^3$, but the missing loop
$\gamma'$ now projects non-trivially to the $S^1$ factor, and is
isotopic to $S^1\times \{pt\}\subset S^1\times S^3$. Therefore, we now have
$X_-\cup W\simeq S^1\times B^3$, and by gluing $X_+=D^2\times S^2$ along the
boundary we obtain $X'\simeq S^4$. Theorem \ref{thm:converse} fails to apply
in this case, because the cohomological assumption fails to hold (the
fibres are homologically trivial).
\medskip

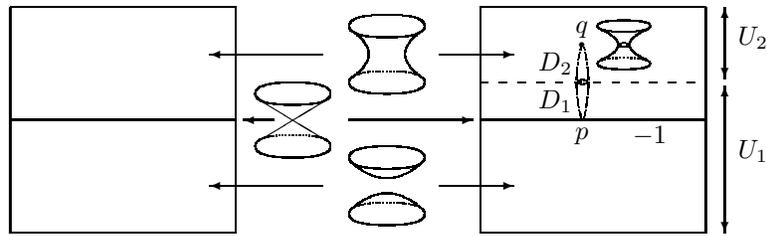
\begin{figure}[ht!]
\centering
\setlength{\unitlength}{5mm}
\begin{picture}(19,6)(-1,0)
\put(-1,0){\line(1,0){6}}
\put(-1,6){\line(1,0){6}}
\put(-1,0){\line(0,1){6}}
\put(5,0){\line(0,1){6}}
\put(-1,3){\line(1,0){6}}
\put(7.3,1.25){\vector(-1,0){3}}
\qbezier[55](8,2)(8,2.3)(9,2.3)
\qbezier[55](10,2)(10,2.3)(9,2.3)
\qbezier[55](8,2)(8,1.7)(9,1.7)
\qbezier[55](10,2)(10,1.7)(9,1.7)
\qbezier[80](8.03,1.9)(9,1)(9.97,1.9)
\qbezier[20](8,0.5)(8,0.8)(9,0.8)
\qbezier[20](10,0.5)(10,0.8)(9,0.8)
\qbezier[55](8,0.5)(8,0.2)(9,0.2)
\qbezier[55](10,0.5)(10,0.2)(9,0.2)
\qbezier[80](8.03,0.6)(9,1.5)(9.97,0.6)
\put(7.3,4.75){\vector(-1,0){3}}
\qbezier[55](8,5.5)(8,5.8)(9,5.8)
\qbezier[55](10,5.5)(10,5.8)(9,5.8)
\qbezier[55](8,5.5)(8,5.2)(9,5.2)
\qbezier[55](10,5.5)(10,5.2)(9,5.2)
\qbezier[20](8,4)(8,4.3)(9,4.3)
\qbezier[20](10,4)(10,4.3)(9,4.3)
\qbezier[55](8,4)(8,3.7)(9,3.7)
\qbezier[55](10,4)(10,3.7)(9,3.7)
\qbezier[80](8.03,5.4)(9,4.75)(8.03,4.1)
\qbezier[80](9.97,5.4)(9,4.75)(9.97,4.1)
\qbezier[55](5.5,3.7)(5.5,4)(6.5,4)
\qbezier[55](7.5,3.7)(7.5,4)(6.5,4)
\qbezier[55](5.5,3.7)(5.5,3.4)(6.5,3.4)
\qbezier[55](7.5,3.7)(7.5,3.4)(6.5,3.4)
\qbezier[20](5.5,2.3)(5.5,2.6)(6.5,2.6)
\qbezier[20](7.5,2.3)(7.5,2.6)(6.5,2.6)
\qbezier[55](7.5,2.3)(7.5,2)(6.5,2)
\qbezier[55](5.5,2.3)(5.5,2)(6.5,2)
\put(5.6,2.45){\line(5,3){1.8}}
\put(5.6,3.55){\line(5,-3){1.8}}
\put(6,3){\vector(-1,0){0.8}}
\put(11.5,0){\line(1,0){6}}
\put(11.5,6){\line(1,0){6}}
\put(11.5,0){\line(0,1){6}}
\put(17.5,0){\line(0,1){6}}
\put(11.5,3){\line(1,0){6}}
\put(10.4,1.25){\vector(1,0){2}}
\put(10.4,4.75){\vector(1,0){2}}
\put(8,3){\vector(1,0){3.3}}
\put(18,2){\vector(0,-1){2}}
\put(18,2){\vector(0,1){1.9}}
\put(18,5){\vector(0,1){1}}
\put(18,5){\vector(0,-1){0.9}}
\multiput(11.5,4)(0.5,0){12}{\line(1,0){0.2}}
\put(14.2,5){\circle*{0.15}}
\qbezier[35](14.6,5.5)(14.6,5.65)(15.3,5.65)
\qbezier[35](16.0,5.5)(16.0,5.65)(15.3,5.65)
\qbezier[35](14.6,5.5)(14.6,5.35)(15.3,5.35)
\qbezier[35](16.0,5.5)(16.0,5.35)(15.3,5.35)
\qbezier[13](14.6,4.5)(14.6,4.65)(15.3,4.65)
\qbezier[13](16.0,4.5)(16.0,4.65)(15.3,4.65)
\qbezier[35](14.6,4.5)(14.6,4.35)(15.3,4.35)
\qbezier[35](16.0,4.5)(16.0,4.35)(15.3,4.35)
\qbezier[50](14.7,5.4)(15.6,5)(14.7,4.6)
\qbezier[50](15.9,5.4)(15.0,5)(15.9,4.6)
\qbezier[10](15.15,5)(15.3,4.9)(15.45,5)
\qbezier[5](15.15,5)(15.3,5.1)(15.45,5)
\put(16,2.85){\makebox(0,0)[ct]{\small $-1$}}
\put(18.3,2){\small $U_1$}
\put(18.3,5){\small $U_2$}
\qbezier[30](14.22,5)(14.5,4)(14.22,3)
\qbezier[30](14.18,5)(13.9,4)(14.18,3)
\qbezier[6](14.07,4)(14.2,3.9)(14.33,4)
\qbezier[6](14.07,4.02)(14.2,4.12)(14.33,4.02)
\put(14,4.5){\makebox(0,0)[rc]{\small $D_2$}}
\put(14,3.5){\makebox(0,0)[rc]{\small $D_1$}}
\put(14.2,2.8){\makebox(0,0)[ct]{\small $p$}}
\put(14.2,5.2){\makebox(0,0)[cb]{\small $q$}}
\end{picture}
\caption{Blowing up near $\Gamma$: $f$ (left) and $f'$ (right)}
\label{fig:example2}
\end{figure}

{\bf Example 2 -- Isotropic blow-up}\qua
There are several different operations that can be performed on a singular
Lefschetz fibration $f\co X^4\to S^2$ in order to modify its total space
by a topological blow-up operation (ie, connected sum with
$\overline{\mathbb{CP}}{}^2$). Keeping symplectic
Lefschetz fibrations in mind, the ``usual'' blow-up construction amounts to
the insertion of an isolated singular fibre with a homotopically trivial
vanishing cycle. The exceptional sphere is then obtained as a component of
the singular fibre, and is hence naturally symplectic with respect to any
2--form compatible with the fibration structure. If we perform the blow-up
near a point $p\in\Gamma$, we can instead modify $f$ according to
the local operation represented on Figure \ref{fig:example2}.

We start from a small ball $B^4$ centred at $p$, over which $f$ is
as shown in the left half of Figure \ref{fig:example2}, and
replace it with the total space of the fibration $f'$ represented in the
right half of the figure.
The map $f'$ differs from $f$ in two respects: (1) it has an additional
isolated critical point $q\in X_-$, where the vanishing cycle $\gamma$ is the
same as at $p$; (2) the relative framing of the link $L\subset \partial X_+$
is modified by $-1$. As explained at the end of the previous Section,
changing the relative framing modifies the lift $\hat{\psi}_+$ of the
monodromy of $\partial X_+$ to $\mathrm{Map}(\Sigma_-,A)$ by the inverse
of the Dehn twist along $\gamma$; this compensates the modification of
the monodromy of $\partial X_-$ by the same Dehn twist due to the new
isolated singular fibre.

The total space of $f'$ is the union of two subsets $U_1$ and $U_2$ (see
figure), both diffeomorphic to 4--balls. Over $U_1$ the map $f'$ is
modelled on $(t,x,y,z)\mapsto (t,x^2+y^2-z^2)$, while over $U_2$ it is
modelled on $(z_1,z_2)\mapsto z_1^2+z_2^2$. The total space of $U_1$ can
be viewed as a disc bundle over a disc $D_1=\{z=t=0\}$, while the
total space of $U_2$ is a disc bundle over a disc $D_2=\{\mathrm{Im}\,z_1=
\mathrm{Im}\,z_2=0\}$. The boundaries of the two discs $D_1$ and $D_2$ match
with each other, so that the total space of $f'$ is a disc bundle over a
sphere $S=D_1\cup D_2$ (dotted in Figure~\ref{fig:example2}). Moreover, it
is easy to check that the normal bundle of $S$ has degree $-1$.

From the near-symplectic point of view, 
this type of blow-up is not equivalent to the usual one.
Indeed, in this setup the exceptional sphere $S$ arises from a matching
pair of vanishing cycles above an arc joining the critical values $f'(p)$
and $f'(q)$, and for a suitable choice of the compatible near-symplectic
form $\omega$ on the total space of $f'$ it will be $\omega$--isotropic
(or ``near-Lagrangian'').

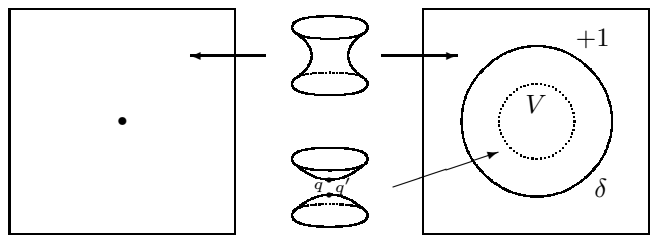
\begin{figure}[ht!]
\centering
\setlength{\unitlength}{5mm}
\begin{picture}(19,6)(-1,0)
\put(0.5,0){\line(1,0){6}}
\put(0.5,6){\line(1,0){6}}
\put(0.5,0){\line(0,1){6}}
\put(6.5,0){\line(0,1){6}}
\qbezier[50](16.5,3)(16.5,3.8)(15.9,4.4)
\qbezier[50](14.5,5)(15.3,5)(15.9,4.4)
\qbezier[50](12.5,3)(12.5,3.8)(13.1,4.4)
\qbezier[50](14.5,5)(13.7,5)(13.1,4.4)
\qbezier[50](16.5,3)(16.5,2.2)(15.9,1.6)
\qbezier[50](14.5,1)(15.3,1)(15.9,1.6)
\qbezier[50](12.5,3)(12.5,2.2)(13.1,1.6)
\qbezier[50](14.5,1)(13.7,1)(13.1,1.6)
\put(10.7,1.25){\vector(3,1){2.8}}
\qbezier[55](8,2)(8,2.3)(9,2.3)
\qbezier[55](10,2)(10,2.3)(9,2.3)
\qbezier[55](8,2)(8,1.7)(9,1.7)
\qbezier[55](10,2)(10,1.7)(9,1.7)
\qbezier[80](8.03,1.9)(9,1)(9.97,1.9)
\qbezier[20](8,0.5)(8,0.8)(9,0.8)
\qbezier[20](10,0.5)(10,0.8)(9,0.8)
\qbezier[55](8,0.5)(8,0.2)(9,0.2)
\qbezier[55](10,0.5)(10,0.2)(9,0.2)
\qbezier[80](8.03,0.6)(9,1.5)(9.97,0.6)
\put(9,1.05){\circle*{0.15}}
\put(9.1,1.05){\makebox(0,0)[lb]{\tiny $q'$}}
\put(9,1.45){\circle*{0.15}}
\put(8.9,1.4){\makebox(0,0)[rt]{\tiny $q$}}
\put(7.3,4.75){\vector(-1,0){2}}
\qbezier[55](8,5.5)(8,5.8)(9,5.8)
\qbezier[55](10,5.5)(10,5.8)(9,5.8)
\qbezier[55](8,5.5)(8,5.2)(9,5.2)
\qbezier[55](10,5.5)(10,5.2)(9,5.2)
\qbezier[20](8,4)(8,4.3)(9,4.3)
\qbezier[20](10,4)(10,4.3)(9,4.3)
\qbezier[55](8,4)(8,3.7)(9,3.7)
\qbezier[55](10,4)(10,3.7)(9,3.7)
\qbezier[80](8.03,5.4)(9,4.75)(8.03,4.1)
\qbezier[80](9.97,5.4)(9,4.75)(9.97,4.1)
\put(15.5,5){\small $+1$}
\put(11.5,0){\line(1,0){6}}
\put(11.5,6){\line(1,0){6}}
\put(11.5,0){\line(0,1){6}}
\put(17.5,0){\line(0,1){6}}
\put(10.4,4.75){\vector(1,0){2}}
\put(3.5,3){\circle*{0.2}}
\qbezier[7](15.5,3)(15.5,3.4)(15.2,3.7)
\qbezier[7](14.5,4)(14.9,4)(15.2,3.7)
\qbezier[7](13.5,3)(13.5,3.4)(13.8,3.7)
\qbezier[7](14.5,4)(14.1,4)(13.8,3.7)
\qbezier[7](15.5,3)(15.5,2.6)(15.2,2.3)
\qbezier[7](14.5,2)(14.9,2)(15.2,2.3)
\qbezier[7](13.5,3)(13.5,2.6)(13.8,2.3)
\qbezier[7](14.5,2)(14.1,2)(13.8,2.3)
\put(14.5,3.7){\makebox(0,0)[ct]{\small $V$}}
\put(16,1){\small $\delta$}
\end{picture}
\caption{Inserting a critical circle: $f$ (left) and $f'$ (right)}
\label{fig:example3}
\end{figure}

\medskip
{\bf Example 3}\qua Consider an isolated Lefschetz-type critical point of a
singular fibration, with vanishing cycle a loop $\gamma$ in the nearby
generic fibre. We can remove a neighbourhood of this singular fibre and
insert in its place a configuration where the critical values form a
simple closed loop $\delta$, with fibre genus decreased by 
$1$ inside $\delta$, and using the same loop $\gamma$ as ``vanishing cycle'',
as shown in Figure \ref{fig:example3}. This adds a new component to
$\Gamma$ 
(this component is not mapped to the equator of $S^2$;
here we consider singular Lefschetz fibrations
more general than those given by Theorem \ref{thm:main}).

The fibres outside
$\delta$ are obtained from those inside by attaching a handle joining
two points $q,q'$ as shown in the figure. Along $\delta$ the points $q,q'$
describe a trivial braid, but the relative framing differs from the trivial
one by $+1$, so that on the outer side the monodromy around
$\delta$ consists of a single positive Dehn twist along $\gamma$ (which
balances the loss of the isolated singular fibre).

The total space of the local model for $f$ given on Figure
\ref{fig:example3} (left) is
simply a $4$--ball. On the other hand, the total space of the new fibration
$f'$ contains a smoothly embedded sphere $S$, obtained by
considering the two points $q$ and $q'$ in each of the fibres inside
$\delta$ (yielding the two hemispheres of $S$), and the singular points
in the fibres above $\delta$ (yielding the equator). Using the fact that
the monodromy around $\delta$ is a positive Dehn twist along $\gamma$, it can be
checked easily that $S$ has self-intersection $+1$. Moreover, the preimage
of the interior region $V$ is the disjoint union of two $D^2\times D^2$'s,
and hence a disc bundle over $S\cap f'{}^{-1}(V)$. On the other hand, the
preimage of the outer region is diffeomorphic to $S^1\times B^3$, and is
again a disc bundle over a neighbourhood of the equator in $S$. Therefore,
the total space of $f'$ is a disc bundle over the sphere $S$, and it is
diffeomorphic to the complement of a ball in $\mathbb{CP}^2$.

It follows that the operation we have described amounts to a connected sum
with $\mathbb{CP}^2$ -- an operation
whose result is never a symplectic 4--manifold unless the original manifold
had $b_2^+=0$, by the work of Taubes.
In particular, if the configuration $f'$ occurs inside a singular Lefschetz
fibration satisfying the assumptions of Theorem \ref{thm:converse},
then its total space has $b_2^+\ge 2$ and splits off a $\mathbb{CP}^2$
summand, and hence does not admit any symplectic structure (more generally,
this also holds for similar configurations with arbitrarily positive
relative framings, since these contain $+n$--spheres which can be blown up
to produce a $\mathbb{CP}^2$ summand).


\vspace{-2mm}
\begin{thebibliography}{99}

\itemsep 1ptplus1pt

\bibitem{Au1}
{\bf D Auroux},
{\it Symplectic 4--manifolds as branched coverings of $\mathbb{CP}^2$},
Invent. Math. 139 (2000) 551--602
  \MR{1738061}
\bibitem{Au2}
{\bf D Auroux},
{\it A remark about Donaldson's construction of symplectic submanifolds},
J. Symplectic Geom. {1} (2002) 647--658
  \MR{1959060}
\bibitem{Do1} 
{\bf S\,K Donaldson}, {\it Symplectic submanifolds and almost-complex
geometry}, J.\ Differential Geom. {44} (1996) 666--705
  \MR{1438190}
\bibitem{Do2}
{\bf S\,K Donaldson},  {\it Lefschetz fibrations in symplectic geometry},
Documenta Math. Extra Volume ICM II, (1998) 309--314
  \MR{1648081}
\bibitem{Do3} 
{\bf S\,K Donaldson}, {\it Lefschetz pencils on symplectic manifolds},
J. Differential Geom. {53} (1999) 205--236
  \MR{1802722}
\bibitem{DS}
{\bf S\,K Donaldson}, {\bf I Smith}, {\it Lefschetz pencils and the canonical class
for symplectic 4--manifolds}, Topology {42} (2003) 743--785
  \MR{1958528}
\bibitem{GK}
{\bf D\,T Gay}, {\bf R Kirby}, {\it Constructing symplectic forms on 4--manifolds
which vanish on circles}, \gtref8{2004}{20}{743}{777}
  \MR{2057780}
\bibitem{Go2}
{\bf R\,E Gompf}, {\it Toward a topological characterization of 
symplectic manifolds},
J.\ Symplectic Geom. {2} (2004) 177--206
  \MR{2108373}
\bibitem{Ho1}
{\bf K Honda}, {\it Transversality theorems for harmonic forms},
Rocky Mountain J.\ Math. {34} (2004) 629--664
  \MR{2072799}
\bibitem{Ho2}
{\bf K Honda}, {\it Local properties of self-dual harmonic 2--forms on a
4--manifold}, J.\ Reine Angew. Math.\ {577} (2004) 105--116
  \MR{2108214}
\bibitem{LB}
{\bf C Le\,Brun}, {\it Yamabe constants and the perturbed Seiberg--Witten
equations}, Comm.\ Anal.\ Geom.\ {5} (1997) 535--553
  \MR{1487727}
\bibitem{Presas}
{\bf F Presas}, 
{\it Submanifolds of symplectic manifolds with contact border},
\arxiv{math.SG/0007037}
\bibitem{Ta1}
{\bf C\,H Taubes}, {\it The geometry of the Seiberg--Witten invariants},
Surveys in Differential Geometry, Vol.\ III (Cambridge, 1996)
Int.\ Press, Boston (1998) 299--339
  \MR{1677891}
\bibitem{Ta2}
{\bf C\,H Taubes},  {\it Seiberg--Witten invariants and pseudo-holomorphic
subvarieties for self-dual, harmonic 2--forms}, \gtref3{1999}{8}{167}{210}
  \MR{1697181}
\bibitem{Th}
{\bf W Thurston}, {\it Some simple examples of symplectic manifolds},
Proc.\ Amer.\ Math.\ Soc.\ {55} (1976) 467--468
  \MR{0402764}

\end{thebibliography}
\end{document}